\documentclass[10pt]{amsart}
\usepackage{amssymb, latexsym}
\usepackage{amsmath}
\usepackage{color}
\usepackage{mathrsfs}
\usepackage{url}
\newtheorem{Lemma}{Lemma}[section]
\newtheorem{theorem}[Lemma]{Theorem}
\newtheorem{proposition}[Lemma]{Proposition}
\newtheorem{corollary}[Lemma]{Corollary}
\newtheorem{definition}[Lemma]{Definition}
\newtheorem{remark}[Lemma]{Remark}
\newtheorem{conjecture}[Lemma]{Conjecture}
\newtheorem{lemma}[Lemma]{Lemma}
\newtheorem{example}[Lemma]{Example}

\newcommand{\prf}[1]{\noindent {\sc Proof.} #1 \hspace*{\fill} $\Box$}

\title{Minimal slope conjecture of $F$-isocrystals}

\author{Nobuo Tsuzuki}
\date{\today}


\begin{document}

\footnote[0]{Mathematical Institute, Graduate School of Science, 
Tohoku University, Aramaki Aza-Aoba 6-3, Aobaku, Sendai, 980-8578, Japan}
\footnote[0]{\textit{E-mail address.} \url{tsuzuki@math.tohoku.ac.jp}}

\begin{abstract}
The minimal slope conjecture, which was proposed by K.S.Kedlaya, asserts that two irreducible overconvergent $F$-isocrystals 
on a smooth variety are isomorphic to each other if both minimal slope constitutions of slope filtrations are isomorphic to each other. 
We affirmatively solve the minimal slope conjecture for overconvergent $F$-isocrystals on curves 
and for overconvergent $\overline{\mathbb Q}_p$-$F$-isocrystals on smooth varieties over finite fields. 
\end{abstract}

\maketitle

\vspace*{-5mm}

\section{Introduction}

In this paper we study the minimal slope conjecture of overconvergent $F$-isocrystals 
which was proposed by K.S.Kedlaya in \cite{Ke16}. At first let us explain the conjecture and 
our main results. 

\subsection{Slope filtrations}
Let us fix the notation as follows:
\begin{itemize}
\item $p$ : a prime number;
\item $k$ : a perfect field of characteristic $p$;
\item $R$ : a complete discrete valuation field of 
mixed characteristic with residue field $k = R/\mathbf m$;
\item $K$ : the field of fractions of $R$;
\item $a \mapsto |a|$ : a multiplicative valuation of $K$; 
\item $\sigma : K \rightarrow K$ : a $q$-Frobenius on $K$ for a positive power $q$ of $p$, that is, 
a $p$-adic continuous homomorphism of fields 
such that $\sigma(a)\, \equiv\, a^q\, (\mathrm{mod}\, \mathbf m)$ for any $a \in R$. 
\end{itemize}

Let $X$ be a scheme separated locally of finite type over $\mathrm{Spec}\, k$. 
Let $\mathcal M^\dagger$ be an overconvergent $F$-isocrystal on $X/K$, 
and denote the convergent $F$-isocrystal on $X/K$ associated to $\mathcal M^\dagger$ by $\mathcal M$. 
We say $\mathcal M$ admits a slope filtration if 
there exists an increasing filtration 
$$
      0 = \mathcal M_0 \subsetneq \mathcal M_1 \subsetneq \cdots \subsetneq 
      \mathcal M_{r-1} \subsetneq\mathcal M_r = \mathcal M
$$
of $\mathcal M$ as convergent $F$-isocrystals on $X/K$ such that 
\begin{list}{}{}
\item[(i)] $\mathcal M_i/\mathcal M_{i-1}$ 
is nonzero and isoclinic of slope $s_i$ and 
\item[(ii)] $s_1 < s_2 < \cdots < s_r$. 
\end{list}
We call $s_1$ (resp. $s_r$) the minimal slope (resp. the maximal slope) of $\mathcal M$ when $\mathcal M \ne 0$. 
If the Newton polygons of the Frobenius structure of a convergent $F$-isocrystal $\mathcal M$ 
are constant on a smooth scheme $X$, 
$\mathcal M$ admits a slope filtration (\cite[Corollary 4.2]{Ke16}, \cite[Corollary 2.6]{Ts19}). It is known that, for any convergent $F$-isocrystal $\mathcal M$ 
on a smooth connected scheme $X$, there exists 
an open dense subscheme $U$ of $X$ such that 
the Newton polygons of $\mathcal M$ are constant by Grothendieck's specialization theorem (see \cite[Theorem 3.1.2]{Ke16} 
and \cite[Proposition 2.2]{Ts19}). Hence $\mathcal M$ always admits a unique slope filtration as convergent $F$-isocrystals after a certain shrinking of $X$.

\subsection{The minimal slope conjecture}

Kedlaya proposed the following problem in \cite[Remarks 5.14, 5.15]{Ke16}, 
which we call the minimal slope conjecture. 

\begin{conjecture}{\mbox{\rm \cite[Remark 5.14]{Ke16}}}\label{Kcj}
Let $X$ be a smooth connected scheme separated of finite type over $\mathrm{Spec}\, k$. 
Let $\mathcal M^\dagger$ and $\mathcal N^\dagger$ be irreducible overconvergent $F$-isocrystals on $X/K$ 
such that both convergent $F$-isocrystals $\mathcal M$ and $\mathcal N$ 
associated to $\mathcal M^\dagger$ and $\mathcal N^\dagger$ 
admit the slope filtrations $\{\mathcal M_i\}$ and $\{ \mathcal N_j\}$, respectively. 
Suppose that there is an isomorphism $h : \mathcal M_1 \rightarrow \mathcal N_1$ 
between the minimal slope constitution of slope filtrations of $\mathcal M$ and $\mathcal N$ as convergent $F$-isocrystals. Then 
there exists a unique isomorphism $g^\dagger : \mathcal M^\dagger \rightarrow \mathcal N^\dagger$ 
of overconvergent $F$-isocrystals such that the induced diagram 
$$
     \begin{array}{ccc}
              \mathcal M_1 &\overset{h}{\rightarrow} &\mathcal N_1 \\
             \cap & &\cap\\
         \mathcal M &\underset{g}{\rightarrow} &\mathcal N \\
   \end{array}
$$
is commutative in the category of convergent $F$-isocrystals on $X/K$. 
\end{conjecture}

If $X$ is proper, then $\mathcal M_1 = \mathcal M = \mathcal M^\dagger$ by the irreducibility. Hence 
the conjecture is trivially true. 
E.Ambrosi and M.D'Addezio proved the conjecture over a finite field $k$ only with the hypothesis 
of the nontriviality of the morphism $h : \mathcal M_1 \rightarrow \mathcal N_1$ 
when $\mathcal N^\dagger$ is of rank one \cite[Theorem 1.1.1]{AD18}. 
The version with the relaxed 
hypothesis on $h$ is called the stronger version of Kedlaya's conjecture.  
They applied the result to a generalized Lang-N\'eron's theorem on a finiteness of torsion points of Abelian varieties \cite[Theorem 1.2.1]{AD18}. 
    
In this paper we will study the dual form of the minimal slope conjecture as follows. 

\begin{conjecture}{\mbox{\rm (Dual form of the stronger version of Conjecture \ref{Kcj})}}\label{Kcjdual}
With the notation in Conjecture \ref{Kcj} and renumbering the slope filtration such as 
$$
      \mathcal M = \mathcal M^0 \supsetneq \mathcal M^1 \supsetneq \mathcal M^2 \supsetneq 
      \cdots  \supsetneq \mathcal M^{r-1}  \supsetneq \mathcal M^r = 0
$$
with the sequence of slopes $s^0 > s^1 > \cdots > s^{r-1}$, suppose that there is a nontrivial morphism 
$h : \mathcal N/\mathcal N^1 \rightarrow \mathcal M/\mathcal M^1$ between the maximal 
slope quotients as convergent $F$-isocrystals. 
Then there exists a unique isomorphism $g^\dagger : \mathcal N^\dagger \rightarrow \mathcal M^\dagger$ 
of overconvergent $F$-isocrystals such that the induced diagram 
$$
     \begin{array}{ccc}
              \mathcal N &\overset{g}{\rightarrow} &\mathcal M \\
             \downarrow & &\downarrow\\
         \mathcal N/\mathcal N^1 &\underset{h}{\rightarrow} &\mathcal M/\mathcal M^1 \\
   \end{array}
$$
is commutative in the category of convergent $F$-isocrystals on $X/K$. 
\end{conjecture}

\subsection{Results and strategies.} In this paper we establish affirmative results for the dual form of the 
minimal slope conjecture. 

\begin{theorem} \mbox{\rm (Corollary \ref{kedcur})} 
If $C$ is a smooth connected curve over $\mathrm{Spec}\, k$, then Conjecture \ref{Kcjdual} holds. 
\end{theorem}

Our method is different from Ambrosi and D'Addezio's monodromy group method in \cite{AD18}. 
We may assume $C$ is affine and will compare the modules of global sections of $\mathcal M^\dagger$ and $\mathcal N^\dagger$ inside 
the modules of global sections of the maximal slope quotient $\mathcal M/\mathcal M^1$. 
The key ingredients are the notion of PBQ (pure of bounded quotient) for overconvergent $F$-isocrystals in both local and global theories, 
the notion of saturated overconvergent $F$-isocrystals in both local and global theories, 
and the opposite filtrations of $\varphi$-modules in the local theory. 

The notion of PBQ was introduced by B.Chiarellotto and the author to give a necessary condition 
of Dwork's conjecture on the comparison between Frobenius slopes and logarithmic growth for $\varphi$-$\nabla$-modules in \cite{CT11} 
(see Remark \ref{logfrob}). In this paper we globalize the notion of PBQ which requires 
the space of bounded solutions on the generic disc of the associated Frobenius-differential module at the generic point is pure of Frobenius slope. 
We introduce the notion of saturation for overconvergent $F$-isocrystals, which requires 
the module of global sections of $\mathcal M^\dagger$ is naturally included in that 
of the maximal slope quotient $\mathcal M/\mathcal M^1$. 
The notion of opposite filtrations of $\varphi$-modules were introduced by J.A.De Jong to study the homomorphisms 
of $p$-divisible groups on local rings of equal characteristic $p$ in \cite{dJ98}. 
Applying the local theory, we prove the rank of the maximal slope quotient of any overconvergent $F$-isocrystal which is included in 
$\mathcal M/\mathcal M^1$ is less than or equal to the rank of $\mathcal M/\mathcal M^1$. 
Then an overconvergent $F$-isocrystal on a curve is irreducible if and only if 
it is PBQ and saturated and the maximal slope quotient $\mathcal M/\mathcal M^1$ is irreducible (Proposition \ref{globalirr}). 
So the given nontrivial morphism $h : \mathcal N/\mathcal N^1 \rightarrow \mathcal M/\mathcal M^1$ is isomorphic and 
we can compare, in the level of global sections, 
the irreducible overconvergent $F$-isocrystals $\mathcal M^\dagger$ and $\mathcal N^\dagger$ in 
the maximal slope quotient $\mathcal M/\mathcal M^1$. Then 
we obtain that $\mathcal M^\dagger$ and $\mathcal N^\dagger$ coincide 
with each other by the PBQ property and the upper bound of ranks of maximal slope quotients. 

Let $\overline{\mathbb Q}_p$ be an algebraic closure of the field $\mathbb Q_p$ of $p$-adic numbers. 
In the case of general dimensions we deal with the case where $k$ is a finite field and study the minimal slope problem for 
overconvergent $\overline{\mathbb Q}_p$-$F$-isocrystals which T.Abe introduced 
to established the celebrated work on the $p$-adic Langlands correspondence and the companion theorem in \cite{Ab18b} \cite{Ab18}. 
Our result for the general dimensions is as follows. 

\begin{theorem}\label{aa} \mbox{\rm (Theorem \ref{rlt2})} 
Let $X$ be a smooth connected scheme separated of finite type over the spectrum $\mathrm{Spec}\, k$ of a finite field $k$, 
and $\mathcal M^\dagger$ and $\mathcal N^\dagger$ irreducible overconvergent $\overline{\mathbb Q}_p$-$F$-isocrystals on $X$ 
which admit slope filtrations as convergent $\overline{\mathbb Q}_p$-$F$-isocrystals. 
If $h : \mathcal N/\mathcal N^1 \rightarrow \mathcal M/\mathcal M^1$ is a nontrivial morphism 
between the maximal slope quotients, then 
there exists an isomorphism $g^\dagger : \mathcal N^\dagger \rightarrow \mathcal M^\dagger$ of 
overconvergent $\overline{\mathbb Q}_p$-$F$-isocrystals. 
\end{theorem}

Abe and H.Esnault proved Lefschetz theorem for $\overline{\mathbb Q}_p$-$F$-isocrystals which asserts an existence of a smooth curve $C_\alpha$ 
passing at any given closed point $\alpha$ in an open dense subscheme of $X$ such that the restriction of the given irreducible $\overline{\mathbb Q}_p$-$F$-isocrystal on $X$ 
is again irreducible on $C_\alpha$ in \cite{AE16}. They also applied Lefschetz theorem to the weight theory after Abe-D.Caro's work in \cite{AC18}. 
Then the coincidence of characteristic polynomials of Frobenius of $\mathcal M^\dagger$ and $\mathcal N^\dagger$ holds 
at each closed point $\alpha$ by our study of the minimal slope conjecture on curves. 
Then applying \v{C}ebatarev's density theorem \cite{Ab18b} we obtain an isomorphism $g^\dagger : \mathcal N^\dagger \rightarrow \mathcal M^\dagger$ 
of overconvergent $\overline{\mathbb Q}_p$-$F$-isocrystals. 
However we can not show the compatibility between $g^\dagger$ and the given morphism $h$ at this moment (see Lemma \ref{suffcond}). 

\begin{remark}
D'Addezio solved Crew's parabolicity conjecture \cite[p.460]{Cr92} and the minimal slope conjecture for arbitrary base field $k$, recently \cite{DA20}. 
\end{remark}

\subsection{Further problems} The minimal slope conjecture seems to be not an $\ell$-adic but a $p$-adic own problem. 
However, the author also expects to find some consequences in $\ell$-adic theory from the $p$-adic theory 
with a view point of the companion theorem. 
In the theory of $p$-adic local systems on a variety of characteristic $p$ one has a naturally wider category that includes the category 
of local systems arising from geometries. The minimal slope constitution expands in its whole overconvergent $F$-isocrystal 
if it is irreducible. So the author asks a question that how one can directly recover the whole 
overconvergent $F$-isocrystal, e.g., its rank, from the minimal slope constitution of the slope filtrations.

Let $X$ be a smooth variety over $\mathrm{Spec}\, k$, and $\mathcal L$ a unit-root convergent $F$-isocrystal 
on $X/K$. Is $\mathcal L$ a minimal slope constitution of slope filtration of 
an irreducible overconvergent $F$-isocrystal? What is the essential image of the natural functor 
$$
     \left(\begin{array}{c} \mbox{\rm overconvergent $F$-isocrystals on $X/K$ } \\ 
           \mbox{\rm admitting slope filtration} \\\mbox{\rm with the minimal slope $0$}\end{array}\right) \\
           \rightarrow 
        \left(\begin{array}{c} \mbox{\rm unit-root convergent} \\ \mbox{\rm $F$-isocrystals on $X/K$}\end{array}\right)
$$
taking the minimal slope constitution of the slope filtrations? A unit-root convergent $F$-isocrystal 
belonging to the essential image is called ``geometric". For example, the overconvergent $F$-isocrystal $\mathcal{KL}_\psi$ 
on $\mathbb G_{\mathrm{m}\, \mathbb F_p}/\mathbb Q_p(\zeta_p)\, 
(\zeta_p: \mbox{\rm $p$-th root of unity})$ of rank $2$ corresponding to the variation of Kloosterman 
sums 
$$
    a \mapsto -\sum_{x \in k^\times}\psi\circ{tr}_{k/\mathbb F_p}(x + a/x)
$$
includes a rank one unit-root convergent $F$-isocrystal $\mathcal L_\psi$ on 
$\mathbb G_{\mathrm{m}\, \mathbb F_p}/\mathbb Q_p(\zeta_p)$ as a minimal slope constitution, where $\psi$ is a nontrivial additive character of $\mathbb F_p$ 
\cite{Dw74} \cite{Ha17} \cite[6.2]{Ts98b}. One can regard $\mathcal L_\psi$ as a unit-root 
convergent $F$-isocrystal on $\mathbb A^1_{\mathbb F_p}$, 
however the Frobenius acts by $1$ at $a = 0$ and it is not the minimal slope constitution of irreducible objects by the weight reason. 
Is there an explicit list of rank one geometric unit-root convergent $F$-isocrystals on $\mathbb G_{\mathrm{m}\, \mathbb F_p}$? 
In the list any positive tensor power $\mathcal L_\psi^{\otimes n}$ of $\mathcal L_\psi$ is included. 
We know such a geometric unit-root convergent $F$-isocrystal satisfies Dwork's conjecture on the meromorphy of unit-root 
$L$-functions \cite{Dw73} which was affirmatively solved by D.Wan \cite{Wa99} \cite{Wa00a} \cite{Wa00b}. 
In one sense our problem is a converse problem of Dwork's conjecture. 
Another direction is to characterize geometric unit-root convergent $F$-isocrystals by asymptotic behaviors of infinite towers of ramifications 
along the points at which the Newton polygons jump. The work of J.Kramer-Miller seems to be in this direction \cite{Kr16} \cite{Kr18} \cite{Kr19}. 

In conclusion the phenomena of slopes and their jumps are mysterious and interesting, which we should study. 

\subsection{Structure of this paper}

In section \ref{dJ} we fix the notation of the local settings and recall the opposite filtration of 
$\varphi$-modules. In section \ref{log} we introduce and study the local and global properties of the PBQ filtrations 
of overconvegent $F$-isocrystals on a curve. In section \ref{locver} we prove the local version of 
the minimal slope conjecture. In section \ref{satur} we introduce the notion of saturated overconvergent $F$-isocrystals 
in general dimension. 
In section \ref{curve} we estimate the rank of the maximal slope quotients applying results in section \ref{dJ} 
and prove the minimal slope conjecture in the case of curves. 
In section \ref{finite} we study the minimal slope conjecture for overconvergent $\overline{\mathbb Q}_p$-$F$-isocrystals on 
a smooth variety over a finite field. We put two appendices in order to understand overconvergent $F$-isocrystals well. 
We study Frobenius endomorphisms of (partially) 
$\dagger$-spaces, and the trace map of $F$-isocrystals with respect to finite base field changes. 

\vspace*{2mm}

\begin{center}
{\sc Acknowledgments}
\end{center}

The author began to take an interest in the minimal slope conjecture when he heard Emiliano Ambrosi's talk
on the minimal slope problem and its application in \cite{AD18}, which was mentioned above, 
in the conference ``$p$-adic cohomology and arithmetic geometry 2018" at Sendai. 
Not merely the talk but a couple of talks were related to the Frobenius slopes 
and the phenomena arising from gaps between overconvergent and convergent $F$-isocrystals. 
The studies stimulated the author.  

The author thanks Marco D'Addezio who explained his study on the minimal slope conjecture and the parabolicity conjecture.
The author also thanks Professor Takao Yamazaki for useful discussions, and an anonymous referee for
many helpful suggestions. 
This work was supported by JSPS KAKENHI Grant Number JP18H03667.

\section{Opposite slope filtrations}\label{dJ}

In this section we review the opposite filtration 
of $\varphi$-modules which was introduced by de Jong in \cite{dJ98}. 
The opposite filtration is one of the key ingredients to prove the minimal slope conjecture. 

\subsection{$\varphi$-$\nabla$-modules}\label{dJnot} 
Let $R$ be a complete discrete valuation ring with the fraction field $K$ of characteristic $0$ and 
the residue field $k$ which is an arbitrary field of characteristic $p > 0$. 
We suppose that there exists a $q$-Frobenius on $K$ for a positive power $q$ of $p$ 
and denote the $\sigma$-fixed subfield by $K_\sigma$ (see Appendix \ref{FrobK} for Frobenius endomorphisms). 
Let us fix notation as follows:
\begin{itemize}
\item $\mathcal E$ : the Amice ring over $K$, that is, 
$$
\mathcal E = \left\{\left. \sum_{n \in \mathbb Z}a_nt^n\, \right|\, \begin{array}{l}
     a_n \in K,\, \underset{n}{\sup}\, |a_n| < \infty, \\
     a_n \rightarrow 0\, \mathrm{as}\, n \rightarrow -\infty
     \end{array}\right\};
$$
\item $\mathcal R$ : the Robba ring over $K$, that is, 
$$
     \mathcal R = \left\{ \left. \sum_{n \in \mathbb Z}a_nt^n\, \right|\, 
     \begin{array}{l} 
     a_n \in K, \\
     0 < \exists \eta < 1\, \mbox{such that}\, 
     |a_n|\eta^n \rightarrow 0\, 
     \mathrm{as}\, n \rightarrow -\infty,  \\ 
     0 < \forall \xi < 1\, \mbox{such that}\, 
     |a_n|\xi^n \rightarrow 0\, 
     \mathrm{as}\, n \rightarrow \infty
     \end{array} \right\};
$$
\item $\mathcal E^\dagger$ : the bounded Robba ring over $K$, that is, 
$$
     \mathcal E^\dagger = \left\{ \left. \sum_{n \in \mathbb Z}a_nt^n \in \mathcal R\, \right|\, \underset{n}{\sup}\, |a_n| < \infty \right\};
$$
\item $K[\vspace*{-0.2mm}[t]\vspace{-0.2mm}]_0$ : the $K$-algebra of bounded functions on the unit disc $D(0, 1^-) = \{|t| < 1\}$, that is, 
$$
K[\vspace*{-0.2mm}[t]\vspace{-0.2mm}]_0 = K\otimes_RR[\vspace*{-0.2mm}[t]\vspace{-0.2mm}] 
= \left\{ \left. \sum_{n=0}^{\infty}a_nt^n \in K[\vspace*{-0.2mm}[t]\vspace{-0.2mm}]\, \right|\, 
\sup_n |a_n| < \infty \right\}. 
$$
\end{itemize}

\noindent
Let $B$ be one of $\mathcal E, \mathcal R, \mathcal E^\dagger$ and $K[\vspace*{-0.2mm}[t]\vspace{-0.2mm}]_0$. Then $B$ is furnished with 
\begin{itemize}
\item $\left|\sum a_nt^n\right| = \sup_n|a_n|$, which is called the Gauss norm for $\sum_na_nt^n \in B$; 
\item $d : B \rightarrow Bdt$, which is a continuous $K$-derivation $d(\sum_n a_nt^n) = \sum_nna_nt^{n-1}dt$; 
\item $\varphi : B \rightarrow B$, which is a $q$-Frobenius such that 
$\varphi|_K = \sigma$ and $\varphi(t) \equiv t^q\, (\mathrm{mod}\, \mathbf m)$. 
If we treat $\mathcal E, \mathcal E^\dagger$ and $K[\vspace*{-0.2mm}[t]\vspace{-0.2mm}]_0$ 
(resp. $\mathcal E$ and $\mathcal E^\dagger$, resp. $\mathcal R$ and $\mathcal E^\dagger$) in the same time, then 
the Frobenius on $\mathcal E$ and $\mathcal E^\dagger$ (resp. $\mathcal E$, resp. $\mathcal R$) 
is an extension of that of $K[\vspace*{-0.2mm}[t]\vspace{-0.2mm}]_0$ (resp. $\mathcal E^\dagger$, resp. $\mathcal E^\dagger$). 
\end{itemize}

\noindent
$\mathcal E$ (resp. $\mathcal E^\dagger$) is a discrete valuation field 
with residue field $k(\hspace*{-0.3mm}(t)\hspace*{-0.3mm})$ under the Gauss norm, and $K[\vspace*{-0.2mm}[t]\vspace{-0.2mm}]_0$ 
is a principal ideal domain. 
We denote the integer ring of $\mathcal E$ (resp. $\mathcal E^\dagger$) by 
$\mathcal O_{\mathcal E}$ (resp. $\mathcal O_{\mathcal E^\dagger}$). 
Then $\mathcal O_{\mathcal E}$ (resp. $\mathcal O_{\mathcal E^\dagger}$) is a complete (resp. Henselian) 
discrete valuation ring and $\mathcal E^\dagger$ is algebraically closed in $\mathcal E$.

\begin{definition} Let $B$ be one of $\mathcal E, \mathcal E^\dagger,$ and $K[\vspace*{-0.2mm}[t]\vspace{-0.2mm}]_0$. 
\begin{enumerate}
\item A $\varphi$-module over $B$
is a free $B$-module $M$ of finite rank which is furnished with 
a $B$-linear bijection $\varphi_M : \varphi^\ast M \rightarrow M$.  
\item A $\nabla$-module over $B$
is a free $B$-module $M$ of finite rank which is furnished with 
a $K$-linear map $\nabla_M : M \rightarrow Mdt$ such that $\nabla_M(am) = mda + a\nabla(m)$ for $a \in B$ and $m \in M$. 
\item A $\varphi$-$\nabla$-module over $B$ 
is a triplet $(M, \nabla_M, \varphi_M)$ such that $(M, \varphi_M)$ is a $\varphi$-module over $B$ and $(M, \nabla_M)$ 
is a $\nabla$-module over $B$ satisfying the diagram 
$$
       \begin{array}{ccc}
             \varphi^\ast M &\overset{\varphi^\ast\nabla}{\rightarrow} &\varphi^\ast(Mdt) \\
             \varphi_M \downarrow \hspace*{7mm} & &\hspace*{10mm} \downarrow \varphi_M\otimes\varphi \\
             M &\underset{\nabla}{\rightarrow} &Mdt
        \end{array}
$$
is commutative. 
\end{enumerate}
$\varphi_M$ is called Frobenius and $\nabla_M$ is called a connection. 
We denote the category of $\varphi$-$\nabla$-modules over $B$ 
by $\mbox{\rm \bf $\Phi$M}^\nabla_B$. 
\end{definition}

\begin{definition} Let $M$ be a $\varphi$-$\nabla$-module over $\mathcal E$. 
\begin{enumerate}
\item An $\mathcal E^\dagger$-lattice of $M$ is a $\varphi$-$\nabla$-module $M^\dagger$ 
over $\mathcal E^\dagger$ 
such that $M \cong \mathcal E\otimes_{\mathcal E^\dagger}M^\dagger$ as 
$\varphi$-$\nabla$-modules $M$ over $\mathcal E$. 
\item A $K[\hspace*{-0.2mm}[t]\hspace*{-0.2mm}]_0$-lattice of $M$ is a $\varphi$-$\nabla$-module $M_0$ 
over $K_\alpha[\hspace*{-0.2mm}[x_\alpha]\hspace*{-0.2mm}]_0$ 
such that $M \cong \mathcal E\otimes_{K[\hspace*{-0.2mm}[t]\hspace*{-0.2mm}]_0}M_0$ as 
$\varphi$-$\nabla$-modules over $\mathcal E$. 
One can also define the notion of $K[\hspace*{-0.2mm}[t]\hspace*{-0.2mm}]_0$-lattices of 
$\varphi$-$\nabla$-modules over $\mathcal E^\dagger$
\end{enumerate}
\end{definition}

The category $\mbox{\rm \bf $\Phi$M}^\nabla_B$ 
is a $K_\sigma$-linear Abelian category which is furnished with tensor products $\otimes_B$ and duals $M \mapsto M^\vee$ (see \cite[Section 3]{Ts96}). 
The category does not depend on the choice of $q$-Frobenius $\varphi$ \cite[Theorem 3.4.10]{Ts98b}. 
The forgetful functor 
$$
         \mbox{\rm \bf $\Phi$M}^\nabla_{\mathcal E^\dagger} \rightarrow  \mbox{\rm \bf $\Phi$M}^\nabla_{\mathcal E} 
         \hspace*{3mm} M^\dagger \mapsto M =\mathcal E\otimes_{\mathcal E^\dagger} M^\dagger
$$
is fully faithful by \cite[Theorem 5.1]{Ke08}. We use the notation $M$ (resp. $M^\dagger$) for a $\varphi$-$\nabla$-module over $\mathcal E$ 
(resp. $\mathcal E^\dagger$), and for a $\varphi$-$\nabla$-module $N^\dagger$ over $\mathcal E^\dagger$ the image by the above 
forgetful functor is denoted by $N$. 
Note that the category of $\varphi$-modules over $\mathcal E^\dagger$ (resp. $\mathcal E$) 
is a $K_\sigma$-linear Abelian category, 
but, except unit-root objects below, it may depends on the choice of Frobenius and the natural functor from the category over $\mathcal E^\dagger$ 
to that over $\mathcal E$ is not fully faithful (see \cite[Remark 2.2.7]{Ts96} for a counter example).

\begin{definition} 
\begin{enumerate}
\item A unit-root $\varphi$-module over $\mathcal E$ is a $\varphi$-module $M$ over $\mathcal E$ 
such that there exists a finite free $\mathcal O_{\mathcal E}$-submodule $L$ of $M$ satisfying  
$\mathcal E \otimes_{\mathcal O_\mathcal E} L = M$ and $\varphi(L)$ generates $L$ as an $\mathcal O_{\mathcal E}$-module. 
A $\varphi$-module over $\mathcal E^\dagger$ (resp. a $\varphi$-$\nabla$-module over $\mathcal E$, 
resp. a $\varphi$-$\nabla$-module over $\mathcal E^\dagger$) is unit-root if so is the associate $\varphi$-module over $\mathcal E$. 
\item A $\varphi$-module $M$ over $\mathcal E$ is purely of slope $\frac{m}{n}\, (m, n \in \mathbb Z, n > 0)$ 
if $M^{\otimes n}(m)$ is a unit-root $\varphi$-module over $\mathcal E$. 
Here $M^{\otimes n}(m)$ means a tensor product of $n$-copies of $M$ with a Frobenius 
$\varphi_{M^{\otimes n}(m)} = q^{-m}\varphi^{\otimes n}$. 
A $\varphi$-module over $\mathcal E^\dagger$ (resp. a $\varphi$-$\nabla$-module over $\mathcal E$, 
resp. a $\varphi$-$\nabla$-module over $\mathcal E^\dagger$) is pure of slope $m/n$ 
if so is the associate $\varphi$-module over $\mathcal E$. 
\end{enumerate}
\end{definition}

\begin{proposition}\label{slfilg} \mbox{\rm (See \cite[Remark 1.7.8]{Ke08},  \cite[Theorem 2.4]{CT09}.)} 
\begin{enumerate}
\item Let $M$ be 
a $\varphi$-module over $\mathcal E$. Then $M$ admits a slope filtration 
$$
     M = M^0 \supsetneq M^1 \supsetneq \cdots \supsetneq M^{r-1} \supsetneq M^r = 0
$$
as $\varphi$-modules over $\mathcal E$, namely, $M^i/M^{i+1}$ has a unique Frobenius slope $s^i$ 
with $s^0 > s^1 > \cdots > s^{r-1}$. We call $s_0$ the maximal slope of $M$, and $M/M^1$ the maximal slope quotient of $M$. 
\item If $M$ is a $\varphi$-$\nabla$-module over $\mathcal E$, then the 
slope filtration as $\varphi$-modules is a filtration as $\varphi$-$\nabla$-modules over $\mathcal E$. 
\end{enumerate}
We use the notation above for the slope filtrations of $\varphi$-$\nabla$-modules over $\mathcal E$. 
\end{proposition}

\subsection{Generalized series}\label{genser}
In this subsection we assume that the residue field $k$ of $K$ is algebraically closed for simplicity, 
and $\sigma$ is a $q$-Frobenius on $K$ for a positive power $q$ of $p$ such that the natural map 
$$
       K_\sigma \otimes_{W(\mathbb F_q)}W(k) \rightarrow K
$$
is an isomorphism of fields, where $W(k)$ is the Witt-vectors ring with coefficients in $k$. 
We also fix a $q$-Frobenius $\varphi$ on $\mathcal E$ and $\mathcal E^\dagger$ defined by 
$$
     \varphi\left(\sum_{n \in \mathbb Z}a_nt^n\right) = \sum_{n \in \mathbb Z}\sigma(a_n)t^{qn}. 
$$
We introduce Hahn-Mal'cev-Neumann's generalized power series rings $\widetilde{\mathcal E}$, $\widetilde{\mathcal E}^\dagger$ 
and give their fundamental properties. 
These rings were studied in \cite[Section 4]{Ke01} \cite[Section 4]{Ke05} \cite[Section 2]{Ke08}. 
Let us regard $\mathbb Q$ as a well-ordered Abelian group with respect to the usual Archimedean order $\leq$, 
and define 
\begin{itemize}
\item $\widetilde{\mathcal E} = \left\{ \sum_{n \in \mathbb Q} a_nt^n\, \left|\, 
\begin{array}{l} a_n \in K, \sup_n|a_n| < \infty, |a_n|\rightarrow 0\, (\mbox{\rm as}\, n \rightarrow - \infty), \\
\mbox{\rm the support set}\, \{ n \in \mathbb Q\, |\, |a_n| \geq \delta \}\, \mbox{\rm is well-ordered for any $\delta \in |K^\times|$.} 
\end{array} \right. \right\}$
\item $\widetilde{\mathcal E}^\dagger = \left\{ \left. \sum_{n \in \mathbb Q} a_nt^n \in \widetilde{\mathcal E}\, \right|\, 
|a_n|\eta^n \rightarrow 0\, (\mbox{\rm as}\, n \rightarrow -\infty)\, 
\mbox{\rm for some $0 < \eta < 1$.} \right\}$, 
\item $k(\hspace*{-0.2mm}(t^{\mathbb Q})\hspace*{-0.2mm}) = \left\{ \sum_{n \in \mathbb Q} a_nt^n\, \left|\, 
\begin{array}{l} a_n \in k, \\
\mbox{\rm the support set}\, \{ n \in \mathbb Q\, |\, a_n \ne 0 \}\, \mbox{\rm is well-ordered.} 
\end{array} \right. \right\}$
\item $\widetilde{\varphi}$ : the $q$-Frobenius on $\widetilde{\mathcal E}$ (resp. $\widetilde{\mathcal E}^\dagger$) defined by 
$$
     \widetilde{\varphi}\left(\sum_{n \in \mathbb Q}a_nt^n\right) = \sum_{n \in \mathbb Q}\sigma(a_n)t^{qn}.
$$
\end{itemize}
Here the sum and the product in $\widetilde{\mathcal E}$ (resp. $\widetilde{\mathcal E}^\dagger$, 
resp. $k(\hspace*{-0.2mm}(t^{\mathbb Q})\hspace*{-0.2mm})$) are defined by 
$$
  \begin{array}{l}
    \sum_na_nt^n + \sum_nb_nt^n = \sum_n(a_n + b_n)t^n \\
        \sum_na_nt^n \times \sum_nb_nt^n = \sum_n\left(\sum_{l+m=n}a_lb_m\right)t^n 
        \end{array}
$$
which are well-defined by the following lemmas since $K$ is a complete discrete valuation field. 
Indeed, the boundedness of elements of $\widetilde{\mathcal E}$ 
and Lemma \ref{word} (1) (resp. (2)) 
below implies the well-definedness of products (resp. the existence of inverses for nonzero elements). 
The case of $k(\hspace*{-0.2mm}(t^{\mathbb Q})\hspace*{-0.2mm})$ is similar and 
it is algebraically closed by \cite[Proposition 1]{Ke01}. 
$\widetilde{\mathcal E}^\dagger$ is a $K$-subalgebra (resp. a subfield) 
of $\widetilde{\mathcal E}$ by  Lemma \ref{ineq} (3), (4) (resp. (5) and Lemma \ref{word} (3)). 

\begin{lemma}\label{word} Suppose $I, J$ are well-ordered subsets of $\mathbb Q$. 
\begin{enumerate}
\item For any $n \in \mathbb Q$, the set $\{(i, j) \in I\times J\, |\, i + j = n \}$ is finite. 
\item The sets $I \cup J$ and $I+J = \{ i + j\, |\, i \in I,\, j \in J \}$ are well-ordered.
\item Suppose $I \subset \mathbb Q_{> 0}$. 
If $mI$ denotes the sum $I + \cdots + I$ of $m$ copies of $I$, then $\cup_{m \geq 1}mI$ is well-ordered. 
\end{enumerate}
\end{lemma}

\begin{lemma}\label{ineq} For an element $a = \sum_n a_nt^n \in \widetilde{\mathcal E}$, let us define a map 
$N_a : \mathbb Z \rightarrow \mathbb Q \cup \{ \infty \}$ by 
$$
     N_a(l) = \min\{ n \in \mathbb Q\, |\, \mathrm{ord}_\pi(a_n) \leq l \} \cup \{ \infty\}. 
$$
Here $\pi$ is a uniformizer of $K$, $\mathrm{ord}_\pi$ is an additive valuation of $K$ normalized by 
$\mathrm{ord}_\pi(\pi) = 1$, 
$\mathbb Q \cup \{ \infty \}$ is furnished with the standard order $<$ such that $n < \infty$ for all $n \in \mathbb Q$, 
and $n + \infty = \infty + n = \infty$ for any $n \in \mathbb Q \cup \{ \infty\}$. 
\begin{enumerate}
\item If $l \geq m$, then $N_a(l) \leq N_a(m)$. 
\item $N_{\varphi(a)}(l) = qN_a(l)$ for any $l \in \mathbb Z$.  
\item For $a, b \in \widetilde{\mathcal E}$, the inequality 
$N_{ab}(l) \geq \underset{i+j = l}{\inf}(N_a(i) + N_b(j))$ holds for any $l \in \mathbb Z$. 
Note that it is sufficient to take only finite numbers of $(i, j) \in \mathbb Z \times \mathbb Z$ to compute the infimum above since the coefficients of 
$a$ and $b$ are bounded and the valuation of $K$ is discrete. 
\item For an element $a \in \widetilde{\mathcal E}$, the following are equivalent. 
\begin{list}{}{}
\item[\mbox{\rm (i)}] $a \in \widetilde{\mathcal E}^\dagger$.
\item[\mbox{\rm (ii)}] There exist $c, d \in \mathbb Q$ with $d > 0$ such that $N_a(l) \geq c - dl$ for any $l \in \mathbb Z$. 
\item[\mbox{\rm (iii)}] $N_a(l)/l$ is lower bounded on $l > 0$. 
\end{list}
\item Suppose $a = \sum_na_nt^n \in \widetilde{\mathcal E}^\dagger$ such that $|a_n| \leq 1$ for all $n \in \mathbb Q$ 
and $|a_n| < 1$ for all $n < 0$. Then there exists $d \in \mathbb Q_{> 0}$ such that 
$N_a(l) \geq -dl$ for all $l \in \mathbb Z$. Moreover, $N_{a^m}(l) \geq -dl$ for any $l$. 
\end{enumerate}
\end{lemma}

\prf{(4) Let $a = \sum_n a_nt^n \in \widetilde{\mathcal E}$. 
Suppose $a \in \widetilde{\mathcal E}^\dagger$. Then 
there exist $C, D \in \mathbb Q$ with $D>0$ such that $\mathrm{ord}_\pi(a_n) \geq C - Dn$ for any $n$. 
Then $N_a(l) \geq C/D - l/D$ for any $l$. The converse also holds. Hence we have (i) $\Leftrightarrow$ (ii). 
(ii) $\Rightarrow$ (iii) is trivial. Since $\mathrm{sup}_n|a_n| < \infty$, the converse (iii) $\Rightarrow$ (ii) holds. 
}

\begin{lemma}\label{Efield}
\begin{enumerate}
\item $\widetilde{\mathcal E}$ is a complete discrete valuation field under the valuation 
$$
      \left|\sum_n a_nt^n\right| = \sup_n|a_n| 
$$
with the residue field $k(\hspace*{-0.2mm}(t^{\mathbb Q})\hspace*{-0.2mm})$. 
\item $\widetilde{\mathcal E}^\dagger$ is a discrete valuation field such that the completion 
under the valuation in (1) is $\widetilde{\mathcal E}$ and the integer ring $\mathcal O_{\widetilde{\mathcal E}^\dagger}$ 
is Henselian. 
\item The following is a commutative diagram of extensions of discrete valuation fields with the same valuation group
$$
       \begin{array}{ccccc}
            &\nearrow &\mathcal E^\dagger &\rightarrow &\widetilde{\mathcal E}^\dagger \\
            K& &\downarrow & &\downarrow \\
            &\searrow &\mathcal E &\rightarrow &\widetilde{\mathcal E}
        \end{array}
$$
such that $\widetilde{\mathcal E}^\dagger \cap \mathcal E = \mathcal E^\dagger$ in $\widetilde{\mathcal E}$ 
and the natural morphism 
$$
        \widetilde{\mathcal E}^\dagger \otimes_{\mathcal E^\dagger}\mathcal E \rightarrow \widetilde{\mathcal E}
$$
is injective. 
\item The $q$-Frobenius endomorphisms act compatibly in the diagram above. 
If $(-)_\varphi$ denotes the $\varphi$-fixed subset of $(-)$, then 
$$
         \widetilde{\mathcal E}_{\widetilde{\varphi}} = (\widetilde{\mathcal E}^\dagger)_{\widetilde{\varphi}} 
         = \mathcal E_\varphi  = (\mathcal E^\dagger)_\varphi = K_\sigma.
$$   
\end{enumerate}
\end{lemma}

\prf{(1) $\widetilde{\mathcal E}$ is complete by Lemma \ref{word} (2).

(2) We will prove $\mathcal O_{\widetilde{\mathcal E}^\dagger}$ is Henselian. The rest is easy. 
Let $f(x) = x^s + a_1x^{s-1} + \cdots + a_s \in \mathcal O_{\widetilde{\mathcal E}^\dagger}[x]$ be a monic polynomial such that 
the image $\overline{f}(x)$ of $f(x)$ in $k(\hspace*{-0.2mm}(x^{\mathbb Q})\hspace*{-0.2mm})$ has a simple root 
$\overline{\alpha} \in k(\hspace*{-0.2mm}(x^{\mathbb Q})\hspace*{-0.2mm})$. 
We may assume that $\overline{\alpha} = 0$. 
Since the residue field $\mathcal O_{\widetilde{\mathcal E}^\dagger}$ is 
the algebraically closed field $k(\hspace*{-0.2mm}(x^{\mathbb Q})\hspace*{-0.2mm})$, 
we may also assume that $\overline{f}'(\overline{\alpha}) = 1$. 
Let us take a unique element $\alpha$ in the integer ring $\mathcal O_{\widetilde{\mathcal E}}$ 
of $\widetilde{\mathcal E}$ such that 
$f(\alpha) = 0$ and $\alpha\, (\mathrm{mod}\, \mathbf m\mathcal O_{\mathcal E}) 
= \overline{\alpha} = 0$. Indeed, such an $\alpha$ exists 
because the integer ring $\mathcal O_{\widetilde{\mathcal E}}$ is complete. 
Let us fix constants $c, d \in \mathbb Q$ with $c \leq 0, d > 0$ 
such that $N_{a_i}(l) \geq c - dl$ for any $l \in \mathbb Z$ and all $i$ 
by Lemma \ref{ineq}, where $a_0=1$ is the coefficient of $x^n$. Our claim is the inequality 
$$
     N_\alpha(l) \geq -c +(2c-d)l 
$$
for any $l \in \mathbb Z$ by induction on $l$. 
If $l \leq 0$, then $N_\alpha(l) = \infty$. 
Suppose the equality holds for $l \leq h$. Then there exist $\alpha_h, \beta \in \mathcal O_{\widetilde{\mathcal E}}$ 
with $\alpha = \alpha_h + \pi^{h+1}\beta$ such that 
$N_{\alpha_h}(l) \geq -c + (2c - d)l$ for any $l \leq h$ and 
$N_{\alpha_h}(h) = N_{\alpha_h}(h+1) = N_{\alpha_k}(h+2) \cdots$. Since 
$$
   \begin{array}{lll}
       N_{f(\alpha_h)}(h+1) &\geq &\underset{i}{\min}\, N_{a_i\alpha_h^{s-i}}(h+1) \\
       &\geq &\underset{i}{\min}\, \underset{\tiny \begin{array}{c} l_1+l_2 = h+1 \\ l_1, l_2 \geq 0 \end{array}}{\min}\left(N_{a_i}(l_1) + N_{\alpha_h^{s-i}}(l_2)\right)\\
              &\geq &\underset{i}{\min}\, \underset{\tiny \begin{array}{c} l_1+l_2 = h+1 \\ l_1, l_2 \geq 0 \end{array}}{\min} 
              \underset{\tiny \begin{array}{c} m_1 + \cdots + m_{s-i} = l_2 \\ m_j \geq 0 \end{array}}{\min}\left(N_{a_i}(l_1) + \sum_{j=1}^{s-i}N_{\alpha_h}(m_j)\right)\\
       &\geq &\min\left\{ c-dl_1-(s-i)c + (2c-d)l_2\, \left|\, \begin{array}{c} 0 \leq i \leq s \\ l_1+l_2 = h+1, l_2 \geq 0 \\ 
       l_1 \left\{\begin{array}{ll} \geq 0\, &\mathrm{if}\,\, i \leq s-2, \\ > 0\, &\mathrm{if}\,\, i = s-1,\\ =h+1\, &\mathrm{if}\,\, i =s \end{array}\right.\end{array}\right. \right\} \\
        &\geq &-c + (2c-d)(h+1)
   \end{array}
$$
by Lemma \ref{ineq}, the congruence 
$f(\alpha_h)+\pi^{h+1}\beta f'(\alpha_h)\, \equiv\, f(\alpha_h+\pi^{h+1}\beta)\, \equiv 0\,  (\mathrm{mod}\, \pi^{h+2})$ implies an inequality 
$$
      N_{\pi^{h+1}\beta}(h+1) \geq - c + (2c-d)(h+1). 
$$
Define $\alpha_{h+1}$ by $\alpha_h +\pi^{h+1}\beta$ removing the terms whose coefficients are $0$ modulo $\pi^{h+2}$. 
Then $\alpha_{h+1}$ satisfies $N_{\alpha_{h+1}}(l) \geq - c + (2c - d)l$ for any $l \leq h + 1$ and 
$N_{\alpha_{h+1}}(h+1) = N_{\alpha_{h+1}}(h+2) = N_{\alpha_{h+1}}(h+3) \cdots$. 
Therefore, $\alpha$ belongs to $\widetilde{\mathcal E}^\dagger$ and $\mathcal O_{\widetilde{\mathcal E}^\dagger}$ is Henselian. 

(3) We prove the injectivity. For $b_1, b_2, \cdots, b_r \in \mathcal E$, 
suppose $a_1b_1 + \cdots +a_rb_r = 0$ in $\widetilde{\mathcal E}$ for $a_1, \cdots, a_r \in \widetilde{\mathcal E}^\dagger$. 
We may assume that the union of the support sets of $a_1, \cdots, a_r$ has a nontrivial intersection with $\mathbb Z$. 
Take subseries $a_i' = \sum_{n \in \mathbb Z}a_{i, n}t^n$ of $a_i = \sum_{n \in \mathbb Q}a_{i, n}t^n$. 
Then $a_i' \in \mathcal E^\dagger$ for any $i$ and $a_1'b_1 + \cdots +a_r'b_r = 0$. 
Hence $b_1, b_2, \cdots, b_r$ are linearly independent over $\widetilde{\mathcal E}^\dagger$ 
if and only if they are linearly independent over $\mathcal E^\dagger$. 

(4) By the universal property of Witt vectors rings, 
there exists a canonical isomorphism 
 $\widetilde{\mathcal E} = K_\sigma\otimes_{W(\mathbb F_q)}W(k(\hspace*{-0.2mm}(t^{\mathbb Q})\hspace*{-0.2mm}))$ 
 with respect to the $q$-Frobenius $\varphi = \mathrm{id}_{K_\sigma}\otimes \mathrm{Frob}^s$ 
 where $\mathrm{Frob}$ is the $p$-Frobenius of the Witt-vectors ring and $q = p^s$. 
 The rest easily follow from it. 
}

\begin{remark} The field $\widetilde{\mathcal E}$ \mbox{\rm (resp. $\widetilde{\mathcal E}^\dagger$)} plays a similar role with $\Gamma_2$ 
\mbox{\rm (resp. $\Gamma_{2, c}$)} in \cite[Section 4]{dJ98}. $\widetilde{\mathcal E}$ \mbox{\rm (resp. $\widetilde{\mathcal E}^\dagger$)} 
is ``larger" than $\Gamma_2$ \mbox{\rm (resp. $\Gamma_{2, c}$)}. 
\end{remark}

The following lemma is a generalization of \cite[Proposition 2.2.2]{Ts96} which asserts 
the similar claim for Frobenius equations over $\mathcal E^\dagger$.  

\begin{proposition}\label{eequat} Let $a_1, \cdots, a_s$ be elements in $\mathcal O_{\widetilde{\mathcal E}^\dagger}$. 
Suppose $y = \sum_ny_nt^n\in \widetilde{\mathcal E}$ satisfies an Frobenius equation 
$$
      \widetilde{\varphi}^s(y) + a_1\widetilde{\varphi}^{s-1}(y) + \cdots + a_sy = 0. 
$$
Then $y$ belongs to $\widetilde{\mathcal E}^\dagger$. 
\end{proposition}

\prf{We may assume that $|y| = 1$. 
Replacing $y$ by $t^my$ for a sufficient large $m \geq 0$, we may assume that $|y_n| < 1$ for any $n < 0$ and that $|a_{i, n}| < 1$ 
for any $n < 0$ where $a_i = \sum_na_{i, n}t^n$. Then there exists $d \in \mathbb Q$ with $d > 0$ such that 
$N_{a_i}(l) \geq -dl$ for any $l \in \mathbb Z$, where $d$ does not depend on $i$. 
Let us estimate $N_{a_i\varphi^{s-i}(y)}(l)/l$ on $l \in \mathbb Z$ when $a_i \ne 0$. 
If $l$ is a sufficiently large positive integer, then 
$$
   N_{a_i\varphi^{s-i}(y)}(l)/l \geq \frac{1}{l}\underset{\tiny \begin{array}{c} l_1 + l_2= l \\ l_1, l_2 \geq 0 \end{array}}{\mathrm{min}}
     \left(N_{a_i}(l_1)+ N_{\varphi^{s-i}(y)}(l_2)\right) \\
    \geq \underset{0 < l_2 \leq l}{\mathrm{min}}\left\{-d + \frac{N_{\varphi^{s-i}(y)}(l_2)}{l_2}\right\} \cup \{ -d\} 
$$
by Lemma \ref{ineq}. Suppose there exists a positive integer $m$ such that 
\begin{list}{}{}
\item[(i)] $N_y(l)/l > N_y(m)/m$ for any $l < m$ and 
\item[(ii)] $N_y(m)/m < -d$. 
\end{list}
If such an $m$ does not exist, then $N_y(l)/l$ is bounded on $l > 0$ and $y \in \widetilde{\mathcal E}^\dagger$ by Lemma \ref{word} (4). 
Since $y$ satisfies the given Frobenius equation, the following estimate 
$$
    q^sN_y(m)/m =  N_{\varphi^s(y)}/m = N_{a_1\varphi^{s-1}(y) + \cdots +a_sy}(m)/m \geq -d + q^{s-1}N_y(m)/m 
$$
holds. Then it implies $N_y(m)/m \geq -d/(q^s - q^{s-1}) \geq -d$ and contradicts the hypothesis on $m$. 
Therefore $y \in \widetilde{\mathcal E}^\dagger$. 
}

\subsection{Opposite filtrations}\label{oppofil} Let us recall the opposite filtration in our context. 

\begin{proposition}{\mbox{\rm (\textit{cf.} \cite[Proposition 5.5 (ii)]{dJ98})}}\label{unitdec} 
\begin{enumerate}
\item Let $M$ be a unit-root $\varphi$-module 
over $\mathcal E$, put $\widetilde{M} =  \widetilde{\mathcal E} \otimes_{\mathcal E} M$ 
to be a $\widetilde{\varphi}$-module over $\widetilde{\mathcal E}$. Then 
$\widetilde{M}$ is isomorphic to the trivial $\widetilde{\varphi}$-module 
$\widetilde{\mathcal E}^{\oplus \mathrm{dim}_{\mathcal E}M}$. 
\item Let $N^\dagger$ be a unit-root $\varphi$-module 
over $\mathcal E^\dagger$, put $\widetilde{N}^\dagger =  \widetilde{\mathcal E}^\dagger \otimes_{{\mathcal E}^\dagger} N^\dagger$ 
to be a $\widetilde{\varphi}$-module over $\widetilde{\mathcal E}^\dagger$. 
Then 
$\widetilde{N}^\dagger$ is isomorphic to the trivial $\widetilde{\varphi}$-module 
$(\widetilde{\mathcal E}^\dagger)^{\oplus \mathrm{dim}_{\mathcal E^\dagger}N^\dagger}$. 
\end{enumerate}
\end{proposition}

\prf{(1) Since the residue field $k(\hspace*{-0.2mm}(t^{\mathbb Q})\hspace*{-0.2mm})$ 
of $\widetilde{\mathcal E}$ is algebraically closed by \cite[Proposition 1]{Ke01}, 
it contains an algebraic closure of the residue field $k(\hspace*{-0.2mm}(t)\hspace*{-0.2mm})$ of $\mathcal E$. 
Since $\widetilde{\mathcal E}$ is $p$-adically complete, the assertion holds by 
Dieudonn\'e-Manin's classification of $F$-spaces. 

(2) Take a cyclic vector $e$ of the $\varphi$-module $N^\dagger$ \cite[Lemma 3.1.4]{Ts96} 
and let $A$ be a representation matrix of the Frobenius $\varphi_{N^\dagger}$ 
with respect to the basis $e, \varphi_{N^\dagger}(e), \varphi_{N^\dagger}^2(e), \cdots, \varphi_{N^\dagger}^{s-1}(e)$, 
that is, 
$$
   \begin{array}{c}
     \varphi_{N^\dagger}\left(e, \varphi_{N^\dagger}(e), \varphi_{N^\dagger}^2(e), \cdots, \varphi_{N^\dagger}^{s-1}(e)\right)
     =(e, \varphi_{N^\dagger}(e), \varphi_{N^\dagger}^2(e), \cdots, \varphi_{N^\dagger}^{s-1}(e))A \\
     A = \left(\begin{array}{ccccc}
     &&& &a_s \\
     1&& & &a_{s-1} \\
     &1& &&a_{s-2} \\
     & &\ddots & &\vdots \\
     & & &1&a_1
     \end{array}\right) \in \mathrm{GL}_s(\mathcal O_{\mathcal E^\dagger}). 
     \end{array}
$$
One can find a matrix 
$Y \in \mathrm{GL}_s(\widetilde{\mathcal E})$ such that $A\varphi(Y) = Y$ by (1). If we put 
$Y^{-1} = (z_{i, j})$, then 
$$
     \varphi(z_{i, j}) = z_{i, j+1}\, (1 \leq j \leq s-1), \varphi(z_{i, s}) = a_1z_{i, s} + a_2z_{i, s-1} + \cdots + a_sz_{i, 1}. 
$$
Hence $z_{i, 1}$ satisfies the Frobenius equation 
$$
    \varphi^s(z_{i, 1}) = a_1\varphi^{s-1}(z_{i, 1}) + a_2\varphi^{s-2}(z_{i, 1}) + \cdots + a_sz_{i, 1}
$$
for any $i$. Since the unit-root condition implies $a_1, \cdots, a_s \in \mathcal O_{\mathcal E^\dagger}$ and $|a_s| = 1$, 
we have $z_{i, s}$ is included in $\widetilde{\mathcal E}^\dagger$ by Proposition \ref{eequat} and 
so that $Y \in \mathrm{GL}_s(\widetilde{\mathcal E}^\dagger)$. Therefore, $\widetilde{N}^\dagger$ 
is a trivial $\widetilde{\varphi}$-module over $\widetilde{\mathcal E}^\dagger$.}

\begin{proposition}{\mbox{\rm (\textit{cf.} \cite[Proposition 5.5]{dJ98})}}\label{opfil} Let $N^\dagger$ be a $\varphi$-module 
over $\mathcal E^\dagger$, put $\widetilde{N}^\dagger =  \widetilde{\mathcal E}^\dagger \otimes_{{\mathcal E}^\dagger} N^\dagger$ 
to be a $\widetilde{\varphi}$-module over $\widetilde{\mathcal E}^\dagger$. 
Then there is a filtration 
$$
      0 = \widetilde{N}_0^\dagger \subsetneq \widetilde{N}_1^\dagger \subsetneq \cdots \subsetneq \widetilde{N}_r^\dagger = \widetilde{N}^\dagger, 
$$
of $\widetilde{N}^\dagger$ satisfying 
\begin{list}{}{}
\item[\mbox{\rm (i)}] $\widetilde{N}_i^\dagger/\widetilde{N}_{i-1}^\dagger$ is pure of slope $\lambda_i$; 
\item[\mbox{\rm (ii)}] $\lambda_1 > \lambda_2 \cdots > \lambda_r$.
\end{list}
Moreover, if $N = N^0 \supsetneq N^1 \supsetneq \cdots \supsetneq N^r = 0$ be a slope filtration of $N = \mathcal E \otimes_{\mathcal E^\dagger} N^\dagger$, 
then $s = r$ and 
$$
    \widetilde{\mathcal E}\otimes_{\widetilde{\mathcal E}^\dagger} \widetilde{N}_i^\dagger/\widetilde{N}_{i-1}^\dagger \cong 
    \widetilde{\mathcal E}\otimes_{\mathcal E} N^{i-1}/N^i
$$
for all $i$ as $\widetilde{\varphi}$-modules over $\widetilde{\mathcal E}$. 
Such an increasing filtration $\{ \widetilde{N}^\dagger_i\}$ is called the opposite slope filtration of $\widetilde{N}^\dagger$. 
\end{proposition}

\prf{Since the residue field $k(\hspace*{-0.2mm}(t^{\mathbb Q})\hspace*{-0.2mm})$ of 
the complete discrete valuation field $\widetilde{\mathcal E}$ is algebraically closed, we have a decomposition 
$\widetilde{N} = \widetilde{\mathcal E}\otimes_{\mathcal E^\dagger}N^\dagger \cong \oplus_i\widetilde{\mathcal E}\otimes_{\mathcal E} N^{i-1}/N^i$ 
as $\widetilde{\varphi}$-modules. 
Then the assertion follows from Lemma \ref{indbl} below.
}

\begin{lemma}{\mbox{\rm (\textit{cf.} \cite[Corollary 5.7]{dJ98} (i))}}\label{indbl} 
Let $\widetilde{N}^\dagger$ be a nonzero $\widetilde{\varphi}$-module over $\widetilde{\mathcal E}^\dagger$, 
$\lambda_1$ the maximal slope of $\widetilde{N}$, 
and $\widetilde{N} = \widetilde{N}_1 \oplus \widetilde{N}'$ such that 
$\widetilde{N}_1$ is the $\widetilde{\varphi}$-submodule of $\widetilde{N}$ exactly of slope $\lambda_1$ 
and that $\widetilde{N}'$ is the $\widetilde{\varphi}$-submodule of $\widetilde{N}$ whose slopes are strictly less than $\lambda_1$. 
If $\eta : \widetilde{N}^\dagger \rightarrow \widetilde{N}'$ is the natural $\widetilde{\mathcal E}^\dagger$-homomorphism defined by 
$\widetilde{N}^\dagger \subset \widetilde{N} \rightarrow \widetilde{N}/\widetilde{N}_1 = \widetilde{N}'$, 
then $\mathrm{Ker}(\eta)$ is a $\widetilde{\varphi}$-submodule of $\widetilde{N}^\dagger$ over $\widetilde{\mathcal E}^\dagger$ 
such that $\widetilde{\mathcal E}\otimes_{\widetilde{\mathcal E}^\dagger} \mathrm{Ker}(\eta) = \widetilde{N}_1$. 
\end{lemma}

\prf{It is sufficient to prove $\mathrm{Ker}(\eta) \ne 0$. 
Since $K' \otimes_K\mathrm{Ker}(\eta) = \mathrm{Ker}(\mathrm{id}_{K'}\otimes \eta)$ for a finite extension $K'$ of $K$ 
with an extension of the $q$-Frobenius $\sigma$, we may assume $|q|^\lambda \in |K|$ by Lemma \ref{ffix} (3) in Appendix \ref{Frob}. 
Hence we may assume that the maximal slope $\lambda$ of $\widetilde{N}$ is $0$. 
Suppose $m \in \widetilde{N} \setminus \{ 0 \}$ satisfies $\varphi_{\widetilde{N}}(m) = m$. 
Such an $m$ exists by Proposition \ref{unitdec} (1) since $\widetilde{N}^\dagger \ne 0$ and $\lambda = 0$. 
Then Lemma \ref{eequat} implies $m \in \widetilde{N}^\dagger$. Indeed, take a cyclic vector $e$ of the dual $(\widetilde{N}^\dagger)^\vee$ 
of $\widetilde{N}^\dagger$ and 
$B \in \mathrm{GL}_s(\widetilde{\mathcal E}^\dagger)$ is the representation matrix of Frobenius $\varphi_{(\widetilde{N}^\dagger)^\vee}$ 
with respect to the basis $e, \varphi_{(\widetilde{N}^\dagger)^\vee}(e), \cdots, \varphi_{(\widetilde{N}^\dagger)^\vee}^{s-1}(e)$. 
Then all entries of $B$ belongs to $\mathcal O_{\widetilde{\mathcal E}^\dagger}$ since 
all slopes of the dual $\widetilde{N}^\vee$ of $\widetilde{N}$ is $\geq 0$. If we use the dual basis of 
$e, \varphi_{(\widetilde{N}^\dagger)^\vee}(e), \cdots, \varphi_{(\widetilde{N}^\dagger)^\vee}^{s-1}(e)$ as a basis of $\widetilde{N}^\dagger$, 
then the representation matrix of $\varphi_{\widetilde{N}^\dagger}$ is ${}^tB^{-1}$ where ${}^tB^{-1}$ means the transposition 
of the matrix $B^{-1}$. Then, by the similar way of the proof of Proposition \ref{unitdec} (2), 
one can prove $m$ belongs to $\widetilde{N}^\dagger$ by Lemma \ref{eequat}. 
Since all slopes of $\widetilde{N}'$ is strictly less than $0$, we have $\eta(m) = 0$. 
}

\begin{theorem}{\mbox{\rm (\textit{cf.} \cite[Corollary 5.7]{dJ98})}}\label{oprk1} Let $N^\dagger$ be a $\varphi$-module 
over $\mathcal E^\dagger$, put $\widetilde{N}^\dagger =  \widetilde{\mathcal E}^\dagger \otimes_{{\mathcal E}^\dagger} N^\dagger$ 
to be a $\widetilde{\varphi}$-module over $\widetilde{\mathcal E}^\dagger$ with the opposite slope filtration 
$$
      0 = \widetilde{N}_0^\dagger \subsetneq \widetilde{N}_1^\dagger \subsetneq \cdots \subsetneq \widetilde{N}_r^\dagger = \widetilde{N}^\dagger, 
$$
and $\widetilde{\eta} : \widetilde{N}^\dagger \rightarrow \widetilde{\mathcal E}$ a nonzero injective $\widetilde{\mathcal E}^\dagger$-homomorphism 
such that $\widetilde{\eta} \circ \varphi_{\widetilde{N}^\dagger} = \widetilde{\varphi} \circ \widetilde{\varphi}^\ast(\widetilde{\eta})$. 
Then the slope of $\widetilde{N}^\dagger_1$ is $0$ and 
$\mathrm{dim}_{\widetilde{\mathcal E}^\dagger}N_1^\dagger = 1$. 
\end{theorem}

\prf{Since there exists a nonzero morphism 
$\widetilde{\mathcal E}\otimes_{\widetilde{\mathcal E}^\dagger} \widetilde{N}_1^\dagger \rightarrow \widetilde{\mathcal E}$ of $\varphi$-modules 
over $\widetilde{\mathcal E}$, the slope of $\widetilde{N}_1^\dagger$ should be $0$. Then 
$\widetilde{N}_1^\dagger$ has a basis $e_1, \cdots, e_{n_1}$ over $\widetilde{\mathcal E}^\dagger$ such that 
$\varphi_{\widetilde{N}^\dagger}(e_i) = e_i$ for $1 \leq i \leq n_1$ by Proposition \ref{unitdec} (2). Then 
$$
     \widetilde{\eta}(1\otimes e_i) \in \widetilde{\mathcal E}_{\varphi} = K_\sigma.
$$
for any $i$ by the commutativity of Frobenius and Lemma \ref{Efield} (2). Since $\widetilde{\eta}$ is injective, 
the equality $\mathrm{dim}_{\widetilde{\mathcal E}^\dagger}\widetilde{N}_1^\dagger = 1$ holds. 
}

\vspace*{3mm}

Now we assume the residue field $k$ of $K$ is a perfect field of characteristic $p$ 
and $\sigma$ is a $q$-Frobenius on $K$ without any extra condition. 

\begin{theorem}\label{rkr} Let $N^\dagger$ be a $\varphi$-$\nabla$-module over $\mathcal E^\dagger$, 
and $M$ a $\varphi$-$\nabla$-module over $\mathcal E$ such that $M$ has a unique slope. Let $\eta : N^\dagger \rightarrow M$ 
be an injective $\mathcal E^\dagger$-homomorphism which is compatible with Frobenius and connections. 
If $\eta(N^\dagger)$ generates $M$ as an $\mathcal E$-space, 
then the maximal slope of $N$ coincides with the slope of $M$ and the equality 
$\mathrm{dim}_{\mathcal E}N/N^1 = \mathrm{dim}_{\mathcal E}M$ holds. 
\end{theorem}

\prf{Since the maximal slope and its dimension are stable under the scalar extension, 
we may assume that $k$ is algebraically closed and the natural map 
$K_\sigma \otimes_{W(\mathbb F_q)}W(k) \rightarrow K$ is an isomorphism by Lemma \ref{ffix} (2) in Appendix \ref{Frob}. 
Since the category $\mbox{\rm \bf $\Phi$M}^\nabla_{\mathcal E^\dagger}$ is independent of the choice of 
$q$-Frobenius $\varphi$ with respect to $\sigma$ \cite[Theorem 3.4.10]{Ts98b}, we may assume that $\varphi(t) = t^q$. 
We may also assume $M$ is unit-root. Indeed, if $s$ is the slope of $M$, 
then one can find a finite extension $K'$ of $K$ 
such that the $q$-Frobenius $\sigma$ extends on $L$ and 
there exists an element $b \in K'_\sigma$ with $|b| = |q|^s$ by Lemma \ref{ffix} (3). 
Put $\widetilde{M} = \widetilde{\mathcal E} \otimes_{\mathcal E}M, 
\widetilde{N}^\dagger=\widetilde{\mathcal E}^\dagger \otimes_{\mathcal E^\dagger}N^\dagger$. 
The injectivity of the natural map $\widetilde{\mathcal E}^\dagger \otimes_{\mathcal E^\dagger} \mathcal E \rightarrow \widetilde{\mathcal E}$ by 
Lemma \ref{Efield} (3) implies 
the homomorphism $\widetilde{\eta} : \widetilde{N}^\dagger \rightarrow \widetilde{M}$ induced by $\eta : N^\dagger \rightarrow M$ 
is injective. Therefore, the assertion follows from Theorem \ref{oprk1} 
since $\widetilde{M}$ is isomorphic to $\widetilde{\mathcal E}^{\otimes \mathrm{dim}_{\mathcal E}M}$ as $\widetilde{\varphi}$-modules by 
Proposition \ref{unitdec} (1) and 
$\mathrm{dim}_{\mathcal E}N/N^1 = \mathrm{dim}_{\widetilde{\mathcal E}^\dagger}\widetilde{N}_1^\dagger$ 
by Proposition \ref{opfil}.}

\section{PBQ overconvergent $F$-isocrystals}\label{log}

We introduce the notion of the PBQ filtration for overconvergent $F$-isocrystals on a smooth curve. 
The notion has been already defined for local objects, $\varphi$-$\nabla$-modules, by B.Chiarellotto and the author in \cite{CT11} 
in order to give a necessary 
condition of $\varphi$-$\nabla$-modules over $\mathcal E$ which satisfies Dwork's conjecture on logarithmic growth. 
We will globalize the notion of PBQ modules in this section. 

Let us keep the notation $k, R, K, \sigma$ in section \ref{dJ} such that the residue field $k$ of $K$ is supposed perfect and 
$\sigma$ is a $q$-Frobenius on $K$. 

\subsection{Settings}\label{set} At first we fix our situation (see Appendix \ref{Frob} for a brief introduction of $\dagger$-algebras, their Frobenius endomorphisms 
and (over)convergenet $F$-isocrystals). 

Let us fix the notation as follows:
\begin{list}{}{}
\item[$C$ :] a smooth connected affine curve over $\mathrm{Spec}\, k$ such that $\overline{C}$ is a smooth completion of $C$ 
with an open immersion $j_C : C \rightarrow \overline{C}$;
\item[$\mathcal C$ :] a smooth affine lift $\mathrm{Spec}\, A_C$ of $C$ over $\mathrm{Spec}\, R$;
\item[$\widehat{\mathcal C}$ :] the $p$-adic formal completion of $\mathcal C$. 
\end{list}
Let us fix a projective smooth lift $\overline{\mathcal C}$ of $\overline{C}$ over $\mathrm{Spec}\, R$ \cite[III, Corollaire 7.4]{Gr}. 
Then one can regard $\widehat{\mathcal C}$ as an open formal subscheme of 
the $p$-adic formal completion $\widehat{\overline{\mathcal C}}$ 
over $\mathrm{Spf}\, R$ by \cite[Proposition 2.4.4 (i)]{vP86}. Let us put 
$$
    \begin{array}{llllll}
     A^\dagger_{C, K} &= A_C^\dagger\otimes_RK &= &\Gamma(]\overline{C}[_{\widehat{\overline{\mathcal C}}}, 
     j_C^\dagger{\mathcal O_{]\overline{C}[_{\widehat{\overline{\mathcal C}}}}}) 
     &\mbox{where $A_C^\dagger$ is the $p$-adically weak completion of $A_C$}  \\
          \widehat{A}_{C, K} &= \widehat{A}_C\otimes_RK&= &\Gamma(]C[_{\widehat{\overline{\mathcal C}}}, 
          j_C^\dagger{\mathcal O_{]\overline{C}[_{\widehat{\overline{\mathcal C}}}}})
     &\mbox{where $\widehat{A}_C$ is the $p$-adically formal completion of $A_C$}. 
     \end{array}
$$
Here $j_C^\dagger\mathcal O_{]\overline{C}[_{\widehat{\overline{\mathcal C}}}}$ is 
the sheaf of functions 
on the rigid analytic space $]\overline{C}[_{\widehat{\overline{\mathcal C}}} = \widehat{\overline{\mathcal C}}^{\mathrm{an}}$
overconvergent  along $\overline{C} \setminus C$. Note that $A_C^\dagger$ and $\widehat{A}_C$ are independent of the choice of 
the lift $\mathcal C$ of $C$ up to $R$-isomorphisms, 
and are Noetherian integral domains by \cite[Theorem]{Fu69}. 
Since $\mathcal C$ is smooth over $\mathrm{Spec}\, R$, 
there exist a $q$-Frobenius 
$$
     \varphi :  A_{C, K}^\dagger \rightarrow A_{C, K}^\dagger
$$
with $\varphi|_K = \sigma$ 
and a continuous derivation
$$
      d : A_{C, K}^\dagger \rightarrow A_{C, K}^\dagger\otimes_{A_C}\Omega_{A_C/R}^1
$$
\cite[Section 2.4]{vP86}. 
Any $q$-Frobenius $\varphi$ is faithfully flat since $C$ is smooth over $\mathrm{Spec}\, k$. 
Let us define a $p$-adic lift $E_\eta$ of the function field $k(C)$ of $C$ by 
$$
      E_\eta = \mbox{\rm the $p$-adic completion of the localization $(A_C)_{\mathbf m}$ and then tensoring $\otimes K$}
$$
where $\eta$ means the generic point of $C$. 
$E_\eta$ is an extension of complete valuation fields over $K$ having the same valuation group with $K$ such that 
the residue field of $E_\eta$ is the function field $k(C)$ of $C$.  
Then the Frobenius $\varphi$ and the derivation $d$ extend uniquely on $E_\eta$.  

For a closed point $\alpha \in \overline{C}$ with the canonical $i_{\alpha,  \overline{C}} : \alpha \rightarrow  \overline{C}$, 
let us put
\begin{list}{}{}
\item[$k_\alpha$ :] the function field of $\alpha$;
\item[$K_\alpha$ :] the finite unramified extension of $K$ such that the residue field is $k_\alpha$;
\item[$R_\alpha$ :] the integer ring of $K_\alpha$; 
\item[\hspace*{3mm} $x_\alpha \in \mathcal O_{\overline{\mathcal C}}$ :] 
a lift of local coordinate of $\overline{C}$ at $\alpha$.  
\end{list}
When $\alpha \in C$, we have two natural commutative diagrams  
$$
           \begin{array}{ccc}
    \widehat{A}_{C, K} &\rightarrow &K_\alpha[\hspace*{-0.2mm}[x_\alpha]\hspace*{-0.2mm}]_0 \\
    \cap & &\cap \\
     E_\eta &\rightarrow &\mathcal E_\alpha  \\
       \end{array}
              \hspace*{10mm} 
    \begin{array}{ccccc}
      A^\dagger_{C, K} &\rightarrow &A^\dagger_{C\setminus \{ \alpha\}, K} &\rightarrow &\mathcal E^\dagger_\alpha  \\
      \cap & &\cap& &\cap   \\
       \widehat{A}_{C, K} &\rightarrow &\widehat{A}_{C\setminus \{\alpha\}, K} &\rightarrow &\mathcal E_\alpha  \\
       \end{array}
$$
where $\mathcal E_\alpha, \mathcal E_\alpha^\dagger, K_\alpha[\hspace*{-0.2mm}[x_\alpha]\hspace*{-0.2mm}]_0$ 
are the rings with coordinate $x_\alpha$ over $K_\alpha$ 
which are introduced in Section \ref{dJnot}. Moreover, 
the $q$-Frobenius $\varphi$ and the derivation $d$ 
on $A^\dagger_{C, K}$ extend uniquely and compatibly on $K_\alpha$ and 
all rings in the diagrams above. 

Let us introduce partially weakly complete finitely generated (w.c.f.g.) algebras 
in order to study the relation of various $\dagger$-algebras and $p$-adically complete algeberas. 
Let $V$ be an open subscheme of $\overline{C}$ over $\mathrm{Spec}\, k$ 
such that $C \subset V \subset \overline{C}$ 
with the open immersion $j_{C, V} : C \rightarrow V$ and consider 
the $K$-algebra 
$$
     \widehat{A}^\dagger_{C, V, K} = \Gamma(]V[_{\widehat{\overline{\mathcal C}}}, 
      j_C^\dagger{\mathcal O_{]\overline{C}[_{\widehat{\overline{\mathcal C}}}}}). 
$$
If $V = \overline{C}$, then $\widehat{A}^\dagger_{C, V, K} = \widehat{A}^\dagger_{C, K}$. 
Suppose $V \subsetneq \overline{C}$ and take an affine smooth lift $\mathcal V = \mathrm{Spec}\, A_V$ 
of $V$ over $\mathrm{Spec}\, R$. 
If there exists an element $y \in A_V$ such that $\mathrm{Spec}\, A_V/(y, \mathbf m)A_V = V \setminus C$ as closed subsets of $V$, 
then 
$$
        \widehat{A}^\dagger_{C, V, K} \cong (\widehat{A}_V[z]^\dagger /(yz-1))\otimes_RK. 
$$

\begin{lemma}\label{flat} With the notation as above, let $\alpha$ be a closed point in $V \setminus C$. 
Consider the natural commutative diagram 
$$
\begin{array}{ccccc}
    \widehat{A}^\dagger_{C \cup \{ \alpha \}, V, K} &\rightarrow 
    &\widehat{A}^\dagger_{C, V, K} &\rightarrow &\widehat{A}^\dagger_{C, V \setminus \{ \alpha \}, K} \\
    \downarrow & &\downarrow & &\downarrow \\
    K_\alpha[\hspace*{-0.2mm}[x_\alpha]\hspace*{-0.2mm}]_0 &\rightarrow &\mathcal E_\alpha^\dagger &\rightarrow &\mathcal E_\alpha. 
       \end{array}
$$
\begin{enumerate}
\item The three vertical morphisms above are flat.
\item $\widehat{A}^\dagger_{C \cup \{ \alpha \}, V, K} \rightarrow \widehat{A}^\dagger_{C, V, K}$ is flat 
and $\widehat{A}^\dagger_{C, V, K} \rightarrow \widehat{A}^\dagger_{C, V \setminus \{ \alpha \}, K}$ is faithfully flat. 
\item The equality 
$\widehat{A}^\dagger_{C \cup \{ \alpha\}, V, K} = K_\alpha[\hspace*{-0.2mm}[x_\alpha]\hspace*{-0.2mm}]_0 \cap \widehat{A}^\dagger_{C, V, K}$ 
(resp. $\widehat{A}^\dagger_{C, V, K} = \mathcal E^\dagger_\alpha \cap \widehat{A}^\dagger_{C, V \setminus \{ \alpha \}, K}$) holds in 
$\mathcal E^\dagger_\alpha$ (resp. $\mathcal E_\alpha$). 
\end{enumerate}
\end{lemma}

\prf{We may assume that $V$ is affine and $\alpha \in V, C = V \setminus \{ \alpha\}$ by glueing after Tate's acyclic theorem 
and the fact that $]V[_{\widehat{\overline{\mathcal C}}}$ is quasi-separated and quasi-compact. 
In this case $\widehat{A}_{C \cup \{\alpha\}, V, K}^\dagger = \widehat{A}_V = \widehat{A}_{C \cup \{\alpha\}}$ and 
$\widehat{A}^\dagger_{C, V \setminus \{ \alpha \}, K} = \widehat{A}_{C, K}$. 
Moreover, we may assume that our local coordinate $x_\alpha$ at $\alpha$ belongs to $A_V$ 
and the closed subscheme defined by $x_\alpha = 0$ in $V$ consists of only $\alpha$ after shrinking $V$. 

(1) Let $I_\alpha$ be the maximal ideal of $A_V$ of the closed point $\alpha$. 
Since $R_\alpha[\hspace*{-0.2mm}[x_\alpha]\hspace*{-0.2mm}]$ is naturally isomorphic to 
the $I_\alpha$-adic completion of $\widehat{A}_V$, 
$R_\alpha[\hspace*{-0.2mm}[x_\alpha]\hspace*{-0.2mm}]$ is flat over $\widehat{A}_V$. 
Since $\widehat{A}^\dagger_{C, V, K}$ (resp. $\widehat{A}^\dagger_{C, V \setminus \{ \alpha \}, K}$) is an integral domain 
and $\mathcal E_\alpha^\dagger$ (resp. $\mathcal E_\alpha$) is a field, the rest is trivial. 

(2) Since $\widehat{A}^\dagger_{C \cup \{ \alpha \}, V, K} = \Gamma(]V[_{\widehat{\overline{\mathcal C}}}, 
j_{C \cup \{ \alpha\}}^\dagger\mathcal O_{]\overline{C}[_{\widehat{\overline{\mathcal C}}}}) \rightarrow 
\Gamma(]C[_{\widehat{\overline{\mathcal C}}}, 
j_{C \cup \{ \alpha \}}^\dagger\mathcal O_{]\overline{C}[_{\widehat{\overline{\mathcal C}}}}) = \widehat{A}_{C, K}$ is flat, 
we have only to prove $\widehat{A}^\dagger_{C, V, K} \rightarrow \widehat{A}_{C, K}$ 
is faithfully flat. The faithful flatness holds because the $\mathbf m$-adic topological ring $\widehat{A}^\dagger_{C, V}$ is a Zariski ring. 
Indeed, any element in $1 + \mathbf m\widehat{A}^\dagger_{C, V}$ is a unit in $\widehat{A}^\dagger_{C, V}$ so that 
any maximal ideal of $\widehat{A}^\dagger_{C, V}$ includes $\mathbf m\widehat{A}^\dagger_{C, V}$. 

(3) We may assume that $\alpha$ is $k$-rational by gluing of sections of rigid analytic spaces 
since $\widehat{A}_{C, K} \cap \widehat{A}^\dagger_{C_l, V_l, L} = \widehat{A}^\dagger_{C, V, K}$ 
by the equality $\widehat{A}^\dagger_{C_l, V_l, L} = \widehat{A}^\dagger_{C, V, K}\otimes_KL$ 
for any base change by a finite unramified extension $L$ of $K$ with the residual extension $l$ of $k$ 
and $C_l = C \times_{\mathrm{Spec}\, k}\mathrm{Spec}\, l, V_l = V \times_{\mathrm{Spec}\, k}\mathrm{Spec}\, l$. 
Then we have only to prove, if $V = C \cup \{ \alpha\}$, then  the equality 
$$
    \widehat{A}_{V, K} = K_\alpha[\hspace*{-0.2mm}[x_\alpha]\hspace*{-0.2mm}]_0 \cap \widehat{A}_{C, K}
$$
holds. Indeed, the equality above implies 
$\widehat{A}^\dagger_{C, V, K} = \mathcal E_\alpha^\dagger \cap \widehat{A}_{C, K}$ 
because, for $h = \sum_n c_nx^n_\alpha$ in $\mathcal E_\alpha$ (resp. $\mathcal E_\alpha^\dagger$), 
the minus part $h^{(-)} = \sum_{n<0} c_nx^n_\alpha$ belongs to $\widehat{A}_{C, K}$ (resp. $\widehat{A}^\dagger_{C, V, K}$) 
and $h - h^{(-)} \in \widehat{A}_{V, K}$. 
The equality follows from the integral version of the identification 
$$
      \widehat{A}_{V} = R[\hspace*{-0.2mm}[x_\alpha]\hspace*{-0.2mm}] \cap \widehat{A}_C
$$
which easily comes from the equality 
$\Gamma(V, \mathcal O_{\overline{C}})  = \{ f \in \Gamma(C, \mathcal O_{\overline{C}})\, |\, v_\alpha(f) \geq 0 \}$,  
where $v_\alpha$ is a valuation of $k(C)$ at $\alpha$, since all rings in the integral version are $p$-adically complete. 
}

\begin{lemma}\label{intersect} 
\begin{enumerate}
\item $\widehat{A}_{C, K} = \underset{\tiny \begin{array}{c} \alpha \in C \\ \mbox{\rm a closed point}\end{array}}{\cap} 
E_\eta \cap K_\alpha[\hspace*{-0.2mm}[x_\alpha]\hspace*{-0.2mm}]_0$. 
\item $A^\dagger_{C, K} = \underset{\tiny \begin{array}{c} \alpha \in \overline{C} \setminus C \\ \mbox{\rm a closed point}\end{array}}{\cap} 
 \mathcal E_\alpha^\dagger \cap \widehat{A}_{C, K} = 
\underset{\tiny \begin{array}{c} \alpha \in C \\ \mbox{\rm a closed point}\end{array}}{\cap}
        E_\eta\cap K_\alpha[\hspace*{-0.2mm}[x_\alpha]\hspace*{-0.2mm}]_0\, \, 
        \bigcap \underset{\tiny \begin{array}{c} \alpha \in \overline{C} \setminus C \\ \mbox{\rm a closed point}\end{array}}{\cap} 
        E_\eta \cap \mathcal E^\dagger_\alpha$. 
\end{enumerate}
\end{lemma}

\prf{(1) Since all the rings are $p$-adically complete, it is sufficient to prove 
$A_C/\mathbf mA_C = \cap_{\alpha \in C} k(C) \cap k_\alpha[\hspace*{-0.2mm}[x_\alpha]\hspace*{-0.2mm}]$. 
Here $k(C)$ is the field of functions of $C$. 
It holds because $A_C/\mathbf mA_C$ is a normal integral domain. 

(2) For any closed point $\alpha \in \overline{C} \setminus C$, there exists an affine open subscheme $V$ of $\overline{C}$ 
such that $\alpha \in V$ and $\alpha$ is defined by a single element of $\Gamma(V, \mathcal O_V)$ in $V$. 
By using the gluing argument, the first equality follows from Lemma \ref{flat} (3). 
The second equality follows from (1). 
}

\subsection{Frobenius slope filtration}\label{notat}
Let $\mathcal M^\dagger$ be an overconvergent $F$-isocrystal on $C/K$, and 
put 
$$
   \begin{array}{l}
   M^\dagger = \Gamma(]\overline{C}[_{\widehat{\overline{\mathcal C}}}, \mathcal M^\dagger), \hspace*{5mm} 
   M=\widehat{A}_{C, K}\otimes_{A^\dagger_{C, K}}M^\dagger = \Gamma(]C[_{\widehat{\overline{\mathcal C}}}, \mathcal M) \\
     M_\eta = E_\eta \otimes_{\widehat{A}_{C, K}} M, \\
     M^\dagger_\alpha = \mathcal E^\dagger_\alpha\otimes_{A^\dagger_{C, K}}M^\dagger, \hspace*{5mm} 
     M_\alpha = \mathcal E \otimes_{\widehat{A}_{C, K}}M 
     \hspace*{3mm} \mbox{\rm for a closed point $\alpha$ of $\overline{C} \setminus C$.}
     \end{array} 
$$
$M^\dagger$ is a projective $A^\dagger_{C, K}$-module of finite rank with a $K$-connection 
$\nabla_{M^\dagger} : M^\dagger \rightarrow M^\dagger \otimes_{A_C}\Omega^1_{A_C/R}$ 
and an isomorphism  $\varphi_{M^\dagger} : \varphi^\ast M^\dagger \rightarrow M^\dagger$ which is called Frobenius 
such that $(\varphi_{M^\dagger}\otimes\varphi) \circ \varphi^\ast\nabla_{M^\dagger} = \nabla_{M^\dagger} \circ \varphi_{M^\dagger}$. 
The similar hold for the $\widehat{A}_{C, K}$-module $M$ and the $E_\eta$-space $M_\eta$. 
There are various natural isomorphisms 
$$
    \begin{array}{l}
     \widehat{A}_{C, K}\otimes_{A^\dagger_{C, K}}M^\dagger \cong M, \hspace*{3mm} E_\eta\otimes_{\widehat{A}_{C, K}}M \cong M_\eta, \\
     \mathcal E^\dagger_\alpha\otimes_{A^\dagger_{C, K}}M^\dagger \cong M_\alpha^\dagger, \hspace*{3mm} 
     \mathcal E_\alpha\otimes_{\widehat{A}_{C, K}}M \cong  \mathcal E_\alpha\otimes_{E_\eta}M_\eta\cong M_\alpha, \\
     \mathcal E_\alpha\otimes_{\mathcal E_\alpha^\dagger}M_\alpha^\dagger \cong M_\alpha
     \end{array}
$$
of $\varphi$-$\nabla$-modules. 

\begin{proposition}\label{frfil} \mbox{\rm (See \cite[Remark 1.7.8]{Ke08}, \cite[Theorem 2.4]{CT09}, and Proposition \ref{slfilg}.)} 
There exists a slope filtration $\{ M_\eta^i \}$ of $M_\eta$ with respect to Frobenius $\varphi_{M_\eta}$ 
as $\varphi$-$\nabla$-modules over $E_\eta$ such that 
$\{ \mathcal E_\alpha \otimes_{E_\eta} M_\eta^i\}$ coincides with the slope filtration $\{ M_\alpha^i \}$ 
of $M_\alpha$ with respect to Frobenius $\varphi_{M_\alpha}$. 
\end{proposition}

\subsection{Generic disc}\label{gendisc}

Let $B$ be one of $E_\eta$ and $\mathcal E_\alpha$ for a closed point $\alpha$ 
in $\overline{C}$, and put $\ast = \eta\, \, \mathrm{or}\, \, \alpha$ respectively. 
We recall the notion of generic open unit disc (see \cite[2.5, 4.1, 4.6]{Ch83} and \cite[0.4]{CT09}). 
Let $B^\tau$ be the complete discrete valuation ring which is isomorphic to $B$, but 
the coordinate $x_\ast$ is replaced by $t_\ast$ in $B^\tau$ for $\ast = \eta\,\, \mathrm{or}\, \, \alpha$. 
In the case where $\ast = \eta$ we fix a generically etale morphism $\mathcal C \rightarrow \mathbb A^1_R$ 
over $\mathrm{Spec}\, R$ and take the coordinate $x_\eta = x$ which is the standard one of the affine line $\mathbb A_R^1$. 
Then $B^\tau$ is also an extension of $K$ as a complete discrete valuation field with the same valuation groups. 
We introduce two $B^\tau$-algebras $B^\tau[\hspace*{-0.2mm}[X-t_\ast]\hspace*{-0.2mm}]_0$ 
and $\mathcal A_{B^\tau}(t_\ast, 1)$, called the ring of bounded functions on the unit open generic disc and 
a ring of analytic functions on the unit open generic disc respectively as follows: 
$$
\begin{array}{cll}
B^\tau[\hspace*{-0.2mm}[X-t_\ast]\hspace*{-0.2mm}]_0 &= &\left\{ \left. \sum_n a_n(X-t_\ast)^n \in B^\tau[\hspace*{-0.2mm}[X-t_\ast]\hspace*{-0.2mm}]\, \right|\, 
\mathrm{sup}_n|a_n| < \infty \right\}, \\
  \mathcal A _{B^\tau}(t_\ast, 1) &= &\left\{ \sum_n a_n(X-t_\ast)^n \in B^\tau[\hspace*{-0.2mm}[X-t_\ast]\hspace*{-0.2mm}]\, \left|\, 
\begin{array}{l} \sum_n a_n(X-t_\ast)^n\, \, \mbox{\rm is convergent} \\ \mbox{\rm on $|X-t_\ast| < 1$.} 
\end{array}\right.\right\}. 
\end{array}
$$
For an element $f \in B$, we put 
$$
     f^\tau = \sum_{n=0}^\infty\, \frac{d^nf}{dx^n}(t_\ast)\frac{(X-t_\ast)^n}{n!}
$$
where $g(t_\ast)$ means the evaluation of $g$ at $x_\ast = t_\ast$ for $g \in B$. 
We define a $q$-Frobenius $\varphi^\tau$ on $B^\tau[\hspace*{-0.2mm}[X-t_\ast]\hspace*{-0.2mm}]_0 \subset \mathcal A _{B^\tau}$ by 
$$
  \begin{array}{l}
  \varphi^\tau|_{B^\tau} = \varphi\, \, \mbox{\rm under the isomorphism $B^\tau \rightarrow B\, \, (t_\ast \mapsto x_\ast)$} \\
  \varphi^\tau(X-t_\ast) = \varphi(x_\ast)^\tau - \varphi^\tau(t_\ast). 
  \end{array}
$$
Since $\varphi(x_\ast)^\tau - \varphi^\tau(t_\ast) \in (X-t_\ast)B^\tau[\hspace*{-0.2mm}[X-t_\ast]\hspace*{-0.2mm}]_0$, $\varphi^\tau$ is well-defined. 

\begin{lemma} The application $f \mapsto f^\tau$ is a $K$-algebra homomorphism such that 
\begin{list}{}{}
\item[\mbox{\rm (i)}] $(\frac{df}{dx_\ast})^\tau = \frac{df^\tau}{dX}$; 
\item[\mbox{\rm (ii)}] $\varphi(f)^\tau = \varphi^\tau(f^\tau)$. 
\end{list}
\end{lemma}

\vspace*{3mm}

Let $M$ be a $\nabla$-module (resp. a $\varphi$-$\nabla$-module) $M$ over $B$. 
We define a $\nabla$-module (resp. a $\varphi$-$\nabla$-module) $M^\tau$ over $B^\tau[\hspace*{-0.2mm}[X-t_\ast]\hspace*{-0.2mm}]_0$ 
associated to $M$ by 
$$
    \begin{array}{ll} 
     &M^\tau = B^\tau[\hspace*{-0.2mm}[X-t_\ast]\hspace*{-0.2mm}]_0 \otimes_BM \\
     &\nabla_{M^\tau}(a \otimes m) = a\otimes\nabla_M(m) + \frac{da}{dX}\otimes mdX \\
     \mbox{\rm (resp.} &\varphi_{M^\tau}(a\otimes m) = \varphi^\tau(a)\otimes\varphi_M(m)\, ). 
     \end{array}
$$
If the matrix representation of the connection $\nabla_M$ of $M$ of arbitrary order is given by 
$$
   \nabla_M\left(\frac{d}{dx_\ast}\right)^n(e_1, \cdots, e_s) = (e_1, \cdots, e_s)C_n\, \hspace*{5mm} C_n \in M_s(B)
$$
for any nonnegative integer $n$, where $C_0$ is a unit matrix, then the connection of $M^\tau$ is given by 
$$
    \nabla_{M^\tau}(1\otimes e_1, \cdots, 1\otimes e_s) = (1\otimes e_1, \cdots, 1\otimes e_s)C_1^\tau dX. 
$$
Hence the solution matrix of $M^\tau$, which is called the generic solution matrix of $M$ at the generic point $t_\ast$, is 
$$
      Y  = \sum_{n = 0}^\infty C_n(t_\ast)\frac{(X-t_\ast)^n}{n!}.
$$

\begin{definition}\label{solv} Let $M$ be a $\nabla$-module over $B$. 
\begin{enumerate}
\item $M$ is said to be solvable if, for any $0 < \eta < 1$, 
$$
    \left|\frac{1}{n!}C_n\right|\eta^n \rightarrow 0 \hspace*{3mm} \mbox{\rm as}\, \, n \rightarrow \infty. 
$$
Equivalently, $M$ is solvable if all entries of the solution matrix $Y$ above belong to $\mathcal A_{B^\tau}(t_\ast, 1)$. 
\item $M$ is bounded if $M$ is solvable and satisfies 
$$
    \sup_n\left|\frac{1}{n!}C_n\right| < \infty. 
$$
Equivalently, $M$ is bounded if all entries of the solution matrix $Y$ above belong to $B^\tau[\hspace*{-0.2mm}[X-t_\ast]\hspace*{-0.2mm}]_0$. 
\end{enumerate}
\end{definition}

\begin{proposition}\label{bmod}
The notion of solvability (resp. boundedness) of $\nabla$-modules over $B$ 
does not depends on the choice of basis $e_1, \cdots, e_s$ of $M$ and the choice of coordinate $x_\ast$. 
\end{proposition}

\prf{The independence of the choice of basis follows from the fact that the map $f \mapsto f^\tau$ 
preserves units. The independence of coordinates follows from 
the fact that the formal lift $\widehat{A}_C$ of $A_C$ is independent of the choices 
up to continuous isomorphisms. Hence, for another choice of coordinate $x_\ast'$, one has a continuous isomorphism   
$B_{x_\ast'} \cong B_{x_\ast}$ so that 
$B_{x_\ast'}^\tau[\hspace*{-0.2mm}[X'-t_\ast']\hspace*{-0.2mm}]_0 \cong B_{x_\ast}^\tau[\hspace*{-0.2mm}[X-t_\ast]\hspace*{-0.2mm}]_0$ 
and $\mathcal A_{B_{x_\ast'}^\tau}(t_\ast', 1) \cong \mathcal A_{B_{x_\ast}^\tau}(t_\ast, 1)$. 
Here $B_{x_\ast'}$ is the differential field with coordinate $x_\ast'$.}

\begin{proposition} The category of solvable (resp. bounded) $\nabla$-modules over $B$ 
is Abelian, and it is closed under tensor products and duals. 
\end{proposition}

\prf{One can easily see the existence of duals follows from the solvable $\nabla$-module over $B$ is bounded. 
The boundedness follows from \cite[Lemma 1.7]{CT09}.
}

\vspace*{2mm}

Because one can choose the coordinate $x_\alpha$ at $\alpha$ as a coordinate of $E_\eta$, 
Proposition \ref{bmod} implies the proposition below. 

\begin{proposition} Let $M_\eta$ be a $\nabla$-module over $E_\eta$. 
$M_\eta$ is solvable (resp. bounded) if and only if so is $\mathcal E_\alpha\otimes_{E_\eta}M_\eta$ 
for a closed point (all closed points) $\alpha$ of $\overline{C}$. 
\end{proposition}

\subsection{Bounded $\varphi$-$\nabla$ modules}\label{boumod}

In this subsection we study properties of bounded $\varphi$-$\nabla$-modules 
over $B = E_\eta$ or $\mathcal E_\alpha$. At first we recall a well-known fact (see \cite[Theorem 6.6]{CT09} for example). 

\begin{proposition}\label{unet} Let $M$ be a $\varphi$-$\nabla$-module over $B$. 
Then the following hold. 
\begin{enumerate}
\item $M$ is solvable. 
\item If $M$ is unit-root, then $M$ is bounded. 
\end{enumerate}
\end{proposition}

\vspace*{1mm}

The following theorem is a characterization of bounded $\varphi$-$\nabla$-modules. Chiarellotto and the author proved 
it in the local case, i.e., $B = \mathcal E_\alpha$ in \cite[Theorem 4.1]{CT11}. 
Here we prove the assertion in the case where $B = E_\eta$ by using the local result. 

\begin{theorem}\label{split} Let $M$ be a $\varphi$-$\nabla$-module over $B$ with the slope filtration $\{ M^i\}$. 
Then $M$ is bounded if and only if $M \cong \oplus_i M^i/M^{i+1}$ as $\varphi$-$\nabla$-modules over $B$. 
\end{theorem}

Let us define the $p$-adic completion $\widehat{E}_\eta^{\mathrm{perf}}$ of the perfection of the residue field of $E_\eta$ by 
$$
    \widehat{E}_\eta^{\mathrm{perf}} = \mbox{\rm the $p$-adic completion of}\, \, \underset{\rightarrow}{\mathrm{lim}}\left(
    E_\eta \overset{\varphi}{\rightarrow} E_\eta \overset{\varphi}{\rightarrow} \cdots\right). 
$$
One can regard $\widehat{E}_\eta^{\mathrm{perf}}$ as a subfield of $\widetilde{\mathcal E}_\alpha$ in Section \ref{genser} by 
the natural embedding of direct systems 
$$
   \left(E_\eta \overset{\varphi}{\rightarrow} E_\eta \overset{\varphi}{\rightarrow} \cdots\right) 
   \rightarrow \left(\widetilde{\mathcal E}_\alpha \overset{\varphi}{\rightarrow} \widetilde{\mathcal E}_\alpha \overset{\varphi}{\rightarrow} \cdots\right). 
$$

\begin{lemma}\label{pperf} With the notation as above, $E_\eta = 
\widehat{E}_\eta^{\mathrm{perf}} \cap \mathcal E_\alpha$ in $\widetilde{\mathcal E}_\alpha$ for any closed point $\alpha \in \overline{C}$. 
\end{lemma}

\prf{Denote the perfection of $(-)$ by $(-)^{\mathrm{perf}}$, 
and the completion of $k(C)$ along $\alpha$ by $k(C)_\alpha$. 
Since all the items are $p$-adically complete discrete valuation fields with the same valuation group, the assertion follows 
from the fact 
$$
     k(C) = k(C)^{\mathrm{perf}} \cap k(C)_\alpha
$$
where the intersection is taken in $(k(C)_\alpha)^{\mathrm{perf}}$. Indeed, let $D : k(C) \rightarrow k(C)$ be the derivation induced 
by the standard derivation of $k(\mathbb P^1_k)$ and the finite separable morphism $C \rightarrow \mathbb P^1_k$. 
Then the $p$-power subfield $k(C)^p$ of $k(C)$ is the kernel of $D$ since $[k(C):k(C)^p] = p$. 
On the other hand, $D$ extends on $k(C)_\alpha$ and, if $u$ belongs to $(k(C)_\alpha)^p$, then $D(u) = 0$. Hence 
$$
           k(C)^p = k(C) \cap (k(C)_\alpha)^p. 
$$
The desired equality easily follows.}

\vspace*{2mm}

\noindent
{\sc Proof of Theorem \ref{split}.} Let $A = \left(\begin{array}{ccc} A_{11} &\cdots &A_{1r} \\ &\ddots &\vdots \\
0 & &A_{rr} \end{array}\right)$ be the representation matrix of Frobenius $\varphi_{M_\eta}$ with respect to the slope filtration of $M_\eta$. 
Since $\widehat{E}_\eta^{\mathrm{perf}}$ 
is $p$-adically complete, one has a unique solution $X$ which is a upper triangle matrix with entries in $\widehat{E}_\eta^{\mathrm{perf}}$ satisfying 
$$
       A\varphi(X) = X\left(\begin{array}{ccc} A_{11} & &0 \\ &\ddots & \\
0 & &A_{rr} \end{array}\right) \hspace*{5mm} X = \left(\begin{array}{ccc} E & &\ast \\ &\ddots & \\
0 & &E \end{array}\right) 
$$
because all the slopes of $A_i$'s are different. 
Here $E$ denotes the unit matrix of certain size. On the other hand we know $X \in M_n(\widehat{E}_\eta^{\mathrm{perf}}) 
\subset M_n(\widetilde{\mathcal E}_\alpha)$ belongs to $M_n(\mathcal E_\alpha)$ by \cite[Theorem 4.1]{CT11}. 
Hence $X$ belongs to $M_n(E_\eta)$ by Lemma \ref{pperf}. 
\hspace*{\fill} $\Box$

\subsection{Maximally bounded quotients} The following theorem was studied using $p$-adic functional 
analysis in \cite[\S 4.3]{Ch83}. In this paper we use the assertion only in the case with Frobenius structures.  

\begin{theorem}\label{bquot} \cite[Theoreme 4.3,5]{Ch83} Let $M$ be a solvable $\nabla$-module over $B$. 
Then there exists a unique minimal $\nabla$-submodule $M^b$ of $M$ over $B$ 
such that $M/M^b$ is nonzero and bounded. Here ``minimal" means, if $N$ 
is a $\nabla$-submodule of $M$ over $B$ such that $M/N$ is bounded, then $M^b \subset N$. Moreover, 
if $M$ is a $\varphi$-$\nabla$-module over $B$, then $M^b$ is a 
$\varphi$-$\nabla$-submodule of $M$ over $B$. 
\end{theorem}

\prf{Let $I$ be a set of $\nabla$-submodules over $B$ of $M$ such that, for any $N \in I$, 
$M/N$ is bounded. Our claim in the first part is 
(i) $\cap_{\mathcal N \in I}\mathcal N \in I$ and (ii) $I \supsetneq \{ M \}$. 
Suppose $N_1, N_2 \in I$. Since the natural morphism $M/(N_1 \cap N_2) \rightarrow M/N_1 \oplus M/N_2$ 
is injective, $M/(N_1 \cap N_2)$ is bounded so that $N_1 \cap N_2 \in I$. 
Hence the finite dimensionality of $M$ implies the claim (i). 
In the case where $M$ is a $\varphi$-$\nabla$-module over $B$ the claim (ii) follows from 
Propositions \ref{frfil} and \ref{unet}. In general case it was proved in \cite[Proposition 4.3.4]{Ch83}. 

If $M$ is $\varphi$-$\nabla$-module over $B$, then $M/\varphi^\ast M^b \cong \varphi^\ast(M/M^b)$ 
is also bounded, hence $\varphi^\ast M^b \in I$. Comparing with the dimensions, we have 
$\varphi^\ast M^b = M^b$. Hence, $M^b$ is a $\varphi$-$\nabla$-submodule of $M$. 
(See also \cite[Proposition 4.8]{CT09}.)
}

\begin{definition}\label{loggr} For a solvable $\nabla$-module $M$ over $B$, 
the unique nontrivial quotient $M/M^b$ is called the maximally bounded quotient of $M$. 
\end{definition}

\vspace*{2mm}

The following theorem is a bounded version of \cite[Proposition 1.10]{CT09} 
which implies the stability of logarithmic growth filtrations by scalar extensions. 

\begin{proposition}\label{bounf} Let $M_\eta$ be a solvable $\nabla$-module over $E_\eta$, $\alpha$ a closed point in $\overline{C}$ 
and choose the coordinates $x$ of $\widehat{A}_C$ and $x_\alpha$ of $\mathcal E_\alpha$ such that $x = x_\alpha$. 
If one defines $E_\eta^\tau$-spaces 
$\mathrm{Sol}(M_\eta^\tau)$ and $\mathrm{Sol}(M_\eta^\tau)$ 
by 
$$
    \begin{array}{lll}
       \mathrm{Sol}(M_\eta^\tau) &= &\left\{ f : M_\eta^\tau \rightarrow \mathcal A_{E_\eta^\tau}(t_\eta, 1)\, \left|\, \begin{array}{l} 
       \mbox{\rm $f$ is $E_\eta^\tau[\hspace*{-0.2mm}[X-t_\ast]\hspace*{-0.2mm}]_0$-linear such that} \\
       f(\nabla(\frac{d}{dX})(m)) = \nabla(\frac{d}{dX})(f(m))\, \mbox{\rm for any $m \in M^\tau_\eta$.}
       \end{array} \right.\right\} \\
              \mathrm{Sol}_0(M_\eta^\tau) &= &\left\{ f : M_\eta^\tau \rightarrow E_\eta^\tau[\hspace*{-0.2mm}[X-t_\eta]\hspace*{-0.2mm}]_0\, \left|\, \begin{array}{l} 
       \mbox{\rm $f$ is $E_\eta^\tau[\hspace*{-0.2mm}[X-t_\ast]\hspace*{-0.2mm}]_0$-linear such that} \\
       f(\nabla(\frac{d}{dX})(m)) = \nabla(\frac{d}{dX})(f(m))\, \mbox{\rm for any $m \in M^\tau_\eta$.}
       \end{array} \right.\right\}  
       \end{array}
$$
and defines $\mathcal E_\alpha^\tau$-spaces $\mathrm{Sol}(M_\alpha^\tau)$ and $\mathrm{Sol}_0(M_\alpha^\tau)$ 
for the solvable $\nabla$-module $M_\alpha = \mathcal E_\alpha\otimes_{E_\eta}M_\eta$ 
over $\mathcal E_\alpha$ similarly, 
then we have an isomorphism 
$$
      \mathrm{Sol}_0(M_\alpha^\tau) \cong \mathcal E_\alpha^\tau\otimes_{E_\eta^\tau}\mathrm{Sol}_0(M_\eta^\tau). 
$$
under the natural isomorphism $\mathrm{Sol}(M_\alpha^\tau) \cong \mathcal E_\alpha^\tau\otimes_{E_\eta^\tau}\mathrm{Sol}(M_\eta^\tau)$. 
\end{proposition}

\vspace*{2mm}

Note that, if the matrix representation of $M_\eta$ is as in Definition \ref{solv}, then there are natural $E_\ast^\tau$-isomorphisms 
$$
    \begin{array}{ccc}
       \mathrm{Sol}(M_\eta^\tau) &\cong &\{ y \in \mathcal A _{E_\eta^\tau}^s(t_\eta, 1)\, |\, \frac{dy}{dX} = yC^\tau \} \\
       \cup & &\cup \\
              \mathrm{Sol}_0(M_\eta^\tau) &= &\{ y \in E_\eta^\tau[\hspace*{-0.2mm}[X-t_\eta]\hspace*{-0.2mm}]_0^s\, |\, \frac{dy}{dX} = yC^\tau \} 
       \end{array}
$$
and the same for $M_\alpha$. 

\vspace*{3mm}

\noindent
{\sc Proof of Proposition \ref{bounf}.} Since the natural map 
$\mathcal E_\alpha^\tau\otimes_{E_\eta^\tau}\mathrm{Sol}_0(M_\eta^\tau) 
\rightarrow \mathrm{Sol}_0((\mathcal E_\alpha\otimes_{E_\eta}M_\eta)^\tau)$ 
is an injection of $\mathcal E_\alpha^\tau$-spaces, it is sufficient to prove that, if any linear combination of $f_1, \cdots, f_l \in \mathrm{Sol}(M_\eta^\tau)$ 
over $E_\eta^\tau$ is not contained in $f_1, \cdots, f_l \in \mathrm{Sol}_0(M_\eta^\tau)$, then 
any linear combination of $f_1, \cdots, f_l$ over $\mathcal E_\alpha^\tau$ is not contained in 
$\mathrm{Sol}_0((\mathcal E_\alpha\otimes_{E_\eta}M_\eta)^\tau)$. Consider a linear sum $c_1f_1 + \cdots +c_lf_l$ 
for some $c_1, \cdots, c_l \in \mathcal E_\alpha^\tau$. We will show that the linear sum 
does not belong to $\mathrm{Sol}_0((\mathcal E_\alpha\otimes_{E_\eta}M_\eta)^\tau)$. 
Let $V$ be an $E_\eta^\tau$-subspace of $\mathcal E_\alpha^\tau$ 
generated by $c_1, \cdots, c_l$, and $d_1, \cdots, d_m$ a basis of $V$ over $E_\eta^\tau$. Then 
$$
   c_1f_1 + \cdots +c_lf_l = d_1g_1 + \cdots + d_mg_m
$$
for some $g_1, \cdots, g_m \in \mathrm{Sol}(M_\eta^\tau)$. By our hypothesis either $g_i$ is not contained in $\mathrm{Sol}_0(M_\eta^\tau)$ 
or $g_i = 0$ for all $i$. Since $V$ is a finite dimensional topological vector space over the complete topological field $E_\eta^\tau$, 
the topology of $V$ induced by the norm of $\mathcal E_\alpha^\tau$ coincides with the product topology induced by 
the isomorphism $V \cong E_\eta^\tau d_1 \oplus \cdots \oplus E_\eta^\tau d_m$. 
If $c_1f_1 + \cdots +c_lf_l $ is bounded, then the sequence $g_{1, n}d_1 + \cdots + g_{m, n}d_m$ in $n$ is bounded in $V$ 
where $g_i = \sum_n\, g_{i, n}(X-t_\eta)^n$. This holds only in the case where $g_i =0$ for all $i$. 
Hence we complete the proof. 
\hspace*{\fill} $\Box$

\begin{corollary}\label{boundfil} With the notation as in Proposition \ref{bounf}, there is a natural isomorphism 
$$
      (\mathcal E_\alpha\otimes_{E_\eta}M_\eta)^b \cong \mathcal E_\alpha\otimes_{E_\eta}M_\eta^b
$$
as $\varphi$-$\nabla$-modules over $\mathcal E_\alpha$. 
\end{corollary}

\subsection{Generic PBQ filtrations} 

\begin{definition}\label{genPBQ} A $\varphi$-$\nabla$-module over $B$ is said to be PBQ (pure of bounded quotient) if the maximally bounded quotient 
$M/M^b$ has a unique Frobenius slope. Note that the slope of $M/M^b$ is the maximal slope of $M$ by Proposition \ref{unet} (2). 
\end{definition}

\begin{lemma}\label{subqut} 
\begin{enumerate}
\item If $M_1$ and $M_2$ are PBQ $\varphi$-$\nabla$-modules over $B$ of a same maximal slope, 
then so is the direct sum $M_1\oplus M_2$. 
\item Let $\theta : M \rightarrow N$ be 
a surjection of $\varphi$-$\nabla$-module over $B$. 
If $M$ is PBQ and $N \ne 0$, then $N$ is also PBQ of the same maximal slope with $M$. 
\end{enumerate}
\end{lemma}

\prf{The assertions easily follow from the definition.
}

\begin{proposition}\label{tenunit} Let $M$ be a $\varphi$-$\nabla$-module over $B$, 
and $N$ a unit-root $\varphi$-$\nabla$-module over $B$. 
Then $M$ is PBQ if and only if $M \otimes_B N$ is PBQ. 
\end{proposition}

\prf{On the generic disc the unit-root object is trivial as differential modules, that is, there exists 
an $F$-space $N^{\tau, 0}$ over $B^\tau$ such that 
$N^\tau \cong N^{\tau, 0}\otimes_{B^\tau}B^\tau[\hspace*{-0.2mm}[X-t_\eta]\hspace*{-0.2mm}]_0$. Since $(M\otimes_B N)^\tau$ is a 
direct sum of $\mathrm{dim}_BN$ copies of $M^\tau$ as differential modules, 
the bounded solutions of $(M\otimes_B N)^\tau$ 
are $\mathrm{dim}_BN$ copies of those of $M^\tau$. Hence 
$(M\otimes_B N)/(M\otimes_B N)^b 
\cong (M/M^b) \otimes_BN$. Therefore, the assertion holds. 
}

\begin{theorem}\label{genPBQfil} \mbox{\rm (\cite[Corollary 5.5]{CT11} in the case of $\mathcal E_\alpha$.)} 
Let $M$ be a $\varphi$-$\nabla$-module over $B$. 
\begin{enumerate}
\item There exists a unique increasing filtration 
$$
      0 = P_0 \subsetneq P_1 \subsetneq P_2 \subsetneq \cdots \subsetneq P_r = M
$$
as $\varphi$-$\nabla$-modules over $B$ such that 
\begin{list}{}{}
\item[\mbox{\rm (i)}] $P_i/P_{i-1}$ is PBQ;
\item[\mbox{\rm (ii)}] if $\lambda_i$ is the maximal slope of $P_i$, then $\lambda_1 > \lambda_2 > \cdots >\lambda_r$. 
\end{list}
Moreover, if an increasing filtration of $M$ satisfies the conditions (i), (ii), 
then it coincides with the filtration $\{ P_i\}$. 
The filtration $\{ P_i\}$ is called the PBQ filtration of $M$. 
\item There is an isomorphism $M/M^b \cong \oplus_{i=1}^r\, P_i/P_i^1$ as $\varphi$-$\nabla$-modules over $B$. 
\end{enumerate}
\end{theorem}

\prf{(1) We prove the assertion by induction on the dimension of $M$. If $M/M^b$ has only one slope, then $M$ is PBQ. 
If $M/M^b$ has several Frobenius slopes, then one obtain the maximal $\varphi$-$\nabla$-submodule $M'$ such that 
the image of $M' \rightarrow M/M^b$ coincides with the direct summand of maximal slope of $M/M^b$ by Proposition \ref{split}. Applying 
the induction hypothesis to $M'$ we has a nontrivial submodule $P_1$ of $M'$ which is PBQ and the image of 
$P_1 \rightarrow M/M^b$ coincides with  the direct summand of maximal slope of $M/M^b$. 
Applying the induction hypothesis to the quotient $M/P_1$ again, we has a desired filtration of $M$. 

Now we prove the uniqueness. Let $P_1$ and $P_1'$ be first steps of two filtrations satisfying the conditions (i) and (ii). 
Consider $Q = P_1 + P_1'$ in $M$ which is the image of $P_1 \oplus P_1' \rightarrow M$. If $Q \ne P_1$, then the nontrivial surjection $P_1' \rightarrow Q/P_1$ 
contradicts to the hypothesis that $P_1'$ is PBQ with the maximal slope $\lambda_1$ by Lemma \ref{subqut} (2) 
since there is an isomorphism $M/M^1 \cong P/P^1$ of the maximal slope quotients. Hence $Q = P_1 = P_1'$. 

(2) We prove the assertion by induction on the length $r$ of the PBQ filtration $\{ P_i\}$ 
of $M$. If $r=1$, then the assertion follows from 
the definition of PBQ. Suppose the assertion holds for $r-1$ and put $M' = M/P_1$. 
Then there exists an exact sequence 
$$
      P_1/P_1^b \rightarrow M/M^b \rightarrow M'/(M')^b \rightarrow 0 
$$
of $\varphi$-$\nabla$-modules over $B$. The left arrow is injective 
since the composite 
$\mathrm{Ker}(M/M^b \rightarrow M'/(M')^b) \rightarrow M/M^1 \overset{\cong}{\leftarrow} P_1/P_1^1$ is surjective 
by the maximality of the maximal slope $\lambda_1$ of $P_1$. Since $P_1$ is PBQ with the maximal slope $\lambda_1$, 
the assertion $M/M^b \cong \oplus_{i=1}^r\, P_i/P_i^1$ 
follows from the induction hypothesis and the maximality $\lambda_1$ of slopes of $M$. 
}

\begin{definition}\label{maxpBQ} The PBQ $\varphi$-$\nabla$-module $P_1$ of the first step of the PBQ filtration 
of $M$ in Theorem \ref{genPBQfil} above is called the maximally PBQ submodule of $M$. 
\end{definition}

The same notion will be defined for $\varphi$-$\nabla$-modules over $\mathcal E^\dagger_\alpha, K_\alpha[\hspace*{-0.2mm}[x_\alpha]\hspace*{-0.2mm}]_0$ in Section \ref{pbqdag} 
and for overconvergent $F$-isocrystals on $C/K$ in Section \ref{globalPBQ}. 

\begin{theorem}\label{coinPBQfil} Let $M_\eta$ be a $\varphi$-$\nabla$-modules over $E_\eta$, and $\alpha$ a closed point in $\overline{C}$. 
\begin{enumerate}
\item $M_\eta$ is PBQ if and only if $\mathcal E_\alpha\otimes_{E_\eta}M_\eta$ is PBQ . 
\item If $\{ P_{\eta, i}\}$ is the PBQ filtration of $M_\eta$,  
then $\{ \mathcal E_\alpha\otimes_{E_\eta}P_{\eta, i} \}$ is the PBQ filtration of $\mathcal E_\alpha\otimes_{E_\eta}M_\eta$.
\end{enumerate}
\end{theorem}

\prf{(1) follows from Corollary \ref{boundfil}.

(2) follows from (1) and the uniqueness of PBQ filtrations by Theorem \ref{genPBQfil} (1).
}

\subsection{PBQ filtrations over $\mathcal E^\dagger_\alpha,$ and $K[\hspace*{-0.2mm}[x]\hspace*{-0.2mm}]_0$}\label{pbqdag}
What we deal are only local objects in this subsection, so we drop the subscript $\alpha$ from the notation 
$x_\alpha, \mathcal E_\alpha, \mathcal E^\dagger_\alpha$ as in Section \ref{dJ}. 
We recall the fact that the PBQ filtrations over $\mathcal E$ descend to those over $K[\hspace*{-0.2mm}[x]\hspace*{-0.2mm}]_0$. 

\begin{definition}
\begin{enumerate}
\item \mbox{\rm \cite[Definition 5.1]{CT11}} A $\varphi$-$\nabla$-module $M_0$ over $K[\hspace*{-0.2mm}[x]\hspace*{-0.2mm}]_0$ 
is said to be PBQ if $\mathcal E\otimes_{K[\hspace*{-0.2mm}[x]\hspace*{-0.2mm}]_0}M_0$ is PBQ.
\item \mbox{\rm \cite[Section 12]{Oh18}}  A $\varphi$-$\nabla$-module $M^\dagger$ over $\mathcal E^\dagger$ 
is said to be PBQ if $\mathcal E\otimes_{\mathcal E^\dagger}M^\dagger$ is PBQ.
\end{enumerate}
\end{definition}

\begin{theorem}\label{PBQfillocal} 
\begin{enumerate}
\item 
\mbox{\rm \cite[Theorem 5.6]{CT11}} Let $M_0$ be a $\varphi$-$\nabla$-module over $K[\hspace*{-0.2mm}[x]\hspace*{-0.2mm}]_0$. 
Then then there exists a unique filtration $\{ P_{0, i} \}$ of $M_0$ as 
$\varphi$-$\nabla$-modules over $K[\hspace*{-0.2mm}[x]\hspace*{-0.2mm}]_0$ such that 
$\{ P_{0, i} \}$ is a $K[\hspace*{-0.2mm}[x]\hspace*{-0.2mm}]_0$-lattice of the PBQ filtration 
of $\mathcal E \otimes_{K[\hspace*{-0.2mm}[x]\hspace*{-0.2mm}]_0}M_0$. 
$\{ P_{0, i} \}$ is also called the PBQ filtration of $M_0$. 
\item \mbox{\rm \cite[Thoerem 12.7]{Oh18}} Let $M^\dagger$ be a $\varphi$-$\nabla$-module over $\mathcal E^\dagger$. 
Then then there exists a unique filtration $\{ P_i^\dagger \}$ of $M^\dagger$ as 
$\varphi$-$\nabla$-modules over $\mathcal E^\dagger$ such that 
$\{ P_i^\dagger \}$ is an $\mathcal E^\dagger$-lattice of the PBQ filtration 
of $\mathcal E \otimes_{\mathcal E^\dagger}M^\dagger$. 
$\{ P_i^\dagger \}$ is also called the PBQ filtration of $M^\dagger$. 
\end{enumerate}
In particular, an irreducible $\varphi$-$\nabla$-module over $K[\hspace*{-0.2mm}[x]\hspace*{-0.2mm}]_0$ 
(resp. $\mathcal E^\dagger$) is PBQ. 
\end{theorem}

\vspace*{1mm}

For a $\nabla$-module $M_0$ over $K[\hspace*{-0.2mm}[x]\hspace*{-0.2mm}]_0$ we define a $K$-space of solutions by 
$$
     \mathrm{Sol}(M_0) = \left\{ f : M_0 \rightarrow \mathcal A_K(0, 1)\, \left|\, \begin{array}{l}
     \mbox{\rm $f$ is $K[\hspace*{-0.2mm}[x]\hspace*{-0.2mm}]_0$-linear such that} \\
     \mbox{\rm $\frac{d}{dx}f(m) = f\left(\nabla(\frac{d}{dx})(m))\right)$ for any $m \in M_0$.}
     \end{array}\right.  \right\}
$$
where $\mathcal A_K(0, 1)$ is analytic on the open unit disk $|x| <1$, and define a $K$-space of bounded solutions by
$$
     \mathrm{Sol}_0(M_0) = \left\{ f : M_0 \rightarrow K[\hspace*{-0.2mm}[x]\hspace*{-0.2mm}]_0\, \left|\, \begin{array}{l}
     \mbox{\rm $f$ is $K[\hspace*{-0.2mm}[x]\hspace*{-0.2mm}]_0$-linear such that} \\
     \mbox{\rm $\frac{d}{dx}f(m) = f(\nabla(\frac{d}{dy})(m))$ for any $m \in M_0$.}
     \end{array}\right.  \right\}
$$
If $\mathcal E\otimes_{K[\hspace*{-0.2mm}[x]\hspace*{-0.2mm}]_0}M_0$ is solvable, 
then $\mathrm{dim}_K\, \mathrm{Sol}(M_0) = \mathrm{rank}_{K[\hspace*{-0.2mm}[x]\hspace*{-0.2mm}]_0}M_0$ 
by Christol transfer theorem \cite[Th\'eor\`eme 2]{Ch84}. Moreover, $\mathrm{Sol}(M_0)$ 
is an $F$-space over $K$ whose Frobenius is defined by 
$$
     F_{\mathrm{Sol}(M_0)}(f) = F_{K[\hspace*{-0.2mm}[x]\hspace*{-0.2mm}]_0}\circ f \circ F_{M_0}^{-1}\, \, \mbox{\rm for}\, \, f \in \sigma^\ast\mathrm{Sol}(M_0)
$$
and $\mathrm{Sol}_0(M_0)$ is an $F$-subspace over $K$ since the $q$-Frobenius $\varphi$ acts on $K[\hspace*{-0.2mm}[x]\hspace*{-0.2mm}]_0$. 
Note that $\mathrm{Sol}(M_0)$ is the dual of the $F$-space $M_0\otimes_{K[\hspace*{-0.2mm}[x]\hspace*{-0.2mm}]_0}K$ over $K$ 
$(K[\hspace*{-0.2mm}[x]\hspace*{-0.2mm}]_0 \rightarrow K,\, \sum_na_nx^n \mapsto a_0)$. 
Concerning to the bounded quotients, the proposition below holds 
which corresponds to a bounded Dwork's conjecture on logarithmic growth v.s. Frobenius slopes (see Remark \ref{logfrob} below). 

\begin{proposition}\label{bDw} Let $M_0$ be a PBQ $\varphi$-$\nabla$-module over $K[\hspace*{-0.2mm}[x]\hspace*{-0.2mm}]_0$ 
such that the maximal slope of $\mathcal E \otimes_{K[\hspace*{-0.2mm}[x]\hspace*{-0.2mm}]_0}M_0$ is $\lambda^{\mathrm{max}}$.
Then, if $\mathrm{Sol}(M_0)^{(-\lambda_{\mathrm{max}})}$ denotes the $F$-subspace of slope $-\lambda_{\mathrm{max}}$ 
in the $F$-space $\mathrm{Sol}(M_0)$ over $K$, then 
we have an equality
$$
     \mathrm{Sol}_0(M_0) = \mathrm{Sol}(M_0)^{(-\lambda_{\mathrm{max}})}. 
$$
In particular, there exists a unique $\varphi$-$\nabla$-module $L_0$ over $K[\hspace*{-0.2mm}[x]\hspace*{-0.2mm}]_0$ 
with a surjection $M_0 \rightarrow L_0$ of $\varphi$-$\nabla$-modules over $K[\hspace*{-0.2mm}[x]\hspace*{-0.2mm}]_0$ 
such that $L_0$ is isoclinic of slope $\lambda_{\mathrm{max}}$ and that 
$\mathrm{rank}_{K[\hspace*{-0.2mm}[x]\hspace*{-0.2mm}]_0}L_0 = \mathrm{Sol}(M_0)^{(-\lambda_{\mathrm{max}})}$. 
\end{proposition}

\prf{The equality follows from the proof of Dwork's conjecture \cite[Theorem A]{Oh18} (see Remark \ref{logfrob}). 
Here we give a sketch of another proof. The natural morphism 
$\mathrm{Sol}_0(M_0)\otimes_KM_0 \rightarrow K[\hspace*{-0.2mm}[x]\hspace*{-0.2mm}]_0\,\, ((f, m) \mapsto f(m))$ 
induces a surjection $M_0 \rightarrow K[\hspace*{-0.2mm}[x]\hspace*{-0.2mm}]_0\otimes_K\mathrm{Sol}_0(M_0)^\vee$ of 
$\varphi$-$\nabla$-modules over $K[\hspace*{-0.2mm}[x]\hspace*{-0.2mm}]_0$. 
If $\mathrm{Sol}_0(M_0)$ has a slope different from $-\lambda_{\mathrm{max}}$, 
then it is contradicts to the fact that $\mathcal E\otimes_{K[\hspace*{-0.2mm}[x]\hspace*{-0.2mm}]_0}M_0$ 
is PBQ with the maximal slope $\lambda_{\mathrm{max}}$. Hence $\mathrm{Sol}_0(M_0) \subset \mathrm{Sol}(M_0)^{(-\lambda_{\mathrm{max}})}$. 
The opposite inclusion follows from \cite[Theorem 6.17 (ii)]{CT09}. 
In particular, $L_0 = K[\hspace*{-0.2mm}[x]\hspace*{-0.2mm}]_0\otimes_K\mathrm{Sol}_0(M_0)^\vee$. 
}

\begin{remark}\label{logfrob}
In {}'70-{}'80 B.Dwork studied the logarithmic growth (log-growth for short) of solutions of $p$-adic 
linearly differential equations and proposed problems on the relation between log-growth and Frobenius slopes. 
After the work of Dwork, Christol and Robba, 
Chiarellotto and the author formulated the relation precisely, which is called Dwork's conjecture \cite[Conjectures 2.4, 2.5]{CT11}. 
Let $\lambda$ be a nonnegative real number, and $K[\hspace*{-0.2mm}[x]\hspace*{-0.2mm}]_\lambda$ 
a $K[\hspace*{-0.2mm}[x]\hspace*{-0.2mm}]_0$-module of analytic functions of log-growth $\leq \lambda$ on the unit open disc 
defined by 
$$
   K[\hspace*{-0.2mm}[x]\hspace*{-0.2mm}]_\lambda = \left\{\left. \sum_na_nx^n \in \mathcal A_K(0, 1)\, \right|\, 
   \underset{n \rightarrow \infty}{\mbox{\rm lim sup}}\frac{|a_n|}{(n+1)^\lambda} < \infty \right\}. 
$$
For a $\varphi$-$\nabla$-module $M_0$ over $K[\hspace*{-0.2mm}[x]\hspace*{-0.2mm}]_0$, we define a $K$-space 
of solutions of $M_0$ of log-growth $\leq \lambda$ by 
$$
\mathrm{Sol}_\lambda(M_0) = \left\{ f : M_0 \rightarrow K[\hspace*{-0.2mm}[x]\hspace*{-0.2mm}]_\lambda\, \left|\, \begin{array}{l}
     \mbox{\rm $f$ is $K[\hspace*{-0.2mm}[x]\hspace*{-0.2mm}]_0$-linear such that} \\
     \mbox{\rm $\frac{d}{dx}f(m) = f(\nabla(\frac{d}{dy})(m))$ for any $m \in M_0$.}
     \end{array}\right.  \right\}, 
$$
and 
$\mathrm{Sol}(M_0)^{(\leq \lambda - \lambda_{\mathrm{max}})}$ 
denotes the $F$-subspace of slope $\leq \lambda - \lambda_{\mathrm{max}}$ in $\mathrm{Sol}(M_0)$.  
Then Dwork's conjecture asserts that, if a $\varphi$-$\nabla$-module $M_0$ over $K[\hspace*{-0.2mm}[x]\hspace*{-0.2mm}]_0$ is PBQ, 
then the equality 
$$
   \mathrm{Sol}_\lambda(M_0) = \mathrm{Sol}(M_0)^{(\leq \lambda - \lambda_{\mathrm{max}})}
$$
holds for all $\lambda \geq 0$. Recently S.Ohkubo affirmatively solved the conjecture in 
\cite[Theorems 0.1, 0.2]{Oh18}. More precisely, Ohkubo generalized Dwork's conjecture for $\varphi$-$\nabla$-modules 
over $\mathcal E^\dagger$ and proved it. 
Since the bounded object is split by Frobenius slopes by Theorem \ref{split} 
and the full faithfulness of the functor $\mbox{\rm \bf $\Phi$M}^\nabla_{K[\hspace*{-0.2mm}[x]\hspace*{-0.2mm}]_0} 
\rightarrow \mbox{\rm \bf $\Phi$M}^\nabla_{\mathcal E}$ by \cite[Theorem 9.1]{dJ98}, 
the converse of Dwork's conjecture also holds 
(see also \cite[Section 14]{Oh18}).
\end{remark}

\subsection{Global PBQ filtration}\label{globalPBQ} Now we introduce a global PBQ filtration 
for overconvergent $F$-isocrystals on a curve. Let $C$ be a smooth connected curve over $\mathrm{Spec}\, k$. 

\begin{definition}\label{globPBQ} An overconvergent $F$-isocrystal $\mathcal M^\dagger$ 
on $C/K$ is said to be PBQ if, for an affine open dense subscheme $U$ of $C$ 
with a smooth lift $\mathcal U = \mathrm{Spec}\, A_U$ 
such that $A_U^\dagger$ admits a $q$-Frobenius $\varphi$ which is compatible with the $q$-Frobenius $\sigma$ on $K$, 
the associated $\varphi$-$\nabla$-module 
$M_\eta = E_\eta\otimes_{A_{U, K}^\dagger}\Gamma(]\overline{C}[_{\widehat{\overline{\mathcal C}}}, j^\dagger_U\mathcal M^\dagger)$ 
over $E_\eta$ is PBQ.
\end{definition}

\vspace*{2mm}

In Definition \ref{globPBQ} the definition of PBQ is independent of the choice of affine open subschemes $U$. 
We will prove the following existence theorem which is a key role to study 
the minimal slope conjecture. 

\begin{theorem}\label{globPBQfil} Let $\mathcal M^\dagger$ be an overconvergent $F$-isocrystal on $C/K$. 
Then there exists a unique filtration $0 = \mathcal P_0^\dagger \subsetneq \mathcal P_1^\dagger 
\subsetneq \cdots \subsetneq \mathcal P_r^\dagger = \mathcal M^\dagger$ 
as overconvergent $F$-isocrystals on $C/K$ such that 
$\{ E_\eta\otimes_{A_{U, K}^\dagger}\Gamma(]\overline{C}[_{\widehat{\overline{\mathcal C}}}, j^\dagger_U\mathcal P^\dagger_i)\}$ 
is the PBQ filtration of 
$M_\eta = E_\eta\otimes_{A_{U, K}^\dagger}\Gamma(]\overline{C}[_{\widehat{\overline{\mathcal C}}}, j^\dagger_U\mathcal M^\dagger)$ 
for any affine open dense subscheme $U$ of $C$. 
The filtration $\{ \mathcal P^\dagger_i\}$ is called the PBQ filtration of $\mathcal M$. 
\end{theorem}

\begin{corollary}\label{PBQirr}
An irreducible overconvergent $F$-isocrystal on $C/K$ is PBQ. 
\end{corollary}

We will give an example of PBQ filtrations at the end of this section.

\begin{proposition}\label{globalbdqt} Let $\mathcal M^\dagger$ be an overconvergent $F$-isocrystal on $C/K$ 
admitting the slope filtration of $\mathcal M$, 
and $0 = \mathcal P_0^\dagger \subsetneq \mathcal P_1^\dagger 
\subsetneq \cdots \subsetneq \mathcal P_r^\dagger = \mathcal M^\dagger$ the PBQ filtration of $\mathcal M^\dagger$.  
Then $\oplus_{i=1}^r\, \mathcal P_i/\mathcal P_i^1$ is a quotient of $\mathcal M$ as convergent $F$-isocrystals on $U/K$. 
If we put $\mathcal M^b = \mathrm{Ker}(\mathcal M \rightarrow \oplus_{i=1}^r\, \mathcal P_i/\mathcal P_i^1)$, 
then the generic fiber $E_\eta \otimes \mathcal M/\mathcal M^b$ of $\mathcal M/\mathcal M^b$ 
is isomorphic to $M_\eta/M^b_\eta$. The convergent $F$-isocrystal 
$\mathcal M/\mathcal M^b$ is called the bounded quotient of $\mathcal M$. 
\end{proposition}

\prf{Since $\mathcal P_1^\dagger$ is the maximal PBQ submodule of $\mathcal M^\dagger$, one has 
$\mathcal P_1/\mathcal P_1^1 = \mathcal M/\mathcal M^1$. Hence, there exists a canonical projection  
$\mathcal M \rightarrow \mathcal M/\mathcal M^1 \oplus  (\mathcal M/\mathcal P_1)/(\mathcal M/\mathcal P_1)^1 \cong \oplus_{i=1}^r\, \mathcal P_i/\mathcal P_i^1$ 
by induction on $r$.  The rest follows from Theorems \ref{split} and \ref{globPBQfil}.
}

\vspace*{2mm}

Let us study several properties of PBQ overconvergent $F$-isocrystals. 

\begin{lemma}\label{subqutsh} 
\begin{enumerate}
\item If $\mathcal M_1^\dagger$ and $\mathcal M_2^\dagger$ are PBQ overconvergent $F$-isocrystals on $C/K$ of a same maximal slope, 
then so is the direct sum $\mathcal M_1^\dagger \oplus \mathcal M_2^\dagger$. 
\item Let $\theta^\dagger : \mathcal M^\dagger \rightarrow \mathcal N^\dagger$ be 
a surjection of overconvergent $F$-isocrystals on $C/K$. 
If $\mathcal M^\dagger$ is PBQ and $\mathcal N^\dagger \ne 0$, then $\mathcal N^\dagger$ is also PBQ. Moreover, 
the maximal slope of $N_\eta$ is equal to that of $M_\eta$. 
\end{enumerate}
\end{lemma}

\prf{The assertions follow from Lemma \ref{subqut}.}

\begin{lemma}\label{inver} Let $f : C' \rightarrow C$ be a finite etale morphism 
of smooth connected curves over $\mathrm{Spec}\, k$. 
Let $\mathcal M^\dagger$ (resp. $\mathcal N^\dagger$) be an overconvergent $F$-isocrystal on $C/K$ (resp. $C'/K$), 
and denotes the direct image of $\mathcal M^\dagger$ on $C/K$ (resp. the inverse image of $\mathcal N^\dagger$ on $C'/K$) 
by $f_\ast\mathcal M^\dagger$ (resp. $f^\ast\mathcal N^\dagger$). 
\begin{enumerate}
\item If $\mathcal N^\dagger$ is PBQ, then so is $f^\ast\mathcal N^\dagger$. 
\item If $f_\ast\mathcal M^\dagger$ is PBQ, the so is $\mathcal M^\dagger$. 
\end{enumerate}
\end{lemma} 

Note that the inverse image of overconvergent isocrystals is defined in \cite[Definitions 2.3.2 (iv)]{Be96} 
and the direct image for finite etale morphisms is studied in \cite[5.1]{Ts02}. 

\vspace*{3mm}

\prf{(1) By the projection formula there is an isomorphism
 $f_\ast f^\ast\mathcal N^\dagger = f_\ast(j_{C'}^\dagger\mathcal O_{]\overline{C'}[}) \otimes_{j_C^\dagger\mathcal O_{]\overline{C}[}} \mathcal N^\dagger$. 
Since the direct image $f_\ast(j_{C'}^\dagger\mathcal O_{]\overline{C'}[})$ of the constant overconvergent $F$-isocrystal 
$j_{C'}^\dagger\mathcal O_{]\overline{C'}[}$ is unit-root, the direct image $f_\ast f^\ast\mathcal N^\dagger$ is PBQ by Proposition \ref{tenunit}. 
Since the constant overconvergent $F$-isocrystal $j_C^\dagger\mathcal O_{]\overline{C}[}$ is a direct summand of 
the direct image $f_\ast(j_{C'}^\dagger\mathcal O_{]\overline{C'}[})$, $\mathcal N^\dagger$ is PBQ by Lemma \ref{subqutsh}. 

(2) Since the adjoint morphism $f^\ast f_\ast\mathcal M^\dagger \rightarrow \mathcal M^\dagger$ is surjective, 
the assertions follows from Lemma \ref{subqutsh} and Lemma \ref{inver}.
}

\vspace*{2mm}

Now we return to prove Theorem \ref{globPBQfil}.  
The following are key lemmas to construct unit-root overconvergent $F$-isocrystals 
and overconvergent $F$-subisocrystals from 
several compatible generic and local data. 

\begin{lemma}\label{rankoneext} Suppose that $C$ is affine with a smooth lift $\mathcal C = \mathrm{Spec}\, A_C$ 
such that $A_C^\dagger$ admits a $q$-Frobenius $\varphi$ which is compatible with the $q$-Frobenius $\sigma$ on $K$. 
Let $L_\eta$ be a unit-root $\varphi$-$\nabla$-module over $E_\eta$ satisfying the conditions as follows: 
\begin{list}{}{}
\item[\mbox{\rm (i)}] For any closed point $\alpha \in C$, there exists 
a $K_\alpha[\hspace*{-0.2mm}[x_\alpha]\hspace*{-0.2mm}]_0$-lattice $L_{0, \alpha}$ of $L_\alpha 
= \mathcal E_\alpha\otimes_{E_\eta}L_\eta$. 
\item[\mbox{\rm (ii)}]  For any closed point $\alpha \in \overline{C} \setminus C$, there exists 
an $\mathcal E^\dagger_\alpha$-lattice $L^\dagger_\alpha$ of $L_\alpha = \mathcal E_\alpha\otimes_{E_\eta}L_\eta$. 
\end{list}
Then there exists a unit-root overconvergent $F$-isocrystal $\mathcal L^\dagger$ on $C/K$ 
such that the  $\varphi$-$\nabla$-module $E_\eta\otimes_{A^\dagger_{C, K}}\Gamma(]\overline{C}[_{\widehat{\overline{\mathcal C}}}, \mathcal L^\dagger)$ 
over $E_\eta$ is isomorphic to the given $L_\eta$. $\mathcal L^\dagger$ is unique up to isomorphisms. 
\end{lemma}

\prf{Let $G_{k(C)}$ be the absolute Galois group of the function field $k(C)$ of $C$. 
By Katz correspondence between unit-root $F$-spaces and $p$-adic representations by \cite[Proposition 4.1.1]{Ka72}, 
one has a continuous representation
$$
      \rho :  G_{k(C)} \rightarrow \mathrm{GL}_s(K_\sigma). 
$$
By the condition (i) we know that $\rho$ is unramified at $\alpha \in C$ and $\rho$ has a finite local monodromy 
at $\alpha \in \overline{C} \setminus C$ by \cite[Theorem 4.2.6]{Ts98}. 
Hence $\rho$ corresponds to an unit-root overconvergent $F$-isocrystals on $C/K$ which has a desired property 
by \cite[Theorem 7.2.3]{Ts98}. The uniqueness follows from the equivalence of Katz correspondence. 
}

\begin{lemma}\label{extext} \mbox{\rm (c.f.} \cite[Remark 5.10]{Ke16}\mbox{\rm )} 
Suppose that $C$ is affine with a smooth lift $\mathcal C = \mathrm{Spec}\, A_C$ 
such that $A_C^\dagger$ admits a $q$-Frobenius $\varphi$ which is compatible with the $q$-Frobenius $\sigma$ on $K$. 
Let $\mathcal M^\dagger$ be an overconvergent $F$-isocrystal on $C/K$, 
put $M^\dagger = \Gamma(]\overline{C}[_{\widehat{\overline{\mathcal C}}}, \mathcal M^\dagger)$, and 
$L_\eta$ a $\varphi$-$\nabla$-submodule of $M_\eta = E_\eta \otimes_{A^\dagger_{C, K}}M^\dagger$ over $E_\eta$. 
Suppose that 
\begin{list}{}{}
\item[\mbox{\rm (i)}] for any closed point $\alpha$ of $C$ 
there exists a $K_\alpha[\hspace*{-0.2mm}[x_\alpha]\hspace*{-0.2mm}]_0$-lattice $L_{0, \alpha}$ of 
$L_\alpha = \mathcal E_\alpha \otimes_{E_\eta}L_\eta$, and 
\item[\mbox{\rm (ii)}] for any closed point $\alpha$ of $\overline{C} \setminus C$ 
there exists an 
$\mathcal E^\dagger_\alpha$-lattice $L_\alpha^\dagger$ of $L_\alpha = \mathcal E_\alpha \otimes_{E_\eta}L_\eta$.
\end{list}
(Note that $L_{0, \alpha}$ (resp. $L^\dagger_\alpha$) is canonically a $\varphi$-$\nabla$-submodule of 
$M_{0, \alpha} = K_\alpha[\hspace*{-0.2mm}[x_\alpha]\hspace*{-0.2mm}]_0\otimes_{A^\dagger_{C, K}}M^\dagger$ 
(resp. $M^\dagger_\alpha = \mathcal E^\dagger_\alpha\otimes_{A^\dagger_{C, K}}M^\dagger$) by full faithfulness). 
Then there exists a unique overconvergent $F$-subisocrystal $\mathcal L^\dagger$ of $\mathcal M^\dagger$ on $C/K$ such that, 
if $L^\dagger=  \Gamma(]\overline{C}[_{\widehat{\overline{\mathcal C}}}, \mathcal L^\dagger)$, 
then $L_\eta = E_\eta \otimes_{A^\dagger_{C, K}}L^\dagger$, 
$L_{0, \alpha} =  K_\alpha[\hspace*{-0.2mm}[x_\alpha]\hspace*{-0.2mm}]_0\otimes_{A^\dagger_{C, K}}L^\dagger$ for any closed point $\alpha$ of $C$, 
and $L_\alpha^\dagger =  \mathcal E_\alpha^\dagger \otimes_{A^\dagger_{C, K}}L^\dagger$ for any closed point $\alpha$ of $\overline{C} \setminus C$ 
as $\varphi$-$\nabla$-modules. 
\end{lemma}

\prf{First we prove the assertion in the case where $\mathrm{dim}_{E_\eta}L_\eta = 1$. When $\mathrm{dim}_{E_\eta}L_\eta = 1$, 
Lemma \ref{rankoneext} implies that 
there exists an overconvergent $F$-isocrystal $\mathcal L^\dagger$ on $C/K$ which is determined by the data $L_\eta, L_\alpha\, (\alpha \in C)$ 
and $L_\alpha^\dagger\, (\alpha \in \overline{C} \setminus C)$. What we want is a 
nontrivial morphism $\mathcal L^\dagger \rightarrow \mathcal M^\dagger$ of overconvergent $F$-isocrystals which 
induces the conditions $L_\eta = E_\eta \otimes_{A^\dagger_{C, K}}L$ and so on. 
Replacing $\mathcal M^\dagger$ by $(\mathcal L^\dagger)^\vee\otimes_{j^\dagger_C\mathcal O_{]\overline{C}[_{\widehat{\overline{\mathcal C}}}}}\mathcal M^\dagger$, there exist 
\begin{list}{}{}
\item[($\eta$){\,\,}] an element $e_\eta \in L_\eta$ such that $\varphi_{M_\eta}(e_\eta) = e_\eta$, 
\item[$(\alpha)_0$] an element $e_\alpha \in L_{0, \alpha}$ such that $\varphi_{M_\alpha}(e_\alpha) = e_\alpha$ 
for any closed point $\alpha$ of $C$, and 
\item[$(\alpha)^\dagger$] an element $e_\alpha \in L_\alpha^\dagger$ such that $\varphi_{M_\alpha}(e_\alpha) = e_\alpha$ 
for any closed point $\alpha$ of $\overline{C} \setminus C$.
\end{list}
Since $\mathcal E_\alpha \otimes_{E_\eta} L_\eta = \mathcal E_\alpha \otimes_{K_\alpha[\hspace*{-0.2mm}[x_\alpha]\hspace*{-0.2mm}]_0}L_{0, \alpha}$ 
for any closed point $\alpha$ of $C$, there exists an element $u_\alpha \in \mathcal E_\alpha$ such that 
$$
         e_\eta = u_\alpha e_\alpha
$$
Since both Frobenius in the equality 
$\mathcal E_\alpha \otimes_{E_\eta} M_\eta = \mathcal E_\alpha \otimes_{K_\alpha[\hspace*{-0.2mm}[x_\alpha]\hspace*{-0.2mm}]_0}M_{0, \alpha}$ commute 
with each other, $u_\alpha$ is contained in $(K_\alpha)_\sigma$. We now change $e_\alpha$ by $u_\alpha e_\alpha \in L_{0, \alpha}$. The similar holds 
for $e_\alpha$ for any closed point $\alpha$ of $\overline{C} \setminus C$. 
If $e_\eta$ belongs to $M^\dagger$, then it determines a morphism 
$j^\dagger_C\mathcal O_{]\overline{C}[_{\widehat{\overline{\mathcal C}}}} 
\rightarrow \mathcal (L^\dagger)^\vee \otimes_{j^\dagger_C\mathcal O_{]\overline{C}[_{\widehat{\overline{\mathcal C}}}}} \mathcal M^\dagger$ 
of overconvergent $F$-isocrystals which satisfies the desired conditions. 
In order to prove $e_\eta \in M^\dagger$ it is sufficient to prove 
$$
        M^\dagger = \underset{\tiny \begin{array}{c} \alpha \in C \\ \mbox{\rm a closed point}\end{array}}{\cap}M_\eta \cap M_{0, \alpha}\, \, 
        \bigcap \underset{\tiny \begin{array}{c} \alpha \in \overline{C} \setminus C \\ \mbox{\rm a closed point}\end{array}}{\cap} M_\eta \cap M_\alpha^\dagger
$$
in $M_\eta$. Since $M^\dagger$ is a finite generated projective $A^\dagger_{C, K}$-module, 
the equality above follows from Lemma \ref{intersect}. 

Now we treat the general case. Put $s = \mathrm{dim}_{E_\eta}L_\eta$. Consider the $s$-th exterior product $\wedge^s \mathcal M^\dagger$ of $\mathcal M^\dagger$ 
and the data $\wedge^sL_\eta, \wedge^sL_{0, \alpha}, \wedge^sL_\alpha^\dagger$. 
Then there exist an overconvergent $F$-isocrystal $\mathcal N^\dagger$ of rank one and 
an injective morphism $\mathcal N^\dagger \rightarrow \wedge^s\mathcal M^\dagger$ which is determined by the data 
$\wedge^sL_\eta, \wedge^sL_{0, \alpha}, \wedge^sL_\alpha^\dagger$ by the former part of this proof. 
We define an overconvergent $F$-isocrystal $\mathcal L^\dagger$ on $C/K$ 
by 
$$
      \mathcal L^\dagger = (\mathcal N^\dagger)^\vee \otimes_{j^\dagger_C\mathcal O_{]\overline{C}[_{\widehat{\overline{\mathcal C}}}}}
      \mathrm{Ker}(\mathcal N^\dagger \otimes_{j^\dagger_C\mathcal O_{]\overline{C}[_{\widehat{\overline{\mathcal C}}}}} \mathcal M^\dagger \rightarrow \wedge^{s+1}\mathcal M^\dagger). 
$$
Here $\mathcal N^\dagger \otimes \mathcal M^\dagger \rightarrow \wedge^{s+1}\mathcal M^\dagger$ is defined by $n\otimes m \mapsto n\wedge m$. 
Then $E_\eta \otimes_{A^\dagger_{C, K}}\Gamma(]\overline{C}[_{\widehat{\overline{\mathcal C}}}, \mathcal L^\dagger) = L_\eta$, 
$K_\alpha[\hspace*{-0.2mm}[x_\alpha]\hspace*{-0.2mm}]_0 \otimes_{A^\dagger_{C, K}}\Gamma(]\overline{C}[_{\widehat{\overline{\mathcal C}}}, \mathcal L^\dagger) 
= L_{0, \alpha}$ and $\mathcal E^\dagger_\alpha \otimes_{A^\dagger_{C, K}}\Gamma(]\overline{C}[_{\widehat{\overline{\mathcal C}}}, \mathcal L^\dagger) 
= L^\dagger_\alpha$ hold by our construction. Hence we obtain a desired 
overconvergent $F$-subisocrystal $\mathcal L^\dagger$ of $\mathcal M^\dagger$. 
The uniqueness follows from Lemma \ref{rankoneext} and our construction. 
}

\vspace*{3mm}

\noindent
{\sc Proof of Theorem \ref{globPBQfil}.} We first prove the uniqueness. It is sufficient to prove the uniqueness 
on the first step $\mathcal P_1^\dagger$ of PBQ filtration. 
Suppose $\mathcal P_1^\dagger$ and $(\mathcal P_1')^\dagger$ are first steps of two 
PBQ filtrations of $\mathcal M^\dagger$. Then the associated $\varphi$-$\nabla$-modules  
$P_{1, \eta}$ and $P_{1, \eta}+P_{1, \eta}'$ over $E_\eta$ coincides with each other in $M_\eta$  
by Theorem \ref{genPBQfil}. Hence $\mathcal P_1^\dagger$ and $(\mathcal P_1')^\dagger$ are equal.

Now we prove the existence. Since the problem is local on $C$ by patching and the full faithfulness of restriction functors \cite[Theorem 6.3.1]{Ts02}, 
we may assume that $C$ is affine with a smooth lift $\mathcal C = \mathrm{Spec}\, A_C$ 
such that $A_C^\dagger$ admits a $q$-Frobenius $\varphi$ which is compatible with the $q$-Frobenius $\sigma$ on $K$. 
The unique existence of PBQ filtration of the $\varphi$-$\nabla$-module  
$K_\alpha[\hspace*{-0.2mm}[x_\alpha]\hspace*{-0.2mm}]_0 \otimes_{A_{C, K}^\dagger}\Gamma(]\overline{C}[_{\widehat{\overline{\mathcal C}}}, \mathcal M^\dagger)$ 
for a closed point $\alpha \in C$ 
(resp. $\mathcal E^\dagger_\alpha \otimes_{A_{C, K}^\dagger}\Gamma(]\overline{C}[_{\widehat{\overline{\mathcal C}}}, \mathcal M^\dagger)$ 
for a closed point $\alpha \in \overline{C} \setminus C$) gives the comparison data 
in Lemma \ref{extext} for each step of PBQ filtration by Theorem \ref{PBQfillocal}. 
Hence we obtain a desired PBQ filtration of $\mathcal M^\dagger$ by Lemma \ref{extext}. 
\hspace*{\fill} $\Box$

\subsection{Example}\label{legtwist} Let us give an example of PBQ filtration with two steps. 

Suppose $p$ is an odd prime number and 
fix an embedding $\overline{\mathbb Q} \subset \overline{\mathbb Q}_p$ 
from an algebraic closure $\overline{\mathbb Q}$ of the field $\mathbb Q$ of rational numbers to an algebraic closure 
$\overline{\mathbb Q}_p$ of the field $\mathbb Q_p$ of $p$-adic numbers.  
Let $\sigma$ be the $p$-Frobenius $\mathrm{id}_{\mathbb Q_p}$ on $\mathbb Q_p$. 
Let $C = \mathbb P^1_{\mathbb F_p} \setminus \{0, 1, \infty\}$ be a smooth curve over $\mathrm{Spec}\, \mathbb F_p$ 
with the standard coordinate $z$, $\overline{C} = \mathbb P^1$ the smooth completion of $C$, 
$X$ an affine scheme defined by the equation 
$$
       y^2 = x(x-1)(x-z)
$$
over $C$, 
and $\overline{X}$ the Legendre family of elliptic curves which is a completion of $X$ over $C$. 
Let $\mathcal L^\dagger$ be an overconvergent $F$-isocrystal on $C/\mathbb Q_p$ of rank $2$ which is defined by 
the first relative rigid cohomology of $X/C$, $\chi$ a nontrivial quadratic character on $C$ 
corresponding to the double cover defined by the equation $w^2 = z(z-1)$, and $\mathcal L^\dagger(\chi)$ 
(resp. $j^\dagger_C\mathcal O_{]\overline{C}[_{\widehat{\overline{\mathcal C}}}}(\chi)$) 
the twist of $\mathcal L^\dagger$ (resp. $j^\dagger_C\mathcal O_{]\overline{C}[_{\widehat{\overline{\mathcal C}}}}$) by the rank $1$ 
overconvergent $F$-isocrystal on $C/\mathbb Q_p$ corresponding to $\chi$. 
Note that $\mathcal L^\dagger(\chi)$ is pure of weight $1$ (recall the definition of weights in Definition \ref{weight}). 
Since $\mathcal L^\dagger$ 
(resp. $\chi$) ramifies tamely at $0, 1, \infty$ (resp. $0, 1$) with exponents 
$1, 1, 1/2$ (resp. $1/2, 1/2$), respectively, over the projective line $\overline{C}$, 
the support forgetting map 
$H^1_{\mathrm{rig}}(C/\mathbb Q_p, \mathcal L^\dagger(\chi)) 
\rightarrow H^1_{\mathrm{rig}}(C/\mathbb Q_p, \mathcal L^\dagger(\chi))$ is isomorphic 
and the first rigid cohomology $H^1_{\mathrm{rig}, c}(C/\mathbb Q_p, \mathcal L^\dagger(\chi))$ is of dimension $2$ 
and pure of weight $2$ by Euler-Poincar\'e formula. 
Moreover we have 
\begin{enumerate}
\item[\mbox{\rm (a)}] if $p \equiv 1\, (\mathrm{mod}\, 4)$, then 
the slopes of Frobenius action on $H^1_{\mathrm{rig}}(C/\mathbb Q_p, \mathcal L^\dagger(\chi))$ are $0$ and $2$; 
\item[\mbox{\rm (b)}] if $p \equiv 3\, (\mathrm{mod}\, 4)$, then 
the slopes of Frobenius action on $H^1_{\mathrm{rig}}(C/\mathbb Q_p, \mathcal L^\dagger(\chi))$ are $1$ and $1$.
\end{enumerate}
Indeed, (a) and (b) follows from the fact that the arithmetic family $\overline{Y}$ 
of singular $K3$ surfaces over $\mathrm{Spec}\, \mathbb Z[1/2]$ 
which is a minimal desingularization of completion of the affine scheme $Y$ over $C$ defined by the equation 
$$
      y^2 = z(z-1)x(x-1)(x-z) 
$$
has a complex multiplication of $\mathbb Q(\sqrt{-1})$. If 
the zeta function of the elliptic curve $E$ defined by $y^2 = x^3-x$ over $\mathbb F_p$ 
is $\mathrm{Z}(E/{\mathbb F_p}, t) 
= (1-\pi_pt)(1-\overline{\pi}_pt)/(1-t)(1-pt)$, then 
the zeta function of $\overline{Y}/{\mathbb F_p} =\overline{Y} \times_{\mathrm{Spec}\, \mathbb Z[1/2]}\mathrm{Spec}\, \mathbb F_p$ 
over $\mathbb F_p$ is given by 
$$
      \mathrm{Z}(\overline{Y}/{\mathbb F_p}, t) 
      = \frac{1}{(1-t)(1-p^2t)(1-pt)^{20}(1-\rho\pi_p^2t)(1-\overline{\rho}\overline{\pi}_p^2t)}
$$
for a root $\rho$ of unity in $\mathbb Q(\sqrt{-1})$ 
by \cite[the case A in pp.290, 291]{SB85} applying the work of \cite{IS}. 
Here $\overline{\pi}_p$ (resp. $\overline{\rho}$) is the complex conjugate of $\pi_p$ (resp. $\rho$). 
Hence the $L$-function of $\mathcal L^\dagger(\chi)$ is 
$$
  L(C/{\mathbb F_p}, \mathcal L^\dagger(\chi), t) = (1-\rho\pi_p^2t)(1-\overline{\rho}\overline{\pi}_p^2t) 
$$
since the Leray spectral sequence induces 
an isomorphism $H^1_{\mathrm{rig}, c}(C/\mathbb Q_p, \mathcal L^\dagger(\chi)) \cong H^2_{\mathrm{rig}, c}(Y/\mathbb Q_p)$ 
of $F$-spaces. 
More precisely, $\rho = 1$ if $p \equiv 1\, (\mathrm{mod}\, 4)$ and $\rho = \sqrt{-1}$ if $p \equiv 3\, (\mathrm{mod}\, 4)$ 
by calculating $\sharp Y(\mathbb F_p)$ modulo $4$. 

Suppose $p\, \equiv 1\, (\mathrm{mod}\, 4)$. If we fix the unit root element $\pi_p$ under the fixed 
embedding $\mathbb Q(\sqrt{-1}) \subset \mathbb Q_p$
then there exists a nontrivial extension 
$$
      0 \rightarrow \mathcal L^\dagger(\chi) \rightarrow \mathcal M^\dagger \rightarrow 
      j^\dagger_C\mathcal O_{]\overline{C}[_{\widehat{\overline{\mathcal C}}}}(\pi_p^2) \rightarrow 0
$$
of overconvergent $F$-isocrystals on $C/\mathbb Q_p$ by (a) where 
$j^\dagger_C\mathcal O_{]\overline{C}[_{\widehat{\overline{\mathcal C}}}}(\pi_p^2)$ denotes 
the twist of $j^\dagger_C\mathcal O_{]\overline{C}[_{\widehat{\overline{\mathcal C}}}}$ by 
$\pi_p^2$ time Frobenius. Since the maximal generic slope of $\mathcal L^\dagger(\chi)$ is $1$ and 
the slope of 
$j^\dagger_C\mathcal O_{]\overline{C}[_{\widehat{\overline{\mathcal C}}}}(\pi_p^2)$ is $0$, 
$\mathcal M^\dagger$ has two steps of the PBQ filtration as above.

\section{A local version of the minimal slope conjecture}\label{locver}

In this section we will study a local version of the minimal slope conjecture, Theorem \ref{loc}. 
Let us keep the notation in section \ref{dJ}. 

\subsection{A local version}

\begin{definition} A $\varphi$-$\nabla$-module $M^\dagger$ over $\mathcal E^\dagger$ 
is saturated if the canonical morphism $M^\dagger \rightarrow M/M^1$ is injective. 
\end{definition}

In this section we will prove the theorem below, which is easily deduced from Theorem \ref{locrel}. 

\begin{theorem}\label{loc} Let $M^\dagger$ and $N^\dagger$ be 
$\varphi$-$\nabla$-modules over $\mathcal E^\dagger$, and $h : N/N^1 \rightarrow M/M^1$ a morphism 
of $\varphi$-$\nabla$-modules over $\mathcal E$. 
Suppose either 
\begin{list}{}{}
\item[\mbox{\rm (i)}] both $M^\dagger$ and $N^\dagger$ are irreducible and $h$ is nontrivial, or 
\item[\mbox{\rm (ii)}] both $M^\dagger$ and $N^\dagger$ are saturated and PBQ and $h$ is an isomorphism. 
\end{list}
Then there exists a unique isomorphism 
$g^\dagger : N^\dagger \rightarrow M^\dagger$ of $\varphi$-$\nabla$-modules over $\mathcal E^\dagger$ 
such that the induced morphism $N/N^1 \rightarrow M/M^1$ 
between the maximal slope quotients by $g^\dagger$ coincides with the given $h$. 
\end{theorem}

\subsection{Properties of saturated $\varphi$-$\nabla$-modules}\label{satloc}

\begin{definition} Let $M$ be a $\varphi$-$\nabla$-module over $\mathcal E$ which has a unique slope of Frobenius $\varphi_M$, 
$N^\dagger$ a $\varphi$-$\nabla$-module over $\mathcal E^\dagger$, and $N^\dagger \rightarrow M$ 
an injection of $\mathcal E^\dagger$-spaces which commutes with connections and Frobenius. 
A $\varphi$-$\nabla$-module $N^\dagger$ over $\mathcal E^\dagger$ is quasi-saturated in $M$ 
if $N^\dagger$ generates $M$ as an $\mathcal E$-space, in other words, the associated morphism $N \rightarrow M$ 
of $\varphi$-$\nabla$-modules over $\mathcal E$ is surjective. 
\end{definition}

\vspace*{2mm}

We gives several properties on the notion of saturated and quasi-saturated. 
Let $M$ be a $\varphi$-$\nabla$-module over $\mathcal E$ which has a unique slope, 
and $N^\dagger$ a nontrivial $\varphi$-$\nabla$-module over $\mathcal E^\dagger$. 
Let $\theta^\dagger : N^\dagger \rightarrow M$ be an $\mathcal E^\dagger$-morphism which commutes with connections and Frobenius, 
and $L^\dagger$ the image of $\theta^\dagger$. Since $L^\dagger$ is of finite dimension over $\mathcal E^\dagger$, 
the following proposition holds. 

\begin{proposition}\label{quasi} With the notation as above, 
suppose $\theta^\dagger$ is nontrivial. Then $L^\dagger$ is a nontrivial $\varphi$-$\nabla$-module over $\mathcal E^\dagger$ 
by the induced connection and Frobenius. In particular, 
any irreducible $\varphi$-$\nabla$-module over $\mathcal E^\dagger$ is saturated. 
\end{proposition}

\begin{proposition}\label{irrtop} With the notation as above, 
suppose $\theta^\dagger$ is injective. Then the maximal slope of $N^\dagger$ coincides with the slope of $M$. 
\end{proposition}

\prf{Let $P^\dagger$ be the maximally PBQ submodule of $N^\dagger$. 
If the maximal slope of $P$, which is equal to that of $N$, is not equal to the slope of $M$, then $\theta^\dagger|_{P^\dagger}$ is a zero map 
by Lemma \ref{subqut} (2). This contradicts the injectivity of $\theta^\dagger$. 
}

\begin{proposition}\label{irrdim} With the notation as above, suppose the induced morphism $\theta : N \rightarrow M$ from $\theta^\dagger$ 
is surjective. 
 \begin{enumerate}
\item $L^\dagger$ is quasi-saturated in $M$ and 
the inequality $\mathrm{dim}_{\mathcal E}L/L^1 \geq \mathrm{dim}_{\mathcal E}M$ holds. 
\item If $\mathrm{dim}_{\mathcal E}L/L^1 = \mathrm{dim}_{\mathcal E}M$, then $L^\dagger$ is saturated. 
\end{enumerate}
\end{proposition}

\prf{(1) We may suppose that $\theta^\dagger$ is injective by Proposition \ref{quasi}, and then the maximal slope of $L$ coincides with that of $M$. 
Since $L^\dagger$ is the image of $\theta^\dagger$, $L^\dagger$ generates $M$ as an $\mathcal E$-space. Moreover, 
the induced morphism $L \rightarrow M$ is surjective as a morphism of $\varphi$-$\nabla$-modules over $\mathcal E$. 

(2) Since $L \rightarrow M$ is surjective by definition of $L^\dagger$, the hypothesis implies the natural morphism $L/L^1 \rightarrow M$ 
is an isomorphism. Hence $L^\dagger$ is saturated. 
}

\begin{proposition}\label{stt} Let $M^\dagger$ be a $\varphi$-$\nabla$-module over $\mathcal E^\dagger$, 
and $L^\dagger$ a $\varphi$-$\nabla$-module over $\mathcal E^\dagger$ which is defined by the image of 
the canonical morphism $\theta^\dagger : M^\dagger \rightarrow M/M^1$. 
\begin{enumerate}
\item The induced morphism $M/M^1 \rightarrow L/L^1$ from $\theta$ is an isomorphism. 
\item $L^\dagger$ is saturated. 
\item If $M^\dagger$ is saturated, then $\theta^\dagger : M^\dagger \rightarrow L^\dagger$ is an isomorphism. 
\end{enumerate}
\end{proposition}

\prf{Since $\theta^\dagger : M^\dagger \rightarrow L^\dagger$ is surjective, the induced morphism $M/M^1 \rightarrow L/L^1$ by $\theta$ is surjective. 
(1) and (2) follows from Proposition \ref{irrdim}. (3) is trivial. 
}

\begin{definition}\label{saturation} For a $\varphi$-$\nabla$-module $M^\dagger$ over $\mathcal E^\dagger$, 
we define the saturation $M^{\mathrm{sat}, \dagger}$ of $M^\dagger$ by the $\varphi$-$\nabla$-module $L^\dagger$ in Proposition \ref{stt}. 
\end{definition}

\subsection{Saturated v.s. quasi-saturated}

\begin{proposition}\label{irrecontainer} Let $M$ be a $\varphi$-$\nabla$-module over $\mathcal E$, 
and $N^\dagger$ a $\varphi$-$\nabla$-module over $\mathcal E^\dagger$. 
If $N^\dagger$ is quasi-saturated in $M$, then the equality below holds:
$$
    \mbox{\rm dim}_{\mathcal E}\, N/N^1 = \mbox{\rm dim}_{\mathcal E}\, M. 
$$
In particular, $N^\dagger$ is saturated. 
\end{proposition}

\prf{The assertion follows from Theorem \ref{rkr} in Section \ref{dJ}.}

\begin{corollary}\label{localirr} Let $N^\dagger$ be a $\varphi$-$\nabla$-module over $\mathcal E^\dagger$. 
Then the following conditions are equivalent.
\begin{list}{}{}
\item[\mbox{\rm (i)}] $N^\dagger$ is irreducible. 
\item[\mbox{\rm (ii)}] $N^\dagger$ is PBQ and saturated and $N/N^1$ is irreducible. 
\end{list}
\end{corollary}

\prf{(i) $\Rightarrow$ (ii) : Suppose $M$ is a nontrivial quotient of $N/N^1$ as a $\varphi$-$\nabla$-module over $\mathcal E$. 
Since $N^\dagger$ is irreducible, the natural homomorphism $N^\dagger \rightarrow M$ is injective and its image 
is a $\varphi$-$\nabla$-module over $\mathcal E^\dagger$ by Proposition \ref{quasi}. 
Since $N^\dagger$ is quasi-saturated in $M$, 
we have $M = N/N^1$ by Proposition \ref{irrecontainer} so that $N/N^1$ is irreducible. The rest hold by Theorem \ref{PBQfillocal}. 

(ii) $\Rightarrow$ (i) : Suppose $L^\dagger$ is a nontrivial subobject of $N^\dagger$. Since $L^\dagger \subset N^\dagger \subset N/N^1$, 
the maximal slope of $L$ is equal to 
that of $N$ by Proposition \ref{irrtop}. Since $N/N^1$ is irreducible, 
the maximal slope of $N/L$ is less than the maximal slope of $N$ if $N/L$ is not a zero object. 
Applying Lemma \ref{subqut} (2) to the surjection $N^\dagger \rightarrow N^\dagger/L^\dagger$, 
we have $L^\dagger = N^\dagger$. Hence $N^\dagger$ is irreducible.}

\vspace*{3mm}

Now we give a version of Theorem \ref{loc} for general $\varphi$-$\nabla$-modules over $\mathcal E^\dagger$. 

\begin{theorem}\label{locrel} Let $M^\dagger$ (resp. $N^\dagger$) be a $\varphi$-$\nabla$-module over $\mathcal E^\dagger$ 
with the PBQ filtration $0 = P_0^\dagger \subsetneq P_1^\dagger \subsetneq \cdots \subsetneq P_r^\dagger = M^\dagger$ 
(resp. $0 = Q_0^\dagger \subsetneq Q_1^\dagger \subsetneq \cdots \subsetneq Q_s^\dagger = N^\dagger$), 
and $M/M^b$ (resp. $N/N^b$) the maximally bounded quotient of $M$ (resp. $N$). If there exists an isomorphism 
$h : N/N^b \rightarrow M/M^b$ of $\varphi$-$\nabla$-modules over $\mathcal E$, then $r = s$ and 
there is a unique isomorphism $g^\dagger : \oplus_{i=1}^rQ_i^{\mathrm{sat}, \dagger} \rightarrow \oplus_{i=1}^rP_i^{\mathrm{sat}, \dagger}$ 
such that the induced diagram 
$$
       \begin{array}{ccc} 
            \oplus_{i=1}^r Q_i^{\mathrm{sat}} 
             &\overset{g}{\rightarrow} &\oplus_{i=1}^r P_i^{\mathrm{sat}} \\
             \downarrow & &\downarrow \\
             N/N^b &\underset{h}{\rightarrow} &M/M^b
     \end{array}
$$
is commutative. 
\end{theorem} 

\prf{Since $M/M^b \cong \oplus_{i=1}^rP_i/P_i^1$ and $N/N^b \cong \oplus_{i=1}^sQ_i/Q_i^1$ by Theorems \ref{genPBQfil} and 
\ref{PBQfillocal} and since $h$ is an isomorphism, 
we have $r=s$ and $h|_{Q_i/Q_i^1} : Q_i/Q_i^1 \rightarrow P_i/P_i^1$ is an isomorphism 
of $\varphi$-$\nabla$-modules over $\mathcal E$ for any $i$. 
Hence it is sufficient to prove the assertion when $r=s=1$. 

Suppose $r=s=1$, that is, both $M^\dagger$ and $N^\dagger$ are PBQ. 
We may assume both $M^\dagger$ and $N^\dagger$ are saturated by Proposition \ref{stt}. 
Let $L^\dagger$ be the image of $M^\dagger\oplus N^\dagger \rightarrow M/M^1\, ((m, n) \mapsto m + h(n))$. Then 
$L^\dagger$ is the PBQ with the same maximal slope of $M$ 
and includes both $M^\dagger$ and $N^\dagger$. 
Since $L^\dagger$ is quasi-saturated in $M/M^1$, we have 
$$
       \mathrm{dim}_{\mathcal E}L/L^1 = \mathrm{dim}_{\mathcal E}M/M^1
$$
by Proposition \ref{irrecontainer} so that 
the maximal slope of the quotient $\mathcal E\otimes_{\mathcal E^\dagger}L^\dagger/M^\dagger$ 
is less than the maximal slope of $L$ if $L^\dagger \ne M^\dagger$. 
Since $L^\dagger$ is PBQ, Lemma \ref{subqut} (2) implies $L^\dagger/M^\dagger = 0$, hence $L^\dagger = M^\dagger$. 
The same holds for $N^\dagger$ and 
we have an identification $N^\dagger = M^\dagger = L^\dagger$ in $M/M^1$. Therefore, there is an isomorphism 
$g^\dagger : N^\dagger \rightarrow M^\dagger$ which makes the given diagram commutative. 
The uniqueness of $g^\dagger$ follows from the PBQ property of $N^\dagger$ by Lemma \ref{subqut} (2). 
}

\vspace*{3mm}

\noindent
{\sc Proof of Theorem \ref{loc}.} In the case (i) the nontrivial morphism $h$ is an isomorphism 
by Corollary \ref{localirr}. Then the assertion easily follows from Theorem \ref{locrel}.
\hspace*{\fill} $\Box$

\section{Saturated overconvergent $F$-isocrystals}\label{satur}

In this section we introduce saturated overconvergent $F$-isocrystals. 

\subsection{Saturated overconvergent $F$-isocrystals}\label{satglo} 
Let $X$ be a smooth connected scheme separated of finite type over $\mathrm{Spec}\, k$, 
and $\overline{X}$ a completion of $X$ over $\mathrm{Spec}\, k$ with the canonical open immersion $j_{X, \overline{X}} : X \rightarrow \overline{X}$. 
Let $\mathcal M^\dagger$ be an overconvergent $F$-isocrystal on $X/K$ with respect to Frobenius $\sigma$, 
$\mathcal M^\dagger_ U = j_{U, \overline{X}}^\dagger\mathcal M^\dagger$ the restriction of $\mathcal M^\dagger$ on $U/K$ 
as an overconvergent $F$-isocrystal, and $\mathcal M_U$ the convergent $F$-isocrystal on $U/K$ associated to $\mathcal M^\dagger$. 
For an affine open dense subscheme $U$ of $X$, let us take an affine smooth lift $\mathcal U = \mathrm{Spec}\, A_U$ over $\mathrm{Spec}\, R$ 
such that $\mathcal U \times_{\mathrm{Spec}\, R}\mathrm{Spec}\, k = U$, $\overline{\mathcal U}$ a completion of $\mathcal U$ over $\mathrm{Spec}\, R$, 
$\overline{U} = \overline{\mathcal U} \times_{\mathrm{Spec}\, R} \mathrm{Spec}\, k$ the reduction of $\overline{\mathcal U}$, 
and $\widehat{\overline{\mathcal U}}$ the $p$-adic completion of $\overline{\mathcal U}$. 
Such a lift $\mathcal U$ exists by \cite[Th\'eor\`eme 6]{El73}. 
Let us put $M^\dagger_U = \Gamma(]\overline{U}[_{\widehat{\overline{\mathcal U}}}, \mathcal M^\dagger_U)$, 
$M_U = \Gamma(]U[_{\widehat{\overline{\mathcal U}}}, \mathcal M_U)$ with the slope filtration $\{ M^i_U\}$ if $\mathcal M_U$ admits 
the slope filtration. 

\begin{definition}\label{satdef} An overconvergent $F$-isocrystal $\mathcal M^\dagger$ is saturated 
if the natural morphism $M^\dagger_U \rightarrow M_U/M_U^1$ is injective for any affine open dense subscheme $U \subset X$ 
such that $\mathcal M$ admits the slope filtration. 
\end{definition}

\begin{proposition}\label{injec} Let $\mathcal M^\dagger$ be an overconvergent $F$-isocrystal on $X/K$. 
Let $V$ be an affine open dense subscheme of $X$ such that the associated convergent $F$-isocrystal 
$\mathcal M_V$ to $\mathcal M_V^\dagger$ admits a slope filtration. 
Then $\mathcal M^\dagger$ is saturated if and only if 
$M^\dagger_V\rightarrow M_V/M^1_V$ is injective. 
\end{proposition}

\prf{Suppose that $M^\dagger_V\rightarrow M_V/M^1_V$ is injective. It is sufficient to prove, 
for any affine open dense subscheme $U$ of $V$, the natural morphism $M^\dagger_U = A_{U, K}^\dagger\otimes_{A^\dagger_{V, K}}M_V^\dagger \rightarrow M_U/M^1_U$ 
is injective since the top horizontal arrow of the commutative diagram 
$$
    \begin{array}{ccc} 
        M^\dagger_W &\rightarrow &M_{V \cap W}^\dagger \\
        \downarrow & &\downarrow \\
        M_W/M^1_W &\rightarrow &M_{V \cap W}/M_{V \cap W}^1
        \end{array}
$$
is injective for any affine open dense subscheme $W$ of $X$. 

\begin{lemma}\label{injob} Suppose that $X$ is affine with an affine smooth lift $\mathcal X = \mathrm{Spec}\, A_X$ over $\mathrm{Spec}\, R$ 
such that $A_X^\dagger$ admits a $q$-Frobenius $\varphi$ which is compatible with the $q$-Frobenius $\sigma$ on $K$. 
Let $\mathcal M^\dagger$ be an overconvergent $F$-isocrystal on $X/K$ such that 
the associated convergent $F$-isocrystal $\mathcal M$ to $\mathcal M^\dagger$ admits the slope filtration. 
Then there exists an overconvergent $F$-isocrystal $\mathcal L^\dagger$ on $X/K$ with an isomorphism 
$$
        \Gamma(]\overline{X}[_{\widehat{\overline{\mathcal X}}}, \mathcal L^\dagger) \cong \mathrm{Im}(M^\dagger \rightarrow M/M^1)
$$
of $A_{X, K}^\dagger$-modules such that the isomorphism commutes with connections and Frobenius. 
\end{lemma}

\prf{Since $A^\dagger_{X, K}$ is Noetherian \cite[Theorem]{Fu69}, the image 
$L^\dagger = \mathrm{Im}(M^\dagger \rightarrow M/M^1)$ is a finitely generated $A^\dagger_{X, K}$-module. 
Since the $q$-Frobenius $\varphi$ on $A_{X, K}^\dagger$ is flat (see Appendix \ref{FrobA}) 
and 
$\mathrm{id}\otimes\varphi : A^\dagger_{X, K}\otimes_{A^\dagger_{X, K}, \varphi} \widehat{A}_{X, K} \rightarrow \widehat{A}_{X, K}$ 
is an isomorphism, 
the natural morphism $A^\dagger_{X, K}\otimes_{A^\dagger_{X, K}, \varphi}L^\dagger \rightarrow 
\widehat{A}_{X, K}\otimes_{\widehat{A}_{X, K}, \varphi}M/M^1$ is injective 
where $A^\dagger_{X, K}\otimes_{A^\dagger_{X, K}, \varphi}-$ (resp. $\widehat{A}_{X, K}\otimes_{\widehat{A}_{X, K}, \varphi}-$) means 
the extension by $\varphi : A^\dagger_{X, K} \rightarrow A^\dagger_{X, K}$ 
(resp. $\varphi : \widehat{A}_{X, K} \rightarrow \widehat{A}_{X, K}$). 
By the compatibility of the $A_{X, K}^\dagger$-linear homomorphism $M^\dagger \rightarrow M/M^1$ 
with the connections and the Frobenius the integrable $K$-connection $\nabla_{M^\dagger}$ 
and Frobenius $\varphi_{M^\dagger}$ induce an integrable $K$-connection $\nabla_{L^\dagger} : L^\dagger \rightarrow L^\dagger\otimes_{A_X}\Omega_{A_X/R}^1$ 
and an isomorphism $\varphi_{L^\dagger} : \varphi^\ast L^\dagger \rightarrow L^\dagger$ 
such that $\varphi_{L^\dagger}$ is horizontal. 
Hence the connection $\nabla_{L^\dagger}$ is overconvergent 
and the sheafification $\mathcal L^\dagger$ of $L^\dagger$ is our desired overconvergent $F$-isocrystal on $X/K$ by \cite[Th\'eor\`eme 2.5.7]{Be96} 
(see Theorem \ref{affov}). 
}

\vspace*{3mm}

Let us continue the proof of Proposition \ref{injec}. 
Suppose $U$ is an open dense subscheme of $V$. Our claim is that, if $\mathcal N^\dagger_U$ is an overconvergent $F$-isocrystal on $U/K$ 
such that 
$$
     \Gamma(]\overline{X}[_{\widehat{\overline{\mathcal X}}}, \mathcal N^\dagger_U) \cong 
     \mathrm{Ker}(M_U^\dagger \rightarrow  M_U/M_U^1), 
$$
then $\mathcal N^\dagger_U = 0$. Here the right hand side is a kernel of the natural homomorphism and it 
determines an overconvergent $F$-isocrystal on $U/K$ by Lemma \ref{injob} since the category of overconvergent $F$-isocrystals is Abelian. 
Since $\mathcal M_U^\dagger$ is a restriction of $\mathcal M^\dagger$ on $U$ and 
$\mathcal N_U^\dagger$ 
is a subobject of $\mathcal M_U^\dagger$, there exists an overconvergent $F$-isocrystal $\mathcal N^\dagger$ on $X/K$ 
whose restriction on $U$ is $\mathcal N^\dagger_U$ by \cite[Proposition 5.3.1]{Ke07}. 
By our construction $\mathcal N_V \rightarrow \mathcal M_V/\mathcal M_V^1$ is a zero map by the full faithfulness 
of the restriction functor from the category of convergent $F$-isocrystals on $V$ to that on $U$ \cite[Theorem 5.2.1]{Ke07} \cite[Theorem 4.2.1]{Ke08b}. 
On the other hand the composition 
$$
N_V^\dagger = \Gamma(]\overline{X}[_{\widehat{\overline{\mathcal X}}}, \mathcal N^\dagger_V) \subset M_V^\dagger \rightarrow M_V/M_V^1
$$
is injective. Hence $\mathcal N_U^\dagger = 0$. 
}

\begin{corollary}\label{injopen} Let $U$ be an open dense subscheme of $X$. 
An overconvergent $F$-isocrystal $\mathcal M^\dagger$ on $X/K$ is saturated 
if and only if so is the restriction $\mathcal M^\dagger_U$ on $U/K$. 
\end{corollary}

\begin{corollary}\label{irrsat} An irreducible overconvergent $F$-isocrystal on $X/K$ is saturated. 
\end{corollary}

\prf{If $\mathcal M^\dagger$ is irreducible, then so is  $\mathcal M^\dagger_U$ for arbitrary affine open dense subscheme $U$ of $X$ by 
\cite[Proposition 5.3.1]{Ke07}. $\mathcal M^\dagger_U$ is saturated by Proposition \ref{injec} and Lemma \ref{injob} 
because one can find an affine open dense subscheme $U$ of $X$ on which the hypothesis in Lemma \ref{injob} holds. 
}

\begin{proposition}\label{inj1} Let $(\mathcal M^\dagger, F_{\mathcal M^\dagger})$ be an overconvergent $F$-isocrystal on $X/K$. 
\begin{enumerate}
\item Let $f : Y \rightarrow X$ be a finite etale morphism 
of connected schemes over $\mathrm{Spec}\, k$. 
Then $\mathcal M^\dagger$ is saturated if and only if so is the inverse image $f^\ast\mathcal M^\dagger$ 
as an overconvergent $F$-isocrystal on $Y/K$. 
\item Suppose $L$ (resp. $R_L$) is a finite extension of $K$ with a residue field $l$ (resp. the integer ring of $L$), 
$\sigma_L : L \rightarrow L$ an extension of the $n$-th power $\sigma^n$ of the $q$-Frobenius $\sigma$ on $K$. 
Let us put $X_l = X \times_{\mathrm{Spec}\, k}\mathrm{Spec}\, l$ to be the base extension of $X$, 
and $\mathcal M^\dagger_{X_l/L}$ the induced overconvergent $F$-isocrystal on $X_l/L$ with respect to $\sigma_L$ 
by the inverse image of $(\mathcal M^\dagger, F_{\mathcal M^\dagger}^n)$. 
Then $\mathcal M^\dagger$ is saturated if and only if so is $\mathcal M^\dagger_{X_l/L}$ on each connected component of $X_l$. 
\end{enumerate}
\end{proposition}

\prf{We may assume that $X$ is affine and $\mathcal M$ admits the slope filtration by Corollary \ref{injopen}.
Take an affine smooth lift $\mathrm{Spec}\, A_X$ of $X$ over $\mathrm{Spec}\, R$.  

(1) Since $Y$ is finite etale over $X$, there exists a finite $A_X$-algebra $A_Y$ such that 
the weak completion $A_Y^\dagger$ is finite etale over $A_X^\dagger$ by Jacobian criterion of etaleness. 
The assertion follows from the fact that $A_Y^\dagger$ is faithfully flat over $A_X^\dagger$ and 
the natural morphism $A_{Y, K}^\dagger \otimes_{A_{X, K}^\dagger}\widehat{A}_{X, K} \rightarrow \widehat{A}_{Y, K}$ is an isomorphism. 

(2) The assertion follows from Proposition \ref{injec} and the fact that the morphism 
$A^\dagger_{X, K} \rightarrow A_{X_l, L}^\dagger \cong L \otimes_KA^\dagger_{X, K}$
is faithfully flat and the natural morphism $A_{X_l, L}^\dagger\otimes_{A^\dagger_{X, K}} \widehat{A}_{X, K} \rightarrow \widehat{A}_{X_l, L}$
is an isomorphism.}

\begin{remark}
In the previous version of this paper we define an saturated overconvergent $F$-isocrystal on a curve $C$ 
if the natural morphism $M^\dagger \rightarrow M_\eta/M_\eta^1$ is injective. 
Here $M_\eta$ is the generic $\varphi$-$\nabla$-module over $E_\eta$ (see section \ref{set} for the notation). 
It is equivalent to that of Definition \ref{satdef}. Indeed, the natural morphism $M_U/M_U^1 \rightarrow M_\eta/M_\eta^1$ is injective 
for any affine dense subscheme $U$ of $C$. 
\end{remark}

\subsection{Saturation}

\begin{proposition}\label{satov} Let $\mathcal M^\dagger$ be an overconvergent $F$-isocrystal on $X/K$. 
Then there exists a saturated overconvergent $F$-isocrystal $\mathcal L^\dagger$ on $X/K$ with a surjective 
morphism $\theta^\dagger : \mathcal M^\dagger \rightarrow \mathcal L^\dagger$ as overconvergent $F$-isocrystals 
such that the induced morphism $\mathcal M/\mathcal M^1 \rightarrow \mathcal L/\mathcal L^1$ between the maximal slope quotient 
is an isomorphism as convergent $F$-isocrystals. 
If a saturated overconvergent $F$-isocrystal $(\mathcal L')^\dagger$ on $X/K$ satisfies the similar properties, then 
there exists a unique isomorphism $\xi : (\mathcal L')^\dagger \rightarrow \mathcal L^\dagger$ of overconvergent $F$-isocrystals. 
\end{proposition} 

\prf{For an affine open dense subscheme $U$ of $X$ in Lemma \ref{injob}, let $\mathcal L^\dagger_U$ 
be the overconvergent $F$-isocrystal on $U/K$ in Lemma \ref{injob}. 
Since $\mathcal L^\dagger_U$ is a quotient of $\mathcal M^\dagger$ and $L_U$ generates $M_U/M_U^1$ 
as an $\widehat{A}_{U, K}$-module, $\mathcal L^\dagger_U$ is saturated and the induced morphism 
$\mathcal M_U/\mathcal M^1_U \rightarrow \mathcal L_U/\mathcal L_U^1$ is an isomorphism. 

If one has another $(\mathcal L')^\dagger$ on $U/K$ with the same properties, then 
$$
      (L')^\dagger_U \rightarrow L'_U/(L')^1_U \cong M/M^1 \cong L_U/L^1_U \leftarrow L_U^\dagger
$$
induces a bijection of $A^\dagger_{U, K}$-modules between $L_U^\dagger$ and $(L')_U^\dagger$ 
since there exist compatible surjections from $M_U^\dagger$ to both $(L')^\dagger_U$ and $L_U^\dagger$. 
Hence we have a canonical isomorphism $\xi_U : (\mathcal L')^\dagger_U \rightarrow \mathcal L^\dagger_U$ 
of overconvergent $F$-isocrystals on $U/K$. 
The glueing method works by the canonical uniqueness of $\mathcal L^\dagger_U$ and Corollary \ref{injopen}, 
and we obtain the desired overconvergent $F$-isocrystal $\mathcal L^\dagger$ on $X/K$. 
}

\begin{definition}\label{saturatover} For an overconvergent $F$-isocrystal $\mathcal M^\dagger$ on $X/K$, 
we define the saturation $\mathcal M^{\mathrm{sat}, \dagger}$ of $\mathcal M^\dagger$ by the 
overconvergent $F$-isocrystal $\mathcal L^\dagger$ on $X/K$ in Proposition \ref{satov}. 
\end{definition}

\section{In the case of curves}\label{curve} 

In this section we prove the minimal slope conjecture in the case of curves. 
Let $C$ be a smooth connected curve over $\mathrm{Spec}\, k$, and $\overline{C}$ a smooth completion of $C$. 

\subsection{Ranks of the maximal slope quotients}

\begin{proposition}\label{ranksheaf} 
Let $\mathcal M$ be an isoclinic convergent $F$-isocrystal on $C/K$, and $\mathcal N^\dagger$ be an overconvergent $F$-isocrystals 
on $C/K$ such that the associated convergent $F$-isocrystal $\mathcal N$ on $C/K$ 
admits the slope filtration $\{ \mathcal N^i\}$. 
Suppose there exists a surjective morphism $\mathcal N \rightarrow \mathcal M$ of convergent $F$-isocrystal 
such that the induced morphism $\Gamma(]\overline{C}[_{\widehat{\overline{\mathcal C}}}, j_U^\dagger \mathcal N^\dagger) 
\rightarrow \Gamma(]U[_{\widehat{\overline{\mathcal C}}}, \mathcal M)$ of global sections is injective 
for an affine open dense subscheme $U$ of $C$. Then the maximal slope of $\mathcal N$ coincides with 
the slope of $\mathcal M$ and the equality below holds:
$$
       \mathrm{rank}\, \mathcal N/\mathcal N^1 = \mathrm{rank}\, \mathcal M. 
$$
\end{proposition}

\prf{Fix a notation as in Section \ref{set}. 
If $C = \overline{C}$, the assertion is trivial. 
Hence we may assume that $C$ is affine with a smooth lift $\mathcal C = \mathrm{Spec}\, A_C$. 
First we prove the coincidence between the slope of $\mathcal M$ and the maximal slope of $\mathcal N$. 
Let $\mathcal P^\dagger$ be the maximally PBQ overconvergent $F$-subisocrystal of $\mathcal N^\dagger$. 
Then $\mathcal P$ admits the slope filtration by \cite[Proposition 2.7]{Ts19}.  
Since the induced morphism $\mathcal P \rightarrow \mathcal M$ is nontrivial, 
the maximal slope of $\mathcal P$ must be equal to the slope of $\mathcal M$  by Lemma \ref{subqutsh}. This yields 
the maximal slope of $\mathcal N$ is equal to the slope of $\mathcal M$. 

By the hypothesis of the surjectivity of $\mathcal N \rightarrow \mathcal M$ the inequality 
$\mathrm{rank}\, \mathcal N/\mathcal N^1 \geq \mathrm{rank}\, \mathcal M$ 
always holds. Suppose $\mathrm{rank}\, \mathcal N/\mathcal N^1 > \mathrm{rank}\, \mathcal M$ holds. 
Let us define a nontrivial convergent $F$-isocrystal $\mathcal H$ on $C/K$ by 
$$
      \mathcal H = \mathrm{Ker}(\mathcal N \rightarrow \mathcal M)
$$
and $\mathcal L$ the maximally PBQ convergent $F$-subisocrystal of $\mathcal H$. 
The convergent version of PBQ filtrations also exists similarly to Theorems \ref{genPBQfil} and \ref{globPBQfil}. 
By our hypothesis the maximal slope of $\mathcal H$ and hence $\mathcal L$ is that of $\mathcal M$. 
We will construct the maximally PBQ overconvergent $F$-subisocrystal of $\mathcal L^\dagger$ 
such that the convergent $F$-isocrystal on $C/K$ associated to $\mathcal L^\dagger$ is the above $\mathcal L$. 

Let $\alpha_1, \alpha_2, \cdots, \alpha_s$ be all distinct closed points of $\overline{C} \setminus C$. 
Let us put $N^\dagger_C = 
\Gamma(]\overline{C}[_{\widehat{\overline{\mathcal C}}}, \mathcal N^\dagger)$, 
$M_C = \Gamma(]C[_{\widehat{\overline{\mathcal C}}}, \mathcal M)$ and $L_C = \Gamma(]C[_{\widehat{\overline{\mathcal C}}}, \mathcal L)$. 
For any closed point $\alpha_i$ we put 
$$
     H_{\alpha_i}^\dagger = \mathrm{Ker}(\mathcal E_{\alpha_i}^\dagger\otimes_{A_{C, K}^\dagger}N_C^\dagger \rightarrow \mathcal E_{\alpha_i} \otimes_{\widehat{A}_{C, K}}M_C). 
$$
Then $H_{\alpha_i}^\dagger$ is a $\varphi$-$\nabla$-module over $\mathcal E_{\alpha_i}^\dagger$. 
Let us take the maximally PBQ $\varphi$-$\nabla$-module $L_{\alpha_i}^\dagger$ of $H_{\alpha_i}^\dagger$ over $\mathcal E_{\alpha_i}^\dagger$. 
Then our claim is as follows.

\vspace*{2mm}

\noindent
{\bf Claim.} {\it Let us define $H_{\alpha_i} = \mathcal E_{\alpha_i}\otimes_{\mathcal E_i^\dagger}H_{\alpha_i}^\dagger$ 
and $L_{\alpha_i} = \mathcal E_{\alpha_i}\otimes_{\mathcal E_i^\dagger}L_{\alpha_i}^\dagger$ 
for $1 \leq i \leq s$. 
\begin{enumerate}
\item $\mathrm{dim}_{\mathcal E_{\alpha_i}}\, H_{\alpha_i}/H_{\alpha_i}^1 = \mathrm{rank}\, \mathcal H/\mathcal H^1$. 
\item The natural commutative diagram 
$$
    \begin{array}{ccccccccc}
         0 &\rightarrow &H^\dagger_{\alpha_i} &\rightarrow &\mathcal E_{\alpha_i}^\dagger\otimes_{A_{C, K}^\dagger}N_C^\dagger 
         &\rightarrow &\mathcal E_{\alpha_i} \otimes_{\widehat{A}_{C, K}}M_C  \\
         & &\downarrow & &\downarrow & &\hspace*{3mm} \downarrow = \\
         0 &\rightarrow &\mathcal E_{\alpha_i} \otimes_{\widehat{A}_{C, K}}H_C &\rightarrow &\mathcal E_{\alpha_i} \otimes_{A_{C, K}^\dagger}N_C^\dagger 
         &\rightarrow &\mathcal E_{\alpha_i} \otimes_{\widehat{A}_{C, K}}M_C &\rightarrow &0
         \end{array}
$$
with exact rows induces an injection
$$
           L_{\alpha_i} 
           \rightarrow \mathcal E_{\alpha_i} \otimes_{\widehat{A}_{C, K}}L_C
$$
of $\varphi$-$\nabla$-modules over $\mathcal E_{\alpha_i}$.
\item The injection $L_{\alpha_i} 
           \rightarrow \mathcal E_{\alpha_i} \otimes_{\widehat{A}_{C, K}}L_C$ is an isomorphism.
\end{enumerate}
}

\vspace*{2mm}

\prf{(1) Since 
$N_{\alpha_i} = \mathcal E_{\alpha_i}\otimes_{A_{C, K}^\dagger}N_C^\dagger \rightarrow \mathcal E_{\alpha_i} \otimes_{\widehat{A}_{C, K}}M_C$ 
is surjective, we have 
$$
      \mathrm{dim}_{\mathcal E_{\alpha_i}}H_{\alpha_i}/H_{\alpha_i}^1 = 
      \mathrm{dim}_{\mathcal E_{\alpha_i}}N_{\alpha_i}/N_{\alpha_i}^1 - \mathrm{rank}\, \mathcal M = \mathrm{rank}\, \mathcal H/\mathcal H^1
$$
applying Proposition \ref{irrecontainer} to $N_{\alpha_i}^\dagger/H_{\alpha_i}^\dagger \subset \mathcal E_{\alpha_i} \otimes_{\widehat{A}_{C, K}}M_C$. 

(2) Since $\mathcal E_{\alpha_i}^\dagger \rightarrow \mathcal E_{\alpha_i}$ is a field extension, the natural morphism   
$H_{\alpha_i} 
     \rightarrow \mathcal E_{\alpha_i} \otimes_{\widehat{A}_{C, K}}H_C$ is injective. 
The image of $L_{\alpha_i}$ is included in $\mathcal E_{\alpha_i} \otimes_{\widehat{A}_{C, K}}L_C$ by 
Lemma \ref{subqut} (2) and Theorem \ref{coinPBQfil} (2). Indeed, the maps from $L_{\alpha_i}$ 
to the PBQ subquotients induced by the PBQ filtration of $\mathcal E_{\alpha_i} \otimes_{\widehat{A}_{C, K}}H_C$ 
is a zero map if the maximal slopes are different. 

(3) Note that the slope filtration is strict for morphisms of $F$-isocrystals and it is compatible of localization. 
Hence there is a sequence of isomorphisms 
$$
     L_{\alpha_i}/L_{\alpha_i}^1 \cong H_{\alpha_i}/H_{\alpha_i}^1 \cong \mathcal E_{\alpha_i} \otimes_{\widehat{A}_{C, K}}H_C/\mathcal E_{\alpha_i} \otimes_{\widehat{A}_{C, K}}H_C^1 
     \cong \mathcal E_{\alpha_i} \otimes_{\widehat{A}_{C, K}}L_C/\mathcal E_{\alpha_i} \otimes_{\widehat{A}_{C, K}}L_C^1
$$
by (1) and the global-local compatibility of PBQ filtrations by Theorems \ref{coinPBQfil} (2) and \ref{globPBQfil}. 
Since both $L_{\alpha_i}/L_{\alpha_i}^1$ and $\mathcal E_{\alpha_i} \otimes_{\widehat{A}_{C, K}}L_C/\mathcal E_{\alpha_i} \otimes_{\widehat{A}_{C, K}}L_C^1$ 
are PBQ with isomorphic maximally slope quotients, the injection $L_{\alpha_i} \rightarrow \mathcal E_{\alpha_i} \otimes_{\widehat{A}_{C, K}}L_C$ is surjective 
by Lemma \ref{subqut} (2). 
Note that one can not directly apply Theorem \ref{loc} by lack of saturation. 
}

\vspace*{3mm}

Now applying Lemma \ref{extext} or \cite[Remark 5.10]{Ke16} to the compatible data $L_C, L_{\alpha_i}^\dagger\, (i =1, \cdots, s)$ 
in the overconvergent $F$-isocrystal $\mathcal N^\dagger$, we have an 
overconvergent $F$-isocrystal $\mathcal L^\dagger$ on $C/K$ such that $\mathcal L^\dagger$ is a 
subobject of $\mathcal N^\dagger$ and that 
the convergent 
$F$-isocrystal $\mathcal L$ on $C/K$ associated to $\mathcal L^\dagger$ is $L_C$. 
Then we have an injection 
$$
     L^\dagger_C = \Gamma(]\overline{C}[_{\widehat{\overline{\mathcal C}}}, \mathcal L^\dagger) \rightarrow N^\dagger_C \rightarrow M_C. 
$$
But the associated morphism $L_C \rightarrow M_C$ is a zero map, and it is a contradiction. 
Therefore, the equality $\mathrm{rank}\, \mathcal N/\mathcal N^1 = \mathrm{rank}\, \mathcal M$ holds. 
This completes a proof of Proposition \ref{ranksheaf}. 
}

\begin{proposition}\label{globalirr} Let $\mathcal M^\dagger$ be an overconvergent $F$-isocrystal on $C/K$. 
Then the following conditions are equivalent.
\begin{list}{}{}
\item[\mbox{\rm (i)}] $\mathcal M^\dagger$ is irreducible. 
\item[\mbox{\rm (ii)}] $\mathcal M^\dagger$ is PBQ and saturated and $\mathcal M_U/\mathcal M_U^1$ is irreducible 
for any open dense subscheme $U$ of $C$ such that $\mathcal M_U$ admits the slope filtration. 
\end{list}
\end{proposition}

\prf{The similar proof works with Corollary \ref{localirr} by 
Corollary \ref{PBQirr}, Lemma \ref{subqutsh}, Corollary \ref{injopen} and Proposition \ref{ranksheaf}. 
}

\begin{example}\label{ntext} Let $0 \rightarrow \mathcal L^\dagger \rightarrow \mathcal M^\dagger \rightarrow \mathcal N^\dagger \rightarrow 0$ 
be an exact sequence of overconvergent $F$-isocrystals on $C/K$. 
Suppose $\mathcal L^\dagger$ and $\mathcal N^\dagger$ are irreducible with the maximal slopes $s_L$ and $s_N$, respectively. 
Then the following hold.
\begin{enumerate}
\item If $s_L < s_N$, then $\mathcal M^\dagger$ is not saturated. 
$\mathcal M^\dagger$ is a nontrivial extension if and only if $\mathcal M^\dagger$ is PBQ.
\item If $s_L=s_N$, then $\mathcal M^\dagger$ is PBQ and saturated, but $\mathcal M/\mathcal M^1$ is not irreducible. 
$\mathcal M^\dagger$ is a nontrivial extension if and only if so is $\mathcal M/\mathcal M^1$. 
\item If $s_L > s_N$, then the exact sequence gives the PBQ filtration of $\mathcal M^\dagger$. 
$\mathcal M^\dagger$ is a nontrivial extension if and only if $\mathcal M^\dagger$ is saturated. 
Moreover, if furthermore $\mathcal M^\dagger$ is of rank $2$, then the sequence is split by \cite[Theorem 4.2.1]{Ke08b} since 
the induced sequence as convergent $F$-isocrystals on $C/K$ is split. 
\end{enumerate}
\end{example}
 
\subsection{The minimal slope conjecture on curves}

\begin{theorem}\label{eigen} Let $\mathcal M^\dagger$ (resp. $\mathcal N^\dagger$) be an overconvergent $F$-isocrystal on $C/K$ 
such that $\mathcal M$ (resp. $\mathcal N$) admits the Frobenius filtration, and  
$0 = \mathcal P_0^\dagger \subsetneq \mathcal P_1^\dagger \subsetneq \cdots \subsetneq \mathcal P_r^\dagger = \mathcal M^\dagger$ 
(resp. $0 = \mathcal Q_0^\dagger \subsetneq \mathcal Q_1^\dagger \subsetneq \cdots \subsetneq \mathcal Q_s^\dagger = \mathcal N^\dagger$) 
the PBQ filtration of $\mathcal M^\dagger$ (resp. $\mathcal N^\dagger$), 
and $\mathcal M/\mathcal M^b = \oplus_i\, \mathcal P_i/\mathcal P_i^1$ 
(resp. $\mathcal N/\mathcal N^b = \oplus_j\, \mathcal Q_j/\mathcal Q_j^1$) 
the maximally bounded quotient of $\mathcal M$ (resp. $\mathcal N$) (Proposition \ref{globalbdqt}). 
If there exists an isomorphism 
$h : \mathcal N/\mathcal N^b \rightarrow \mathcal M_\eta/\mathcal M^b$ of overconvergent $F$-isocrystals, then $r = s$ and 
there is a unique isomorphism $g^\dagger : \oplus_{i=1}^r\mathcal Q_i^{\mathrm{sat}, \dagger} \rightarrow \oplus_{i=1}^r\mathcal P_i^{\mathrm{sat}, \dagger}$ 
such that the induced diagram 
$$
       \begin{array}{ccc} 
            \oplus_{i=1}^r \mathcal Q_i^{\mathrm{sat}} 
             &\overset{g}{\rightarrow} &\oplus_{i=1}^r \mathcal P_i^{\mathrm{sat}} \\
             \downarrow & &\downarrow \\
             \mathcal N/\mathcal N^b  &\underset{h}{\rightarrow} &\mathcal M/\mathcal M^b
     \end{array}
$$
is commutative. 
\end{theorem} 

\prf{The proof is similar to that of the local case, Theorem \ref{locrel}, 
by applying Proposition \ref{ranksheaf} and Lemma \ref{subqutsh} (2) instead of Proposition \ref{irrecontainer} 
and Lemma \ref{subqut} (2), respectively. 
}

\begin{corollary}\label{kedcur} Let $\mathcal M^\dagger$ and $\mathcal N^\dagger$ be overconvergent $F$-isocrystals on $C/K$ 
such that the associated convergent $F$-isocrystal $\mathcal M$ and $\mathcal N$ admit the slope filtrations, 
and $h : \mathcal N/\mathcal N^1 \rightarrow \mathcal M/\mathcal M^1$ a morphism between 
the maximal slope quotients as convergent $F$-isocrystals. 
Suppose either
\begin{list}{}{}
\item[\mbox{\rm (i)}] $\mathcal M^\dagger$ and $\mathcal N^\dagger$ are irreducible and $h$ is nontrivial, or 
\item[\mbox{\rm (ii)}] $\mathcal M^\dagger$ and $\mathcal N^\dagger$ are PBQ and saturated and $h$ is an isomorphism. 
\end{list}
Then there is a unique isomorphism  
$$
      g^\dagger : \mathcal N^\dagger \rightarrow \mathcal M^\dagger
$$
of overconvergent $F$-isocrystals on $C/K$ such that 
the induced morphism $\mathcal N/\mathcal N^1 \rightarrow \mathcal M/\mathcal M^1$ by $g^\dagger$ 
coincides with the given $h$. 
\end{corollary}

\prf{In the case (i) the nontrivial morphism $h$ is an isomorphism 
by Proposition \ref{globalirr}. The assertion follows from Theorem \ref{eigen}. 
}

\begin{remark}\label{rmmk}
\begin{enumerate}
\item Suppose $C$ is projective. Then $\mathcal M^\dagger = \mathcal M$. 
If furthermore $\mathcal M$ is irreducible and admits a slope filtration, then $\mathcal M^1 = 0$. 
Therefore the minimal slope conjecture is trivially true in this case. 
\item The hypothesis that $\mathcal M$ and $\mathcal N$ admit slope filtrations in Theorem \ref{eigen} and Corollary \ref{kedcur} 
is redundant. It is enough to assume that there exists a convergent $F$-subisocrystal $\mathcal M^1$ 
such that $\mathcal M/\mathcal M^1$ is isoclinic and all slopes of $\mathcal M^1$ at any point of $C$ is less than the slope of $\mathcal M/\mathcal M^1$, 
and the same for $\mathcal N$. However, if one takes a sufficiently small open dense subscheme $U$ of $C$, 
then both $\mathcal M$ and $\mathcal N$ admit slope filtrations \cite[Theorem 3.1.2, Corollary 4.2]{Ke16}, \cite[Proposition 2.2, Corollary 2.6]{Ts19}. 
The existence of $g^\dagger$ follows from the full faithfulness of the restriction functor of overconvergent $F$-isocrystals \cite[Theorem 6.3.1]{Ts02}. 
\end{enumerate}
\end{remark}

\subsection{Lefschetz condition}\label{LC}

Let $X$ be a scheme separated of finite type over $\mathrm{Spec}\, k$, 
$\alpha = \mathrm{Spec}\, k_\alpha$ a closed point of $X$ with the closed immersion $i_{\alpha, X} : \alpha \rightarrow X$, 
and $\mathcal M^\dagger$ an irreducible overconvergent $F$-isocrystal on $X/K$. 
Let us consider the following condition (LC), called Lefschetz condition, for $(X, \alpha, \mathcal M^\dagger)$ :
\begin{flushleft}
\begin{tabular}{ll}
\hspace*{3mm} (LC) &\begin{tabular}{l} There exist a smooth curve $C_\alpha$ over $\mathrm{Spec}\, k$ and a morphism 
$i_{C_\alpha, X} : C_\alpha \rightarrow X$ over \\ 
$\mathrm{Spec}\, k$ such that there exists a $k_\alpha$-rational point of $C_\alpha$ 
which maps to $\alpha$ by $i_{C_\alpha, X}$ \\
(say $C_\alpha$ is passing at $\alpha$ and denote the $k_\alpha$-rational point in $C_\alpha$ also by $\alpha$) \\
and that the restriction $i_{C_\alpha, X}^\ast\mathcal M^\dagger$ on $C_\alpha$ is irreducible. 
\end{tabular}
\end{tabular} 
\end{flushleft}
The question on an existence of such curves for $\mathcal M^\dagger$ 
has an affirmative answer more generally in the case where $k$ is finite and $\mathcal M^\dagger$ 
is an overconvergent $\overline{\mathbb Q}_p$-$F$-isocrystal 
in \cite[Theorem 0.3]{AE16} (see Theorem \ref{Lef}). 
It will be used in order to prove our main theorem (Theorem \ref{rlt2}). See \cite[Conjecture 5.19]{Ke16} 
for the detail of the problem. The next proposition follows from Theorem \ref{eigen}. 

\begin{proposition}\label{LCI} With the notation as above, let $\mathcal M^\dagger$ and $\mathcal N^\dagger$ 
be overconvergent $F$-isocrystals on $X/K$ such that the associated convergent $F$-isocrystals $\mathcal M$ and $\mathcal N$ 
admit the slope filtrations, and $h : \mathcal N/\mathcal N^1 \rightarrow \mathcal M/\mathcal M^1$ a nontrivial 
morphism as convergent $F$-isocrystals. 
Suppose that $\mathcal M^\dagger$ (resp. $\mathcal N^\dagger$) is irreducible and the triplet $(X, \alpha, \mathcal M^\dagger)$ 
(resp. $(X, \alpha, \mathcal N^\dagger)$) satisfies the condition \mbox{\rm (LC)} with respect to a morphism 
$i_{C_{\alpha, \mathcal M}, X} : C_{\alpha, \mathcal M} \rightarrow X$ 
(resp. $i_{C_{\alpha, \mathcal N}, X} : C_{\alpha, \mathcal N} \rightarrow X$). 
Let $\mathcal Q_{C_{\alpha, \mathcal M}}^\dagger$ (resp. $\mathcal P_{C_{\alpha, \mathcal N}}^\dagger$) 
be the maximally PBQ overconvergent $F$-subisocrystal 
of $i_{C_{\alpha, \mathcal M}, X}^\ast\mathcal N^\dagger$ (resp. $i_{C_{\alpha, \mathcal N}, X}^\ast\mathcal M^\dagger$). 
Then there exists a unique surjective (resp. injective) morphism 
$$
    \begin{array}{ll}
        &g_{C_{\alpha, \mathcal M}}^\dagger : \mathcal Q_{C_{\alpha, \mathcal M}}^\dagger 
        \rightarrow i_{C_{\alpha, \mathcal M}, X}^\ast\mathcal M^\dagger \\
        (\mathrm{resp.} &g_{C_{\alpha, \mathcal N}}^\dagger : i_{C_{\alpha, \mathcal N}, X}^\ast\mathcal N^\dagger 
        \rightarrow \mathcal P_{C_{\alpha, \mathcal N}}^\dagger)
        \end{array}
$$
of overconvergent $F$-isocrystals on $C_{\alpha, \mathcal M}/K$ 
(resp. $C_{\alpha, \mathcal N}/K$) such that the induced morphism between the maximal slope quotients 
from $g^\dagger_{C_{\alpha, \mathcal M}}$ (resp. $g^\dagger_{C_{\alpha, \mathcal N}}$) 
is the given $i_{C_{\alpha, \mathcal M}, X}^\ast(h)$ (resp. $i_{C_{\alpha, \mathcal N}, X}^\ast(h)$). 
\end{proposition}

\prf{Suppose $\mathcal M^\dagger$ (resp. $\mathcal N^\dagger$) is irreducible. Then $\mathcal M/\mathcal M^1$ 
(resp. $\mathcal N/\mathcal N^1$) is irreducible by Proposition \ref{globalirr} so that $h$ is surjective (resp. injective). 
Hence there exists a unique surjective (resp. injective) morphism $g_{C_{\alpha, \mathcal M}}^\dagger$ 
(resp. $g_{C_{\alpha, \mathcal N}}^\dagger$) which is compatible with 
the given $i_{C_{\alpha, \mathcal M}, X}^\ast(h)$ (resp. $i_{C_{\alpha, \mathcal N}, X}^\ast(h)$) by Theorem \ref{eigen}. 
}

\begin{corollary}\label{isom} Suppose furthermore that both $\mathcal M^\dagger$ and $\mathcal N^\dagger$ are irreducible 
and that both $(X, \alpha, \mathcal M^\dagger)$ and $(X, \alpha, \mathcal N^\dagger)$ satisfy the condition \mbox{\rm (LC)}. 
If $K_\alpha$ is an unramified extension of $K$ with residue field $k_\alpha$, then there exists an isomorphism 
$$
         g_\alpha^\dagger : i_{\alpha, X}^\ast\mathcal N^\dagger \rightarrow  i_{\alpha, X}^\ast\mathcal M^\dagger
 $$
 of $F$-spaces over $K_\alpha$. 
\end{corollary}

\prf{Note that, if one can take the same curve $C_\alpha \rightarrow X$ for both $\mathcal M^\dagger$ 
and $\mathcal N^\dagger$, then the assertion follows just from Corollary \ref{kedcur} (1). 
Keep the notation as in Proposition \ref{LCI}. Both inequalities 
$\mathrm{rank}\, \mathcal M^\dagger \leq \mathrm{rank}\, 
\mathcal Q_{C_{\alpha, \mathcal M}}^\dagger \leq \mathrm{rank}\, \mathcal N^\dagger$ 
and $\mathrm{rank}\, \mathcal N^\dagger \leq \mathrm{rank}\, \mathcal P_{C_{\alpha, \mathcal N}}^\dagger 
\leq \mathrm{rank}\, \mathcal M^\dagger$ 
hold by Proposition \ref{LCI}. Hence $\mathrm{rank}\, \mathcal M^\dagger = \mathrm{rank}\, \mathcal N^\dagger$. 
We have an isomorphism 
$g_{C_{\alpha, \mathcal M}}^\dagger : i_{C_{\alpha, \mathcal M}, X}^\ast\mathcal N^\dagger 
\rightarrow i_{C_{\alpha, \mathcal M}, X}^\ast\mathcal N^\dagger$ by 
$$
     i_{C_{\alpha, \mathcal M}, X}^\ast\mathcal N^\dagger\, \overset{\supset}{\longleftarrow}\, \mathcal Q_{C_{\alpha, \mathcal M}}^\dagger\, 
     \overset{\cong}{\longrightarrow}\, i_{C_{\alpha, \mathcal M}, X}^\ast\mathcal M^\dagger 
$$
such that it induces $i_{C_{\alpha, \mathcal M}, X}^\ast(h)$. It induces an isomorphism 
$g_\alpha^\dagger : i_{\alpha, X}^\ast\mathcal N^\dagger \rightarrow  i_{\alpha, X}^\ast\mathcal M^\dagger$ as
$F$-spaces over $K_\alpha$. 
}

\begin{corollary}\label{eigeneta} With the hypothesis in Proposition \ref{LCI}, 
suppose furthermore that $\mathcal M^\dagger$ is irreducible and there is a closed point $\alpha$ of $X$ 
such that $(X, \alpha, \mathcal M^\dagger)$ satisfies the condition \mbox{\rm (LC)}. 
If $\mathcal N^\dagger$ is unit-root, then $\mathcal M^\dagger$ is unit-root and there exists a unique surjection   
$$
      g^\dagger : \mathcal N^\dagger \rightarrow \mathcal M^\dagger
$$
of overconvergent $F$-isocrystals on $X/K$ such that 
the induced morphism of converegent $F$-isocrystals from $g^\dagger$ coincides with the given $h$. 
\end{corollary}

\prf{By the hypothesis $i^\ast_{\alpha, X}\mathcal M^\dagger$ is unit-root by Proposition \ref{LCI}. Hence $\mathcal M^\dagger$ is unit-root (i.e., $\mathcal M^1 = 0$) and 
the assertion follows from the full faithfulness of the functor from the category of overconvergent $F$-isocrystals to that of convergent 
$F$-isocrystals \cite[Theorem 4.2.1]{Ke08b}.}

\section{The case of finite fields}\label{finite}

In this section we study the minimal slope conjecture on varieties of arbitrary dimension over finite fields. 
Let $\overline{\mathbb Q}_p$ be an algebraic closure of the field $\mathbb Q_p$ of $p$-adic numbers. 
Let $q$ be a positive power of $p$, $k$ the field of $q$ elements, 
$R = W(k)$ the ring of Witt vectors with coefficients in $k$, 
and $K$ the field of fractions of $R$ with the $q$-Frobenius $\sigma = \mathrm{id}_K$. 
We regard $K$ as a subfield of $\overline{\mathbb Q}_p$ 
and identify an algebraic closure $\overline{K}$ of $K$ with $\overline{\mathbb Q}_p$. 

\subsection{$F$-isocrystals with $\overline{\mathbb Q}_p$-structures} 
We recall the definition of $\overline{\mathbb Q}_p$-$F$-isocrystals which are introduced by Abe in \cite[1.4, 4.2.1]{Ab18b} 
(see also \cite[1.1]{AE16} and \cite[9.2]{Ke16}). 
Let $X$ be a scheme locally of finite type over $\mathrm{Spec}\, k$. 

\begin{definition}\label{defQpF} Let $L$ be a finite extension of $K$ in $\overline{\mathbb Q}_p$. 
An overconvergent $L$-$F_q$-isocrystal $\mathcal M^\dagger$ on $X$ 
is an overconvergent $F$-isocrystal $\mathcal M^\dagger$ 
on $X/K$ with respect to the $q$-Frobenius $\sigma$ such that $\mathcal M^\dagger$ is furnished with a $K$-algebra homomorphism 
$$
   L \rightarrow \mathrm{End}_{\footnotesize \mbox{\rm $F$-Isoc}}(\mathcal M^\dagger), 
$$
where 
$\mathrm{End}_{\footnotesize \mbox{\rm $F$-Isoc}}(\mathcal M^\dagger)$ is the $K$-space of endomorphisms 
of $\mathcal M^\dagger$ as overconvergent $F$-isocrystals on $X/K$. Note that the $q$ of $F_q$ indicates $q$-Frobenius. 
The $K$-algebra homomorphism $L \rightarrow \mathrm{End}_{\footnotesize \mbox{\rm $F$-Isoc}}(\mathcal M^\dagger)$ is called an $L$-structure. 
A morphism of overconvergent $L$-$F_q$-isocrystals on $X$ is a morphism of overconvergent $F$-isocrystals on $X/K$ 
which commutes with $L$-structures. The category of overconvergent $L$-$F_q$-isocrystals on $X$ 
is denoted by $F_q\mbox{\rm -Isoc}^\dagger(X)\otimes L$. 

We also define the notion of convergent $L$-$F_q$-isocrystals on $X$ similarly and denote their category by 
$F_q\mbox{\rm -Isoc}(X)\otimes L$. 
\end{definition}

\begin{example} 
\begin{enumerate}
\item An overconvergent $F$-isocrystal on $X/K$ is furnished with the natural $K$-structure 
since the category of overconvergent $F$-isocrystals on $X/K$ is $K$-linear by $\sigma = \mathrm{id}_K$. 
\item Let $L$ be an extension 
of $K$ in $\overline{\mathbb Q}_p$. We regard $L$ as an $F$-space over $K$ with respect to $\sigma$ by $F_L=\mathrm{id}_L$. 
Then, for an overconvergent $F$-iaocrystal $\mathcal M^\dagger$ on $X/K$, 
the $j^\dagger_X\mathcal O_{]\overline{X}[}\otimes_KL$-module $\mathcal M^\dagger\otimes_KL$ is 
furnished with an $L$-structure which is defined by 
$$
      a = \mathrm{id}_{\mathcal M^\dagger}\otimes a\mathrm{id}_L : \mathcal M^\dagger\otimes_KL \rightarrow \mathcal M^\dagger\otimes_KL
$$
for $a \in L$. The category $F_q\mbox{\rm -Isoc}^\dagger(X)\otimes L$ is an $L$-linear Abelian category. 
We will show that it is furnished with tensor products, duals and the unit object $j^\dagger_X\mathcal O_{]\overline{X}[}\otimes_KL$ 
in the propositions below. 
\end{enumerate}
\end{example}

\vspace*{3mm}

Let us study tensor products, duals, and extensions of $F$-isocrystals with $\overline{\mathbb Q}_p$-structures, and introduce 
an invariant $\overline{\mathbb Q}_p\mbox{\rm -rank}$ below. 
For an overconvergent $F$-isocrystal $\mathcal M^\dagger$, 
$\mathrm{rank}(\mathcal M^\dagger)$ means the rank of locally free $j^\dagger_X\mathcal O_{]X[}$-module $\mathcal M^\dagger$ 
on each connected component of $X$.

\begin{proposition} Let $\mathcal M^\dagger_i$ be an overconvergent $L_i$-$F_q$-isocrystal on $X$, 
and $L = L_1L_2$ the field of composite between $L_1$ and $L_2$ in $\overline{\mathbb Q}_p$. 
Let us put $\mathcal M^\dagger = (\mathcal M^\dagger_1\otimes\mathcal M_2^\dagger)\otimes_{L_1\otimes_KL_2}L$ 
where $\mathcal M^\dagger_1\otimes_{j^\dagger\mathcal O_{]X[}}\mathcal M_2^\dagger$ is the tensor product as an overconvergent $F$-isocrystal on $X/K$ 
with $L_1\otimes_KL_2$-actions and $L_1\otimes_KL_2 \rightarrow L\, (a\otimes b \mapsto ab)$. 
\begin{enumerate}
\item If we define $\nabla_{\mathcal M^\dagger} 
= (\nabla_{\mathcal M^\dagger_1}\otimes\mathrm{id}_{\mathcal M^\dagger_2} 
+ \mathrm{id}_{\mathcal M^\dagger_1}\otimes\nabla_{\mathcal M^\dagger_2})\otimes\mathrm{id}_L$, 
then the pair $(\mathcal M^\dagger, \nabla_{\mathcal M^\dagger})$ is an overconvergent isocrystal on $X/K$ 
such that 
$$
    \frac{\mathrm{rank}(\mathcal M^\dagger)}{\mathrm{deg}(L/K)} 
    = \frac{\mathrm{rank}(\mathcal M^\dagger_1)}{\mathrm{deg}(L_1/K)}\times\frac{\mathrm{rank}(\mathcal M^\dagger_2)}{\mathrm{deg}(L_2/K)} 
$$
on each connected component of $X$. 
\item If we define $F_{\mathcal M^\dagger} : F^\ast\mathcal M^\dagger \rightarrow \mathcal M^\dagger$ by 
$F_{\mathcal M^\dagger} = F_{\mathcal M_1^\dagger\otimes\mathcal M_2^\dagger}\otimes\mathrm{id}_L$ where 
$F^\ast\mathcal M^\dagger = F^\ast(\mathcal M_1^\dagger\otimes\mathcal M_2^\dagger)\otimes_{L_1\otimes_KL_2}L$, 
then $F_{\mathcal M^\dagger}$ is horizontal with respect to the connections. 
In particular,  the triplet $(\mathcal M^\dagger, \nabla_{\mathcal M^\dagger}, F_{\mathcal M^\dagger})$ is an overconvergent $F$-isocrystal 
on $X/K$. 
\item For an element $a \in L$, we define $a : \mathcal M^\dagger \rightarrow \mathcal M^\dagger$ 
by $a = \mathrm{id}_{\mathcal M_1^\dagger\otimes\mathcal M^\dagger_2}\otimes a\mathrm{id}_L$, then 
it defines an $L$-structure on $\mathcal M^\dagger$.
\end{enumerate}
In particular, $\mathcal M^\dagger$ is an overconvergent $L$-$F_q$-isocrystal on $X$. 
\end{proposition}

\begin{proposition}
Let $\mathcal M^\dagger$ be an overconvergent $L$-$F_q$-isocrystal on $X$. We define an action of $L$ on $(\mathcal M^\dagger)^\vee 
= \mathrm{Hom}_{j^\dagger_X\mathcal O_{]\overline{X}[}}(\mathcal M^\dagger, j^\dagger_X\mathcal O_{]\overline{X}[})$ by
$$
    (a, \eta) \mapsto a\eta \hspace*{5mm} (a\eta)(m) = \eta(am)\, \, \, \mbox{\rm for}\, \,\,  m \in \mathcal M^\dagger, 
$$
for $a \in L$ and $\eta \in (\mathcal M^\dagger)^\vee$, then it is an $L$-structure on $(\mathcal M^\dagger)^\vee$. 
In particular, $(\mathcal M^\dagger)^\vee$ is an overconvergent $L$-$F_q$-isocrystal on $X$ such that 
$$
    \frac{\mathrm{rank}((\mathcal M^\dagger)^\vee)}{\mathrm{deg}(L/K)} 
    = \frac{\mathrm{rank}(\mathcal M^\dagger)}{\mathrm{deg}(L/K)} 
$$
on each connected component of $X$. Moreover there is a perfect pairing 
$$
    (\mathcal M^\dagger \otimes_KL) \otimes ((\mathcal M^\dagger)^\vee\otimes_KL) \rightarrow j^\dagger_X\mathcal O_{]\overline{X}[}\otimes_KL
$$
of overconvergent $L$-$F_q$-isocrystals on $X$. In other words, 
$(\mathcal M^\dagger)^\vee$ is the dual of $\mathcal M^\dagger$ in $F_q\mbox{\rm -Isoc}^\dagger(X)\otimes L$. 
\end{proposition}

\begin{proposition} Let $L_1 \subset L_2$ be extensions of $K$ in $\overline{\mathbb Q}_p$. 
We define a functor $F_q\mbox{\rm -Isoc}^\dagger(X)\otimes L_1 \rightarrow F_q\mbox{\rm -Isoc}^\dagger(X)\otimes L_2$ of scalar extensions by 
$$
      \theta_{L_1, L_2} : \mathcal M^\dagger \mapsto 
      (\mathcal M^\dagger \otimes_KL_2)\otimes_{L_1\otimes_KL_2}L_2 = \mathcal M^\dagger \otimes_{L_1}L_2. 
$$
Then the equality 
$$
   \frac{\mathrm{rank}(\theta_{L_1, L_2}(\mathcal M^\dagger))}{\mathrm{deg}(L_2/K)}
   = \frac{\mathrm{rank}(\mathcal M^\dagger)}{\mathrm{deg}(L_1/K)}. 
$$
holds on each connected component of $X$. 
Moreover, $\theta_{L_1, L_2} = \theta_{L_2, L_3}\circ\theta_{L_1, L_2}$ for $L_1 \subset L_2 \subset L_3$ in $\overline{\mathbb Q}_p$. 
\end{proposition}

\begin{definition}
We define the category $F_q\mbox{\rm -Isoc}^\dagger(X)\otimes \overline{\mathbb Q}_p$ 
by the $2$-colimit of $F_q\mbox{\rm -Isoc}^\dagger(X)\otimes L$ over all finite extensions $L$ in $\overline{\mathbb Q}_p$ by $\theta_{L_1, L_2}$ 
in the previous proposition. An object of $F_q\mbox{\rm -Isoc}^\dagger(X)\otimes \overline{\mathbb Q}_p$ is said to be an 
overconvergent $\overline{\mathbb Q}_p$-$F_q$-isocrystals on $X$. 
We define the $\overline{\mathbb Q}_p\mbox{\rm -rank}$ of an overconvergent $\overline{\mathbb Q}_p$-$F_q$-isocrystal $\mathcal M^\dagger$ 
on $X$ by
$$
      \overline{\mathbb Q}_p\mbox{\rm -rank}(\mathcal M^\dagger) = \frac{\mathrm{rank}(\mathcal M^\dagger)}{\mathrm{deg}(L/K)} 
$$
on each connected component of $X$. 
We also define the category $F_q\mbox{\rm -Isoc}(X)\otimes \overline{\mathbb Q}_p$ of convergent $\overline{\mathbb Q}_p$-$F_q$-isocrystals on $X$ similarly. 
\end{definition}

\begin{lemma}\label{irredQver} 
Let $\mathcal M^\dagger$ be an 
irreducible overconvergent $\overline{\mathbb Q}_p$-$F_q$-isocrystal on $X$ 
which is represented by an overconvergent $L$-$F_q$-isocrystal on $X$ for an extension $L$ of $K$ in 
$\overline{\mathbb Q}_p$. Then there exists an 
irreducible overconvergent $F$-isocrystal $\mathcal N^\dagger$ on $X/K$ such that 
$\mathcal M^\dagger$ is a direct sum of a finite number of copies of $\mathcal N^\dagger$ 
as an overconvergent $F$-isocrystal on $X/K$. 
\end{lemma}

\prf{Let $\mathcal N^\dagger$ be a nontrivial irreducible subobject of $\mathcal M^\dagger$ 
as an overconvergent $F$-isocrystal on $X/K$ with respect to the $q$-Frobenius $\sigma$. 
The irreducibility $\mathcal M^\dagger$ implies that 
the natural morphism $\mathcal N^\dagger\otimes_KL \rightarrow 
\mathcal M^\dagger$ as overconvergent $F$-isocrystals on $X/K$ is surjective 
where $L$ is regarded as an $F_q$-space over $K$ such that $F_L = \mathrm{id}_L$. 
Hence $\mathcal M^\dagger$ is a direct sum of a finite number of copies of $\mathcal N^\dagger$.
}

\vspace*{3mm}

Let $k_0$ be a subfield of $k$ with $\mathrm{deg}(k/k_0) = n$, and $K_0$ the unramified extension 
of $\mathbb Q_p$ with the residue field $k_0$ and the $q_0$-Frobenius $\sigma_0 = \mathrm{id}_{K_0}$ 
where $q_0^n = q$. 
Since $K$ is finite and unramified over $K_0$, the smooth formal scheme over $R$ 
can be regarded as a smooth formal scheme over the integer ring $R_0$ of $K_0$. 
Hence, for an overconvergent $L$-$F_q$-isocrystal $\mathcal M^\dagger$ on $X$, 
$$
  \begin{array}{l}
(\mathcal M^\dagger_0, \nabla_{\mathcal M_0^\dagger}) = (\oplus_{i=0}^{n-1}(F_{q_0}^i)^\ast\mathcal M^\dagger, 
\oplus_{i=0}^{n-1}(F_{q_0}^i)^\ast\nabla_{\mathcal M^\dagger}) \\
F_{\mathcal M^\dagger_0}(a_0\otimes_{\sigma_0}m_0, \cdots, a_{n-1}\otimes_{\sigma_0}m_{n-1}) 
= (F_{\mathcal M^\dagger}(a_{n-1}\otimes_{\sigma_0}m_{n-1}), a_0\otimes_{\sigma_0}m_0, \cdots, a_{n-2}\otimes_{\sigma_0}m_{n-2}) 
\end{array}
$$
has a structure of overconvergent $L$-$F_{q_0}$-isocrystals on $X$. 

\begin{lemma}\label{aex} Suppose that $L$ admits an extension $\sigma_0$ 
of $q_0$-Frobenius $\sigma_0$ on $K_0$. The $L$-structure on $(\mathcal M^\dagger_0, \nabla_{\mathcal M_0^\dagger}, F_{\mathcal M^\dagger_0})$ 
which is induced from that on $\mathcal M^\dagger$ is given by
$$
   ((a_0\otimes_{\sigma_0}m_0, \cdots, a_{l-1}\otimes_{\sigma_0}m_{l-1}), b) \mapsto (ba_0\otimes_{\sigma_0}m_0, \cdots, \sigma_0^{n-1}(b)a_{n-1}\otimes_{\sigma_0}m_{n-1}). 
$$
for $b \in L$. In particular, the $K_0$-structure induced by the $q_0$-Frobenius coincides with 
$K_0 \subset L \rightarrow \mathrm{End}_{\footnotesize \mbox{\rm $F$-Isoc}}(\mathcal M_0^\dagger)$. 
\end{lemma}

\begin{proposition}\label{abequiv} \cite[Corollary 1.4.11]{Ab18b} With the notation as above, 
the functor 
$$
   \Pi_{q/q_0} : F_q\mbox{\rm -Isoc}^\dagger(X)\otimes \overline{\mathbb Q}_p 
   \rightarrow F_{q_0}\mbox{\rm -Isoc}^\dagger(X)\otimes \overline{\mathbb Q}_p
   \hspace*{5mm} \mathcal M^\dagger \mapsto \mathcal M^\dagger_0
$$
is an equivalence of categories. 
Moreover, the equalities hold
$$
    \overline{\mathbb Q}_p\mbox{\rm -rank}(\mathcal M^\dagger) 
    = \frac{\mathrm{rank}_{j^\dagger_X\mathcal O_{]\overline{X}[}}(\mathcal M^\dagger)}{\mathrm{deg}(L/K)} 
    = \frac{\mathrm{rank}_{j^\dagger_X\mathcal O_{]\overline{X}[}}(\mathcal M^\dagger_0)}{\mathrm{deg}(L/K_0)} 
    = \overline{\mathbb Q}_p\mbox{\rm -rank}(\mathcal M^\dagger_0). 
$$ 
on each connected component of $X$.
\end{proposition}

\prf{For any finite extension $L'$ of $K$, there exists a finite extension $L$ of $L'$ which admits a $q_0$-Frobenius $\sigma_0$ 
which is an extension of $\sigma_0$ on $K_0$ by Lemma \ref{ffix} (3) in Appendix \ref{Frob}. 
For an object $\mathcal N^\dagger$ of $F_{q_0}\mbox{\rm -Isoc}^\dagger(X)\otimes L$, 
if $\pi_{(K, \sigma)/(K_0, \sigma_0)}^\ast$ is the pull back functor defined in Appendix \ref{coeffrob}, then 
$\pi^\ast_{(K, \sigma)/(K_0, \sigma_0)}\mathcal N^\dagger$ has a natural $K$-structure induced by 
the $q$-Frobenius so that it admits a $K\otimes_{K_0}L$-structure. 
Let us define an overconvergent $L$-$F_q$-isocrystal $\mathcal M^\dagger$ on $X$ by 
$$
     \mathcal M^\dagger = (\pi^\ast_{(K, \sigma)/(K_0, \sigma_0)}\mathcal N^\dagger)\otimes_{K\otimes_{K_0}L}L
$$
with Frobenius $F_{\mathcal M^\dagger} = (\pi^\ast_{(K, \sigma)/(K_0, \sigma_0)}F_{\mathcal N^\dagger})^n\otimes\mbox{\rm id}_L$, where 
$K\otimes_{K_0}L \rightarrow L\, (b \otimes a \mapsto ba)$. 
One can easily verify this correspondence $\mathcal N^\dagger \mapsto \mathcal M^\dagger$ 
is a quasi-inverse of the given functor. 
}

\vspace*{3mm}

Let $X, Y$ be schemes locally of finite type over $\mathrm{Spec}\, k$ and $f : X \rightarrow Y$ be a morphism of schemes. Then 
the usual pull back functor $f^\ast$ of overconvergent $F$-isocrystals induces the pull back functor
$$
      f^\ast : F_q\mbox{\rm -Isoc}^\dagger(Y)\otimes \overline{\mathbb Q}_p 
      \rightarrow F_q\mbox{\rm -Isoc}^\dagger(X)\otimes \overline{\mathbb Q}_p. 
$$
Then, for an overconvergent $\overline{\mathbb Q}_p$-$F_q$-isocrystal $\mathcal N^\dagger$ on $Y$, we have
$$
      \overline{\mathbb Q}_p\mbox{\rm -rank}(f^\ast\mathcal N^\dagger) = \overline{\mathbb Q}_p\mbox{\rm -rank}(\mathcal N^\dagger)
$$
on each connected component.

\subsection{Changes of base fields}

Let $k_n$ be a finite extension of $k$ such that $\mathrm{deg}(k_n/k) = n$, 
$K_n$ the unramified extension of $K$ in $\overline{K}$ with the residue field $k_n$ 
and denote the $q$-Frobenius of $K_n$ by the same symbol $\sigma$. 
For a scheme $X$ locally of finite type over $\mathrm{Spec}\, k$, 
let us put $X_n = X \times_{\mathrm{Spec}\, k}\mathrm{Spec}\, k_n$ with the projection 
$\pi_n : X_n \rightarrow X$. We define the push forward functor by
$$
     \pi_{n , \overline{\mathbb Q}_p\, \ast} : F_{q^n}\mbox{\rm -Isoc}^\dagger(X_n)\otimes \overline{\mathbb Q}_p \rightarrow 
     F_q\mbox{\rm -Isoc}^\dagger(X)\otimes \overline{\mathbb Q}_p \hspace*{5mm}
     \mathcal M^\dagger \mapsto \pi_{(K_n, \sigma^n)/(K, \sigma)\, \ast}\mathcal M^\dagger = \pi_{(K_n, \sigma)/(K, \sigma)\, \ast}\mathcal M^\dagger_0
$$
where $\pi_{(K_n, \sigma^n)/(K, \sigma)\, \ast}$ and $\pi_{(K_n, \sigma)/(K, \sigma)\, \ast}$ are 
the push forward functor defined in Appendix \ref{coeffrob} and 
$\mathcal M_0^\dagger$ is the $\overline{\mathbb Q}_p$-$F_\sigma$-isocrystal on $X_n$  
introduced just before Lemma \ref{aex}. If $\pi_{(K_n, \sigma^n)/(K, \sigma)}^\ast$ is the pull back functor 
defined in Appendix \ref{coeffrob}, then $\pi_{(K_n, \sigma^n)/(K, \sigma)}^\ast\mathcal M^\dagger$ is furnished 
with an extra $K_n$-structure arising from the $q^n$-Frobenius structures. 
We define the pull back functor by 
$$
     \pi_{n , \overline{\mathbb Q}_p}^\ast : F_q\mbox{\rm -Isoc}^\dagger(X)\otimes \overline{\mathbb Q}_p \rightarrow 
     F_{q^n}\mbox{\rm -Isoc}^\dagger(X_n)\otimes \overline{\mathbb Q}_p \hspace*{5mm}
     \mathcal N^\dagger \mapsto (\pi_{(K_n, \sigma^n)/(K, \sigma)}^\ast\mathcal N^\dagger)\otimes_{K_n\otimes_KL}L
$$
in which two $K_n$-structures from the $q^n$-Frobenius and from the $\overline{\mathbb Q}_p$-structure are identified 
by tensoring $\otimes_{K_n\otimes_KL}L$. One can easily see the 
definition of the pull back functor $\pi_{n, \overline{\mathbb Q}_p}^\ast$ does not depends on the choice of $L$. 
Then $\pi_{n , \overline{\mathbb Q}_p}^\ast$ is a left adjoint of $\pi_{n , \overline{\mathbb Q}_p\, \ast}$. 
Moreover the equalities 
$$
   \begin{array}{c}
        \displaystyle{\overline{\mathbb Q}_p\mbox{\rm -rank}(\pi_{n , \overline{\mathbb Q}_p\, \ast}\mathcal M^\dagger) 
      = \frac{\mathrm{rank}( \pi_{n, \overline{\mathbb Q}_p\, \ast}\mathcal M^\dagger)}{\mathrm{deg}(L/K)} 
      = \frac{\mathrm{rank}(\mathcal M^\dagger)\mathrm{deg}(K_n/K)^2}{\mathrm{deg}(L/K)} = \overline{\mathbb Q}_p\mbox{\rm -rank}(\mathcal M^\dagger)}\times n \\
     \displaystyle{\overline{\mathbb Q}_p\mbox{\rm -rank}(\pi_{n , \overline{\mathbb Q}_p}^\ast\mathcal N^\dagger) 
      = \frac{\mathrm{rank}( \pi_{n, \overline{\mathbb Q}_p}^\ast\mathcal N^\dagger)}{\mathrm{deg}(L/K_n)} 
      = \frac{\mathrm{rank}(\mathcal N^\dagger)}{\mathrm{deg}(L/K)} = \overline{\mathbb Q}_p\mbox{\rm -rank}(\mathcal N^\dagger)}
      \end{array}
$$
hold on each connected component. 

\begin{proposition} Suppose $X$ is separated of finite type over $\mathrm{Spec}\, k$. 
Let $\mathcal M^\dagger$ be an overconvergent $L$-$F_q$-isocrystal $\mathcal M^\dagger$ on $X$ 
for a finite extension $L$ over $K_n$ in $\overline{\mathbb Q}_p$. 
Then the rigid cohomology $H^i_{\mathrm{rig}}(X/K, \mathcal M^\dagger)$ is an $F_q$-space over $K$ with an $L$-structure 
(an $L$-$F_q$-space over $K$ for short)
$$
        L \rightarrow \mathrm{End}_{\mbox{\footnotesize \rm $F$-sp}}(H^i_{\mathrm{rig}}(X/K, \mathcal M^\dagger))
$$
induced from the $L$-structure on $\mathcal M^\dagger$. The base change homomorphism of rigid cohomology 
induces an isomorphism 
$$
      H^i_{\mathrm{rig}}(X_n/K_n, \pi_{n, \overline{\mathbb Q}_p}^\ast\mathcal M^\dagger) 
      \rightarrow (K_n \otimes_KH^i_{\mathrm{rig}}(X/K, \mathcal M^\dagger))\otimes_{K_n \otimes_KL}L 
      \cong H^i_{\mathrm{rig}}(X/K, \mathcal M^\dagger)
$$
of $L$-$F_{q^n}$-spaces over $K_n$. Here the right hand side has a $K_n$-space structure by $K_n \subset L$, 
and the Frobenius on the right hand side is $F_{H^i_{\mathrm{rig}}(X/K, \mathcal M^\dagger)}^n$.  
The results do not depend on the choice of $L$ in the category of $\overline{\mathbb Q}_p$-$F_q$-spaces over $K$. 
The same hold for the rigid cohomology $H^i_{\mathrm{rig}, c}(X/K, \mathcal M^\dagger)$ 
with compact supports. 
\end{proposition}

\subsection{Lefschetz trace formula} Let $X$ be a connected scheme separated of finite type over $\mathrm{Spec}\, k$ which is pure of dimension $d$, 
and $\mathcal M^\dagger$ an object in $F\mbox{\rm -Isoc}^\dagger(X)\otimes \overline{\mathbb Q}_p$ 
such that $\mathcal M^\dagger$ is represented as an overconvergent $L$-$F_q$-isocrystal on $X$ 
for an extension $L$ of $K$ in $\overline{\mathbb Q}_p$. 
Then, for any closed point $\alpha$ in $X$ with $\mathrm{deg}(k_\alpha/k) = n$ such that $K_\alpha \subset L$, 
the inverse image $i_{\alpha, X}^\ast\mathcal M^\dagger$ is an $L$-$F_q$-space over $K_\alpha$, that is, 
an $F$-isocrystal on $\alpha/K$ with respect to $q$-Frobenius $\sigma$ and with an $L$-structure. 
For the linearization of Frobenius, we define a $\overline{\mathbb Q}_p$-$F_{q^n}$-space 
$(\mathcal M^\dagger_{\alpha, \overline{\mathbb Q}_p}, F_{\mathcal M^\dagger_{\alpha, \overline{\mathbb Q}_p}})$ 
over $K_\alpha$ by the colimit of 
$$
(\mathcal M^\dagger_{\alpha, L}, F_{\mathcal M^\dagger_{\alpha, L}}) 
= (\pi_{(K_\alpha, \sigma^n)/(K, \sigma)}^\ast(i_{\alpha, X}^\ast\mathcal M^\dagger)\otimes_{K_n\otimes_KL}L, 
(\pi_{(K_\alpha, \sigma^n)/(K, \sigma)}^\ast(i_{\alpha, X}^\ast(F_{\mathcal M^\dagger}))\otimes\mathrm{id}_L)^n), 
$$
over all finite extensions $L$ of $K$ in $\overline{\mathbb Q}_p$, that is, 
the object of the quasi-inverse functor of $(i_{\alpha, X}^\ast\mathcal M^\dagger, i_{\alpha, X}^\ast(F_{\mathcal M^\dagger}))$ in Proposition \ref{abequiv}. 
By definition we have 
$$
\overline{\mathbb Q}_p\mbox{\rm -rank}(\mathcal M^\dagger_{\alpha, \overline{\mathbb Q}_p}) 
= \overline{\mathbb Q}_p\mbox{\rm -rank}(i_{\alpha, X}^\ast\mathcal M^\dagger) = \overline{\mathbb Q}_p\mbox{\rm -rank}(\mathcal M^\dagger). 
$$
The pair $(\mathcal M^\dagger_{\alpha, L}, F_{\mathcal M^\dagger_{\alpha, L}})$ is 
regarded both as an $L$-space of dimension $\overline{\mathbb Q}_p\mbox{\rm -rank}(\mathcal M^\dagger_{\alpha, \overline{\mathbb Q}_p}) $ 
with an $L$-linear endomorphism $F_{\mathcal M^\dagger_{\alpha, L}}$ 
and as a $K$-space of dimension $\overline{\mathbb Q}_p\mbox{\rm -rank}(\mathcal M^\dagger_{\alpha, \overline{\mathbb Q}_p})\mathrm{deg}(L/K)$ 
with a $K$-linear endomorphism $F_{\mathcal M^\dagger_{\alpha, L}}$. 
On the other hand the pair $(\mathcal M^\dagger_{\alpha, \overline{\mathbb Q}_p}, F_{\mathcal M^\dagger_{\alpha, \overline{\mathbb Q}_p}})$ 
is just a $\overline{\mathbb Q}_p$-space with a $\overline{\mathbb Q}_p$-linear endomorphism $F_{\mathcal M^\dagger_{\alpha, \overline{\mathbb Q}_p}}$. 

\begin{lemma}\label{pt} With the notation as above, the following hold. 
\begin{enumerate}
\item The natural homomorphism 
$$
      (i_{\alpha, X}^\ast\mathcal M^\dagger,i_{\alpha, X}^\ast(F_{\mathcal M^\dagger})^n) \rightarrow 
      (\mathcal M^\dagger_{\alpha, \overline{\mathbb Q}_p}, F_{\mathcal M^\dagger_{\alpha, \overline{\mathbb Q}_p}})
      \hspace*{5mm} m \mapsto m \otimes 1
$$
induces an isomorphism of $L$-spaces with an $L$-linear endomorphism. 
\item For any $k_n$-rational point $\beta$ mapping to the closed point $\alpha$ with 
the composite $i_{\beta, X} = i_{\alpha, X}\circ i_{\beta, \alpha}: \beta \overset{\cong}{\rightarrow} \alpha \rightarrow X$ of morphisms, we have 
$$
     \begin{array}{ccc}
       \mathrm{Tr}_L\left(F_{i_{\beta, X}^\ast\mathcal M^\dagger}^n; i_{\beta, X}^\ast\mathcal M^\dagger\right) 
       &= &\mathrm{Tr}_L\left(F_{\mathcal M^\dagger_{\alpha, L}}; \mathcal M^\dagger_{\alpha, L}\right) \\
       \mathrm{Tr}_K\left(F_{i_{\beta, X}^\ast\mathcal M^\dagger}^n; i_{\beta, X}^\ast\mathcal M^\dagger\right) 
       &= &\mathrm{Tr}_K\left(F_{\mathcal M^\dagger_{\alpha, L}}; \mathcal M^\dagger_{\alpha, L}\right)
     \end{array}
$$
where $\mathrm{Tr}_E$ means the trace of $E$-linear homomorphisms. 
\item If $I_K(L, \overline{\mathbb Q}_p)$ is a set of $K$-algebra 
homomorphisms from $L$ to $\overline{\mathbb Q}_p$, then the equality 
$$
   \mathrm{Tr}_K\left(F_{\mathcal M^\dagger_{\alpha, \overline{\mathbb Q}_p}}; \mathcal M^\dagger_{\alpha, \overline{\mathbb Q}_p}\right) = 
     \sum_{\tau \in I_K(L, \overline{\mathbb Q}_p)} 
     \tau\left(\mathrm{Tr}_L\left(F_{\mathcal M^\dagger_{\alpha, L}}; \mathcal M^\dagger_{\alpha, L}\right)\right)
$$
holds
\end{enumerate}
\end{lemma}

\prf{(2) Let $\tau : K_\alpha \rightarrow K_\beta$ be the $K$-algebra isomorphism which is the lift of the $k$-algebra isomorphism 
$K_\alpha \rightarrow K_\beta$. Then $i_{\beta, X}^\ast\mathcal M^\dagger = i_{\beta, \alpha}^\ast(i_{\alpha, X}^\ast\mathcal M^\dagger)$ and 
$$
      F_{i_{\beta, X}^\ast\mathcal M^\dagger} = 
      (\mathrm{id}_{i_{\alpha, X}^\ast\mathcal M^\dagger}\otimes\tau)\circ F_{i_{\alpha, X}^\ast\mathcal M^\dagger} \circ 
      (\mathrm{id}_{i_{\alpha, X}^\ast\mathcal M^\dagger}\otimes\tau^{-1}).
$$
Hence it implies the assertion.

(3) follows from (2) and the isomorphism $L\otimes_KL \cong \prod_{\tau \in I_K(L, \overline{\mathbb Q}_p)}\tau(L)$.}

\vspace*{3mm}

The Lefschetz trace formula for overconvergent $\overline{\mathbb Q}_p$-$F$-isocrystals follows from the Lefschetz trace formula 
of usual overconvergent $F$-isocrystals \cite[Th\'eor\`eme 6.2]{EL93}. 

\begin{proposition}\label{isoA0} \cite[Theorem A.3.2]{Ab18b}
With the notation as above, the trace formula
$$
      \sum_{\beta \in X(k_n)} 
      \mathrm{Tr}_L\left(F_{i_{\beta, X}^\ast\mathcal M^\dagger}^n; i_{\beta, X}^\ast\mathcal M^\dagger\right) 
      = \sum_{i=0}^{2d}\, (-1)^i\mathrm{Tr}_L(F_{H^i_{\mathrm{rig}, c}}^n; H_{\mathrm{rig}, c}^i(X/K, \mathcal M^\dagger))
$$
holds where $X(k_n)$ is the set of $k_n$-rational points of $X$. 
\end{proposition}

\prf{The idea of this proof is same with that of \cite[Theorem A.3.2]{Ab18b}. 
When $\mathrm{dim}(X) = 0$, the assertion is trivial. In general dimensional cases, 
it is sufficient to prove that the right hand side is $0$ when $X(k_n) = \emptyset$ by applying the excision sequence 
$$
    \cdots \rightarrow H_{\mathrm{rig}, c}^i((X \setminus Z)/K, \mathcal M^\dagger) \rightarrow H_{\mathrm{rig}, c}^i(X/K, \mathcal M^\dagger) 
    \rightarrow H_{\mathrm{rig}, c}^i(Z/K, \mathcal M^\dagger) \rightarrow H_{\mathrm{rig}, c}^{i+1}((X \setminus Z)/K, \mathcal M^\dagger) \rightarrow \cdots
$$
of rigid cohomology with compact supports, where $Z$ is a closed subscheme of $X$ containing the finite set $X(k_n)$. The exact sequence above is 
a sequence of $L$-$F_q$-spaces over $K$. 

Suppose $X(k_n) = \emptyset$. The right hand side of the equality
$$
      \sum_{\tau \in I_K(L, \overline{\mathbb Q}_p)}\tau\left(
      \sum_{i=0}^{2d}\, (-1)^i\mathrm{Tr}_L\left(F_{H^i_{\mathrm{rig}, c}}^n; H_{\mathrm{rig}, c}^i(X/K, \mathcal M^\dagger)\right)\right) 
      = \sum_{i=0}^{2d}\, (-1)^i\mathrm{Tr}_K\left(F_{H^i_{\mathrm{rig}, c}}^n; H_{\mathrm{rig}, c}^i(X/K, \mathcal M^\dagger)\right). 
$$
vanishes by the Lefschetz trace formula for usual overconvergent $F$-isocrystals \cite[Th\'eor\`eme 6.2]{EL93}, and hence 
$$
      \sum_{\tau \in I_K(L, \overline{\mathbb Q}_p)}\tau\left(
      \sum_{i=0}^{2d}\, (-1)^i\mathrm{Tr}_L\left(F_{H^i_{\mathrm{rig}, c}}^n; H_{\mathrm{rig}, c}^i(X/K, \mathcal M^\dagger)\right)\right) = 0. 
$$
For a nonzero element $\lambda$ in $L$, we define an $F_q$-space $L(\lambda)$ over $K$ with an $L$-structure by 
$$
     \left\{\begin{array}{l}
        L(\lambda) = L : \mbox{\rm as an $K$-space}\\
       F_{L(\lambda)}(m) = \lambda m\, \, \mbox{\rm for}\, \, m \in L \\
       L \times L(\lambda) \rightarrow L(\lambda)\hspace*{5mm} (a, m) \mapsto am. 
    \end{array}\right.
$$
Then $\mathcal M^\dagger\otimes_LL(\lambda)$ is an overconvergent $L$-$F_q$-isocrystal on $X$, and 
in particular the action of Frobenius is given by 
$$
        F_{\mathcal M^\dagger\otimes_LL(\lambda)} = \lambda F_{\mathcal M^\dagger}. 
$$
From the assumption $X(k_n) = \emptyset$ we also have 
$$
     \sum_{\tau \in I_K(L, \overline{\mathbb Q}_p)}
     \tau\left( \sum_{i=0}^{2d}\, (-1)^i\lambda\mathrm{Tr}_L\left(F_{H^i_{\mathrm{rig}, c}}^n; H_{\mathrm{rig}, c}^i(X/K, \mathcal M^\dagger)\right)\right) 
     = 0 
$$
for any $\lambda \in L^\times$. 
Since $\sum_{\tau \in I_K(L, \overline{\mathbb Q}_p)}\tau$ is a trace map of the extension $L$ over $K$ of characteristic $0$ 
and $\lambda$ is arbitrary, we have the desired vanishing
$$
     \sum_{i=0}^{2d}\, (-1)^i\mathrm{Tr}_L\left(F_{H^i_{\mathrm{rig}, c}}^n; H_{\mathrm{rig}, c}^i(X/K, \mathcal M^\dagger)\right) = 0. 
$$
}

\vspace*{3mm}

Let us define the $L$-function $L(X, \mathcal M^\dagger, t)$ of $\mathcal M^\dagger$ by 
$$
   \begin{array}{lll}
      L(X, \mathcal M^\dagger, t) &= &\displaystyle{\prod_{n=1}^\infty\prod_{\footnotesize\begin{array}{c} \alpha : \mbox{\rm closed points of $X$} 
      \\ \mathrm{deg}(k_\alpha/k) = n \end{array}}\, 
      \mathrm{det}_{\overline{\mathbb Q}_p}\left(1 - t^nF_{\mathcal M^\dagger_{\alpha, \overline{\mathbb Q}_p}} ; \mathcal M^\dagger_{\alpha, \overline{\mathbb Q}_p}\right)^{-1}} \\
      &=&\displaystyle{\mathrm{exp}\left(\sum_{n=1}^\infty \sum_{\beta \in X(k_n)}
      \frac{\mathrm{Tr}_{\overline{\mathbb Q}_p}\left(F_{i_{\beta, X}^\ast\mathcal M^\dagger}^n; i_{\beta, X}^\ast\mathcal M^\dagger\right)}{n}t^n\right) 
      \in \overline{\mathbb Q}_p[\hspace*{-0.2mm}[t]\hspace*{-0.2mm}]}. 
      \end{array}
$$
where $\mathrm{det}_E(1-tf) \in E[t]$ means a characteristic polynomial of an $E$-linear endomorphism $f$. 
Note that the second equality holds by Lemma \ref{pt}.

\begin{corollary}\label{isoA1} \cite[Corollary A.3.3]{Ab18b}
With the notation as above, we have 
$$
    L(X, \mathcal M^\dagger, t) = \prod_{i=0}^{2d} 
    \mathrm{det}_{\overline{\mathbb Q}_p}\left(1 - tF_{H^i_{\mathrm{rig}, c}}; H_{\mathrm{rig}, c}^i(X/K, \mathcal M^\dagger)\right)^{(-1)^{i+1}}. 
$$
\end{corollary}

\subsection{\v{C}ebatarev density theorem} 
Now we recall weights of overconvergent $\overline{\mathbb Q}_p$-$F$-isocrystals (see \cite[Section 10]{Ke16} for a brief introduction of weights). 
Let $\iota :\overline{\mathbb Q}_p \rightarrow \mathbb C$ be an isomorphism of fields. 

\begin{definition}\label{weight} Let $X$ be a connected scheme separated of finite type over $\mathrm{Spec}\, k$. 
\begin{enumerate}
\item An overconvergent $\overline{\mathbb Q}_p$-$F_q$-isocrystal $\mathcal M^\dagger$ on $X$ is $\iota$-pure of weight $w \in \mathbb Z$ if, 
for any closed point $\alpha$ in $X$ with $\mathrm{deg}(k_\alpha/k) = n$, any reciprocal root of the polynomial 
$$
   \mathrm{det}_{\overline{\mathbb Q}_p}\left(1 - tF_{\mathcal M^\dagger_{\alpha, \overline{\mathbb Q}_p}} ; 
   \mathcal M^\dagger_{\alpha, \overline{\mathbb Q}_p}\right) \in \overline{\mathbb Q}_p[t]
$$
has a complex absolute value $q^{nw/2}$ under the isomorphism $\iota$. 
An overconvergent $\overline{\mathbb Q}_p$-$F_q$-isocrystal on $X$ is $\iota$-pure if it is $\iota$-pure of weight $w$ for some $w \in \mathbb Z$. 
\item An overconvergent $\overline{\mathbb Q}_p$-$F_q$-isocrystal $\mathcal M^\dagger$ on $X$ is $\iota$-mixed if it is a successive extension of 
overconvergent $\overline{\mathbb Q}_p$-$F_q$-isocrystals of $\iota$-pure. 
\end{enumerate}
\end{definition}

\vspace*{3mm}

The following proposition is the \v{C}ebatarev density theorem for overconvergent $\overline{\mathbb Q}_p$-$F$-isocrystals \cite[Proposition A.4.1]{Ab18b}.

\begin{proposition}\label{ceb} Let $X$ be a smooth connected scheme separated of finite type over $\mathrm{Spec}\, k$, 
and $\mathcal M^\dagger$ and $\mathcal N^\dagger$ irreducible 
overconvergent $\overline{\mathbb Q}_p$-$F_q$-isocrystals on $X$ such that $(\mathcal N^\dagger)^\vee \otimes \mathcal M^\dagger$ 
is $\iota$-pure. 
Suppose that there is an open dense subscheme $U$ of $X$ such that, 
for any closed point $\alpha$ of $U$ with $\mathrm{deg}(k(\alpha)/k)=n$, 
there is an isomorphism 
$$
     (\mathcal M^\dagger_{\alpha, \overline{\mathbb Q}_p}, F_{\mathcal M^\dagger_{\alpha, \overline{\mathbb Q}_p}})
     \cong (\mathcal N^\dagger_{\alpha, \overline{\mathbb Q}_p}, F_{\mathcal N^\dagger_{\alpha, \overline{\mathbb Q}_p}})
$$
as $\overline{\mathbb Q}_p$-$F_{q^n}$-spaces over $K_\alpha$. Then the following hold. 
\begin{enumerate}
\item There exists an isomorphism 
$g^\dagger : \mathcal N^\dagger \rightarrow \mathcal M^\dagger$ 
of overconvergent $\overline{\mathbb Q}_p$-$F_q$-isocrystals on $X$. 
\item Suppose furthermore that both $\mathcal M^\dagger$ and $\mathcal N^\dagger$ admit slope filtrations as 
convergent $\overline{\mathbb Q}_p$-$F_q$-isocrystals and 
that there exists an isomorphism $h : \mathcal N/\mathcal N^1 \rightarrow \mathcal M/\mathcal M^1$ 
between the maximal slope quotients as convergent $\overline{\mathbb Q}_p$-$F_q$-isocrystals, and that either 
the canonical homomorphism 
$$
  \begin{array}{cl}
   \mbox{\rm End}_{\footnotesize \overline{\mathbb Q}_p\mbox{\rm -}F\mbox{\rm -Isoc}}(\pi_{m, \overline{\mathbb Q}_p}^\ast\mathcal M^\dagger) \rightarrow 
   \mbox{\rm End}_{\footnotesize \overline{\mathbb Q}_p\mbox{\rm -}F\mbox{\rm -Isoc}}(\pi_{m, \overline{\mathbb Q}_p}^\ast(\mathcal M/\mathcal M^1)) &\mbox{or}\, \\
      \mbox{\rm End}_{\footnotesize \overline{\mathbb Q}_p\mbox{\rm -}F\mbox{\rm -Isoc}}(\pi_{m, \overline{\mathbb Q}_p}^\ast\mathcal N^\dagger) 
      \rightarrow \mbox{\rm End}_{\footnotesize \overline{\mathbb Q}_p\mbox{\rm -}F\mbox{\rm -Isoc}}(\pi_{m, \overline{\mathbb Q}_p}^\ast(\mathcal N/\mathcal N^1)) 
   \end{array}
$$
is surjective for some positive integer $m$. Then there exists a unique isomorphism 
$g^\dagger : \mathcal N^\dagger \rightarrow \mathcal M^\dagger$ such that 
the induced morphism of the maximal slope quotients by $g^\dagger$ coincides with the given $h$. 
Here 
$\mbox{\rm End}_{\footnotesize \overline{\mathbb Q}_p\mbox{\rm -}F\mbox{\rm -isoc}}(\pi_{m, \overline{\mathbb Q}_p}^\ast\mathcal M^\dagger)$ 
\mbox{\rm (resp. $\mbox{\rm End}_{\footnotesize \overline{\mathbb Q}_p\mbox{\rm -}F\mbox{\rm -Isoc}}(\pi_{m, \overline{\mathbb Q}_p}^\ast(\mathcal M^\dagger/\mathcal M^\dagger_1))$)} is 
the $\overline{\mathbb Q}_p$-space of endomorphisms of $\pi_{m, \overline{\mathbb Q}_p}^\ast\mathcal M^\dagger$ 
\mbox{\rm (resp. $\pi_{m, \overline{\mathbb Q}_p}^\ast(\mathcal M^\dagger/\mathcal M^\dagger_1)$)} as an overconvergent 
\mbox{\rm (resp.} a convergent\mbox{)} $\overline{\mathbb Q}_p$-$F_q$-isocrystal on $X_m$, 
and the same for $\mathcal N^\dagger$. 
\end{enumerate}
\end{proposition}

\prf{(1) This is a proof in \cite[Proposition A.4.1]{Ab18b}. Denote the dimension of $X$ by $d$. Note that 
$(\mathcal N^\dagger)^\vee \otimes \mathcal M^\dagger$ 
is $\iota$-pure of weight $0$ by the hypothesis. Since there is an isomorphism 
$$
     \mathrm{Hom}_{\footnotesize \overline{\mathbb Q}_p\mbox{\rm -}F\mbox{\rm -Isoc}}(\mathcal N^\dagger, \mathcal M^\dagger) 
     \cong\, H^0_{\mathrm{rig}}(X/K, (\mathcal N^\dagger)^\vee \otimes \mathcal M^\dagger)^{F_{H^0_{\mathrm{rig}}}=1} 
     =\, H^0_{\mathrm{rig}}(U/K, (\mathcal N^\dagger)^\vee \otimes \mathcal M^\dagger)^{F_{H^0_{\mathrm{rig}}}=1} 
$$
and both $\mathcal M^\dagger$ and $\mathcal N^\dagger$ is irreducible on $X$, 
we have only to prove the nonvanishing of the right hand side. 
Here $\mathrm{Hom}_{\footnotesize \overline{\mathbb Q}_p\mbox{\rm -}F_q\mbox{\rm -Isoc}^\dagger}$ is a group of morphisms 
as overconvergent $\overline{\mathbb Q}_p$-$F_q$-isocrystals on $X$, 
the $F_{H^0_{\mathrm{rig}}}=1$ means the Frobenius invariant subspace of the $0$-th rigid cohomology,  
and the second equality follows from the full faithfulness of the restriction functor of overconvergent $F$-isocrystals \cite[Theorem 6.3.1]{Ts02}. 
By Poincar\'e duality \cite[Theorem 1.2.3]{Ke06} it is sufficient to prove an existence of a nontrivial cocycle in 
$H^{2d}_{\mathrm{rig}, c}(U/K, (\mathcal M^\dagger)^\vee \otimes \mathcal N^\dagger)$ on which Frobenius $F_{H^{2d}_{\mathrm{rig}, c}}$ acts by $q^d$. 
Because $(\mathcal M^\dagger)^\vee \otimes \mathcal N^\dagger$ is $\iota$-pure of weight $0$, 
the weights of Frobenius $F_{H^j_{\mathrm{rig}, c}}$ on $H^j_{\mathrm{rig}, c}(U/K, (\mathcal M^\dagger)^\vee \otimes \mathcal N^\dagger)$ 
are $\leq j$ \cite[Theorem 5.3.2]{Ke06b}. Hence what we want is that $L(U, (\mathcal M^\dagger)^\vee \otimes \mathcal N^\dagger, u)$ has a factor 
$1-q^du$ in the denominator. This holds since 
$$
     L(U, (\mathcal M^\dagger)^\vee \otimes \mathcal N^\dagger, u) = L(U, (\mathcal M^\dagger)^\vee \otimes \mathcal M^\dagger, u)
$$
by the hypothesis of coincidence of $\overline{\mathbb Q}_p$-$F_q$-spaces at each closed point $\alpha \in U$. 

(2) Since there exists an isomorphism $(g')^\dagger : \mathcal N^\dagger \rightarrow \mathcal M^\dagger$ in (1), 
the condition of endomorphisms for $\pi_{m, \overline{\mathbb Q}_p}^\ast\mathcal M^\dagger$ is equivalent to that for $\pi_{m, \overline{\mathbb Q}_p}^\ast\mathcal N^\dagger$. 
Suppose that 
$$
\mbox{\rm End}_{\footnotesize \overline{\mathbb Q}_p\mbox{\rm -}F\mbox{\rm -Isoc}}(\pi_{m, \overline{\mathbb Q}_p}^\ast\mathcal N^\dagger) \rightarrow 
   \mbox{\rm End}_{\footnotesize \overline{\mathbb Q}_p\mbox{\rm -}F\mbox{\rm -Isoc}}(\pi_{m, \overline{\mathbb Q}_p}^\ast(\mathcal N/\mathcal N^1))
$$
is surjective. By the hypothesis we can find a morphism 
$\tau_m : \pi_{m, \overline{\mathbb Q}_p}^\ast\mathcal N^\dagger \rightarrow \mathcal \pi_{m, \overline{\mathbb Q}_p}^\ast\mathcal M^\dagger$ 
which induces $\pi_{m, \overline{\mathbb Q}_p}^\ast((g')^{-1}\circ h)$ on $\pi_{m, \overline{\mathbb Q}_p}^\ast(\mathcal N/\mathcal N^1)$. 
By our construction $\pi_{m, \overline{\mathbb Q}_p}^\ast((g')^{-1})\circ \tau_m : 
\pi_{m, \overline{\mathbb Q}_p}^\ast\mathcal N^\dagger \rightarrow \pi_{m, \overline{\mathbb Q}_p}^\ast\mathcal M^\dagger$ 
is a morphism of overconvergent $\overline{\mathbb Q}_p$-$F_{q^m}$-isocrystals on $X_m$ such that 
the induced morphism of the maximal slope quotients by $\pi_{m, \overline{\mathbb Q}_p}^\ast((g')^{-1})\circ \tau_m$ coincides 
with the given $\pi_{m, \overline{\mathbb Q}_p}^\ast(h)$.                          
Let us consider the composite 
$$
   \begin{array}{ccccc}
      g^\dagger : \mathcal N^\dagger &\overset{\mathrm{ad}_{(K_m, \sigma^m)/(K, \sigma)}}{\rightarrow} 
      &\pi_{m, \overline{\mathbb Q}_p\, \ast}\pi_{(K_m, \sigma^m)/(K, \sigma)}^\ast\mathcal N^\dagger 
            &\rightarrow &\pi_{m, \overline{\mathbb Q}_p\, \ast}\pi_{m, \overline{\mathbb Q}_p}^\ast\mathcal N^\dagger \\
      &\overset{\pi_{m, \overline{\mathbb Q}_p\, \ast}(\pi_{m, \overline{\mathbb Q}_p}^\ast((g')^{-1})\circ \tau_m)}{\rightarrow}
       &\pi_{m, \overline{\mathbb Q}_p\, \ast}\pi_{m, \overline{\mathbb Q}_p}^\ast\mathcal M^\dagger 
       &\overset{\frac{1}{m}\mathrm{Tr}_{(K_m, \sigma^m)/(K, \sigma)}}{\rightarrow} &\mathcal M^\dagger
       \end{array}
$$
of the morphisms of overconvergent $L$-$F_q$-isocrystals on $X$ for a sufficiently large finite extension $L$ 
of $K$ in $\overline{\mathbb Q}_p$, where 
$\mathrm{ad}_{(K_m, \sigma^m)/(K, \sigma)}$ and $\mathrm{Tr}_{(K_m, \sigma^m)/(K, \sigma)}$ are defined in Appendix \ref{coeffrob}. 
Since the induced morphism of the maximal slope quotients by $g^\dagger$ coincides with the given $h$ by Lemma \ref{adtr}, 
$g^\dagger$ is the desired isomorphism by the irreducibility of $\mathcal M^\dagger$ and $\mathcal N^\dagger$. 
}

\vspace*{2mm}

The corollary below follows from Corollary \ref{isom}. 

\begin{corollary}\label{cebiso} Let $\mathcal M^\dagger$ and $\mathcal N^\dagger$ be irreducible 
overconvergent $\overline{\mathbb Q}_p$-$F_q$-isocrystals on $X$ such that  $(\mathcal N^\dagger)^\vee \otimes \mathcal M^\dagger$ 
is $\iota$-pure and 
admit slope filtrations, and $h : \mathcal N/\mathcal N^1 \rightarrow \mathcal M/\mathcal M^1$ 
a nontrivial morphism of convergent $\overline{\mathbb Q}_p$-$F_q$-isocrystals. 
Suppose that there is an open dense subscheme $U$ of $X$ such that, for any closed point $\alpha$ of $U$, 
both $(X, \alpha, \mathcal M^\dagger)$ and $(X, \alpha, \mathcal N^\dagger)$ satisfy the condition 
$(\mbox{\rm LC}_{\overline{\mathbb Q}_p})$ :  
\begin{flushleft}
\begin{tabular}{ll}
\hspace*{5mm} $(\mbox{\rm LC}_{\overline{\mathbb Q}_p})$ &\begin{tabular}{l} 
There exist a smooth curve $C_\alpha$ over $\mathrm{Spec}\, k$ and 
a morphism $i_{C_\alpha, X} : C_\alpha \rightarrow X$ \\
over $\mathrm{Spec}\, k$ such that $C_\alpha$ is passing at $\alpha$ and that the restriction $i_{C_\alpha, X}^\ast\mathcal M^\dagger$ on $C_\alpha$ \\
is an irreducible overconvergent $\overline{\mathbb Q}_p$-$F_q$-isocrystal on $C_\alpha$. 
\end{tabular}
\end{tabular} 
\end{flushleft}
Then the following hold. 
\begin{enumerate}
\item There exists an isomorphism 
$g^\dagger : \mathcal N^\dagger \rightarrow \mathcal M^\dagger$ 
of overconvergent $\overline{\mathbb Q}_p$-$F_q$-isocrystals on $X$. 
\item Suppose furthermore that either the canonical homomorphism 
$$
  \begin{array}{cl}
   \mbox{\rm End}_{\footnotesize \overline{\mathbb Q}_p\mbox{\rm -}F\mbox{\rm -Isoc}}(\pi_{m, \overline{\mathbb Q}_p}^\ast\mathcal M^\dagger) \rightarrow 
   \mbox{\rm End}_{\footnotesize \overline{\mathbb Q}_p\mbox{\rm -}F\mbox{\rm -Isoc}}(\pi_{m, \overline{\mathbb Q}_p}^\ast(\mathcal M/\mathcal M^1)) &\mbox{or}\, \\
      \mbox{\rm End}_{\footnotesize \overline{\mathbb Q}_p\mbox{\rm -}F\mbox{\rm -Isoc}}(\pi_{m, \overline{\mathbb Q}_p}^\ast\mathcal N^\dagger) 
      \rightarrow \mbox{\rm End}_{\footnotesize \overline{\mathbb Q}_p\mbox{\rm -}F\mbox{\rm -Isoc}}(\pi_{m, \overline{\mathbb Q}_p}^\ast(\mathcal N/\mathcal N^1)) 
   \end{array}
$$
is surjective for some positive integer $m$. 
Then there exists a unique isomorphism $g^\dagger : \mathcal N^\dagger \rightarrow \mathcal M^\dagger$ such that 
the induced morphism of the maximal slope quotients by $g^\dagger$ coincides with the given $h$. 
\end{enumerate}
\end{corollary}

\prf{Applying Corollary \ref{isom}, there exists an isomorphism 
$$
     (\mathcal M^\dagger_{\alpha, \overline{\mathbb Q}_p}, F_{\mathcal M^\dagger_{\alpha, \overline{\mathbb Q}_p}})
     \cong (\mathcal N^\dagger_{\alpha, \overline{\mathbb Q}_p}, F_{\mathcal N^\dagger_{\alpha, \overline{\mathbb Q}_p}})
$$
of $\overline{\mathbb Q}_p$-$F_{q^n}$-spaces over $K_\alpha$ for any closed point $\alpha \in U$ such that $\mathrm{deg}(k_\alpha/k) = n$. 
Hence the assertions follow from Proposition \ref{ceb}.}

\vspace*{3mm}

We give a sufficient condition of the bijectivity of endomorphisms in Corollary \ref{cebiso} (2). 

\begin{lemma}\label{suffcond} Let $X, Y$ be schemes separated of finite type over $\mathrm{Spec}\, k$ 
such that $X$ is connected and $Y$ is nonempty, and $f : Y \rightarrow X$ a morphism. 
For an overconvergent $\overline{\mathbb Q}_p$-$F_q$-isocrystal $\mathcal M^\dagger$ on $X$ which admits a slope filtration, 
both canonical homomorphisms 
$$
   \begin{array}{cl}
     \mbox{\rm End}_{\footnotesize {\overline{\mathbb Q}_p}\mbox{\rm -}F\mbox{\rm -Isoc}}(\mathcal M^\dagger) \rightarrow 
     \mbox{\rm End}_{\footnotesize {\overline{\mathbb Q}_p}\mbox{\rm -}F\mbox{\rm -Isoc}}(f^\ast\mathcal M^\dagger)\, \, \, &\mbox{\rm and} \\
     \mbox{\rm End}_{\footnotesize {\overline{\mathbb Q}_p}\mbox{\rm -}F\mbox{\rm -Isoc}}(\mathcal M/\mathcal M^1) 
     \rightarrow \mbox{\rm End}_{\footnotesize {\overline{\mathbb Q}_p}\mbox{\rm -}F\mbox{\rm -Isoc}}(f^\ast(\mathcal M/\mathcal M^1))
     \end{array}
$$
are injective. In particular, the canonical homomorphism 
$\mbox{\rm End}_{\footnotesize {\overline{\mathbb Q}_p}\mbox{\rm -}F\mbox{\rm -Isoc}}(\mathcal M^\dagger) 
\rightarrow \mbox{\rm End}_{\footnotesize {\overline{\mathbb Q}_p}\mbox{\rm -}F\mbox{\rm -Isoc}}(\mathcal M/\mathcal M^1)$ 
is bijective if the bijectivity 
$\mbox{\rm End}_{\footnotesize {\overline{\mathbb Q}_p}\mbox{\rm -}F\mbox{\rm -Isoc}}(\mathcal M^\dagger) \overset{\cong}{\rightarrow} 
\mbox{\rm End}_{\footnotesize {\overline{\mathbb Q}_p}\mbox{\rm -}F\mbox{\rm -Isoc}}(i^\ast_{C_\alpha, X}\mathcal M^\dagger)$ 
holds for a closed point $\alpha \in X$ and a smooth curve $i_{\alpha, C} : C_\alpha \rightarrow X$ passing at $\alpha$ 
such that $i_{C_\alpha, X}^\ast\mathcal M^\dagger$ is irreducible. 
\end{lemma}

\prf{The former assertion follows from the fact that the induced morphisms 
$$
\begin{array}{l} 
H_{\mathrm{rig}}^0(X/K, (\mathcal M^\dagger)^\vee\otimes\mathcal M^\dagger) \rightarrow 
H_{\mathrm{rig}}^0(Y/K, (f^\ast\mathcal M^\dagger)^\vee \otimes f^\ast\mathcal M^\dagger), \hspace*{2mm} \mbox{\rm and} \\
H_{\mathrm{conv}}^0(X/K, (\mathcal M/\mathcal M^1)^\vee \otimes \mathcal M/\mathcal M^1) 
\rightarrow H_{\mathrm{conv}}^0(Y/K, (f^\ast(\mathcal M/\mathcal M^1))^\vee \otimes f^\ast(\mathcal M/\mathcal M^1))
\end{array}
$$
of $L$-$F_q$-spaces over $K$ are injective by the hypothesis of connectedness. 
Here $H_{\mathrm{conv}}^0(X/K, (\mathcal M/\mathcal M^1)^\vee \otimes \mathcal M/\mathcal M^1)$ is the $0$-th convergent cohomology, that is, 
the space of horizontal sections of the convergent isocrystal $(\mathcal M/\mathcal M^1)^\vee \otimes \mathcal M/\mathcal M^1$. 
The second assertion follows from the former assertion and the bijectivity of 
$\mbox{\rm End}_{\footnotesize {\overline{\mathbb Q}_p}\mbox{\rm -}F\mbox{\rm -Isoc}}(i_{C_\alpha, X}^\ast\mathcal M^\dagger) 
\rightarrow \mbox{\rm End}_{\footnotesize {\overline{\mathbb Q}_p}\mbox{\rm -}F\mbox{\rm -Isoc}}(i_{C_\alpha, X}^\ast(\mathcal M/\mathcal M^1))$ 
by Corollary \ref{kedcur}. 
}

\subsection{Consequences from the companion theorem} After Abe's celebrated work on $p$-adic Langlands correspondence and the companion 
theorem in the case of curves \cite{Ab18b} \cite{Ab18}, 
Deligne's companion conjecture between $\ell$-adic coefficients and $p$-adic 
coefficients in general dimension \cite[Conjecture 1.2.10]{De81} 
is one of the hottest problems in arithmetic geometry. (See \cite{Cr92}, \cite{AE16} and \cite{Ke19} for details.)

Using Lefschetz theorem for isocrystals with tame ramifications along boundaries, Abe and Esnault obtained a weight thereom. 

\begin{theorem}\label{iotmix} \cite[Theorem 2.7]{AE16} Let $X$ be a smooth connected scheme separated of finite type over $\mathrm{Spec}\, k$. 
Then any $\overline{\mathbb Q}_p$-$F_q$-isocrystal on $X$ is $\iota$-mixed. 
In particular, any irreducible $\overline{\mathbb Q}_p$-$F_q$-isocrystal on $X$ is $\iota$-pure. 
\end{theorem}

Then Abe and Esnault established the Lefschetz theorem below after the companion theorem. 
Note that one can relax the condition on finite determinant of $\overline{\mathbb Q}_p$-$F_q$-isocrystals by twisting a character \cite[Theorem 6.1]{Ab18} 
and the Lefschetz theorem in \cite{AE16} 
is more general which asserts an existence of a curve passing at given finite closed points and on which 
the pull back is irreducible. 

\begin{theorem}\label{Lef} \cite[Theorem 3.10]{AE16} Let $X$ be a smooth connected scheme separated of finite type over $\mathrm{Spec}\, k$, 
and $\mathcal M^\dagger$ an irreducible $\overline{\mathbb Q}_p$-$F_q$-isocrystal on $X$. 
There exists a dense open subscheme $U$ of $X$ such that, for any closed point $\alpha \in U$, 
the triplet $(X, \alpha, \mathcal M^\dagger)$ satisfies the condition $(\mbox{\rm LC}_{\overline{\mathbb Q}_p})$ in Corollary \ref{cebiso}. 
\end{theorem}

As a consequence we have an affirmative answer to Kedlaya's question ``Minimal slope conjecture" \cite[Remark 5.14]{Ke16} (Conjecures \ref{Kcj}, \ref{Kcjdual}) 
for $\overline{\mathbb Q}_p$-$F$-isocrystals on smooth varieties over finite fields by Theorems \ref{Lef}, \ref{iotmix}, Corollary \ref{cebiso} and Lemma \ref{suffcond} : 

\begin{theorem}\label{rlt2} Let $X$ be a smooth connected scheme separated of finite type over $\mathrm{Spec}\, k$, 
and $\mathcal M^\dagger, \mathcal N^\dagger$ irreducible overconvergent $\overline{\mathbb Q}_p$-$F_q$-isocrystals on $X$ 
admitting slope filtrations of $\mathcal M$ and $\mathcal N$. 
Suppose that there exists a nontrivial morphism $h : \mathcal N/\mathcal N^1 \rightarrow \mathcal M/\mathcal M^1$ 
between the maximal slope quotients as convergent $\overline{\mathbb Q}_p$-$F_q$-isocrystals. 
Then the following hold. 
\begin{enumerate}
\item There exists an isomorphism 
$g^\dagger : \mathcal N^\dagger \rightarrow \mathcal M^\dagger$ of overconvergent $\overline{\mathbb Q}_p$-$F_q$-isocrystals on $X$. 
\item Suppose furthermore that, after a finite extension of $k$, there exists a smooth curve $i_{\alpha, C} : C_\alpha \rightarrow X$ passing at a closed point $\alpha \in X$ 
such that $i_{C_\alpha, X}^\ast\mathcal M^\dagger$ is irreducible and that the restriction map 
$$
\mbox{\rm End}_{\footnotesize {\overline{\mathbb Q}_p}\mbox{\rm -}F\mbox{\rm -Isoc}}(\mathcal M^\dagger) 
\rightarrow \mbox{\rm End}_{\footnotesize {\overline{\mathbb Q}_p}\mbox{\rm -}F\mbox{\rm -Isoc}}(i^\ast_{C_\alpha, X}\mathcal M^\dagger)
$$
is bijective. Then the isomorphism $g^\dagger$ in (1) induces a natural commutative diagram
$$
      \begin{array}{ccc}
          \mathcal N &\overset{g}{\rightarrow} &\mathcal M \\
          \downarrow & &\downarrow \\
          \mathcal N/\mathcal N^1 &\underset{h}{\rightarrow} &\mathcal M/\mathcal M^1
        \end{array}
$$
of convergent $\overline{\mathbb Q}_p$-$F$-isocrystals. 
\end{enumerate}
\end{theorem}

\appendix 

\section{Frobenius}\label{Frob}

\subsection{Frobenius on a complete discrete valuation field}\label{FrobK}
Let $R$ be a complete discrete valuation ring of mixed characteristic $(0, p)$ with the residue field $k = R/\mathbf m$ 
which is not necessarily perfect, and $K$ the field of fractions of $R$. 
For a positive power $q$ of $p$, a $q$-Frobenius $\sigma$ on $K$ is a continuous endomorphism 
$$
     \sigma : K \rightarrow K
$$
such that $\sigma(a)\, \equiv\, a^q\, (\mathrm{mod}\, \mathbf m)$ for $a \in R$. 
We define the $\sigma$-invariant subfield $K_\sigma$ by
$$
      K_\sigma = \{ a \in K\, |\, \sigma(a) = a \}.
$$
Remark that the same letter $\sigma$ is used for the $q$-Frobenius 
which acts on an extension $K'$ of $K$. 

\begin{lemma}\label{ffix} \mbox{\rm (\textit{c.f.} \cite[Remark 2.1]{Ts19})} Suppose $K$ admits a $q$-Frobenius $\sigma$. 
\begin{enumerate}
\item $K_\sigma$ is a complete discrete valuation field with a finite residue field of cardinal $\leq q$. 
In particular, $K_\sigma$ is a finite extension of the field $\mathbb Q_p$ of $p$-adic numbers. 
\item There exists a finite unramified extension $K'$ of $K$ such that $K'$ admits a $q$-Frobenius $\sigma$ 
which is an extension of $\sigma$ on $K$, the residue field $(K')_\sigma$ is a field of $q$-elements and 
the valuation group of $(K')_\sigma$ is same with that of $K$. 
If furthermore $k$ is perfect, then the natural map $(K')_\sigma \otimes_{W(\mathbb F_q)}W(k') \rightarrow K'$ is bijective 
where $k'$ is the residue field of $K'$ and 
$\sigma = \mathrm{id}_{(K')_\sigma}\otimes \mathrm{Frob}_p^s = \mathrm{id}_{(K')_\sigma}\otimes \mathrm{Frob}_q$ for $q = p^s$. 
Here $W(k')$ means the Witt vector ring with coefficients in $k'$ and $\mathrm{Frob}_p$ (resp. $\mathrm{Frob}_q$) is the canonical $p$-Frobenius 
(resp. $q$-Frobenius) on $W(k')$. 
\item For a finite extension $L$ of $\mathbb Q_p$, 
there exist a finite extension $K'$ of $K$ and a positive integer $n$ which satisfy the following: 
\begin{list}{}{}
\item[\mbox{\rm (i)}] $K'$ admits a $q$-Frobenius $\sigma$ which is an extension of $\sigma$ on $K$, 
\item[\mbox{\rm (ii)}] the residue field of $(K')_{\sigma^n}$ is a field of $q^n$-elements, 
\item[\mbox{\rm (iii)}] $L \subset (K')_{\sigma^n}$. 
\end{list}
\end{enumerate}
\end{lemma}

\vspace*{2mm}

In order to prove Lemma \ref{ffix} we introduce extensions $\widehat{K}^{\mathrm{perf}} \subset \widehat{K}^{\mathrm{perf, ur}}$ of $K$ 
as complete discrete valuation fields with the same valuation group to $K$ 
such that the residue field of $\widehat{K}^{\mathrm{perf}}$ (resp. $\widehat{K}^{\mathrm{perf, ur}}$) is a perfection of $k$ 
(resp. an algebraic closure of $k$). Moreover, $\widehat{K}^{\mathrm{perf}}$ and $\widehat{K}^{\mathrm{perf, ur}}$ admit a $q$-Frobenius $\sigma$ 
which is compatible with extensions of $K$. 
The field $\widehat{K}^{\mathrm{perf}}$ is defined by the $p$-adic completion of the inductive limit of the inductive system 
$$
       K \overset{\sigma}{\rightarrow} K \overset{\sigma}{\rightarrow} K \overset{\sigma}{\rightarrow} \cdots. 
$$
where $K = \underset{\rightarrow}{\mathrm{lim}}(K \overset{\sigma}{\rightarrow} \sigma(K) \overset{\sigma}{\rightarrow} \cdots)$ 
and the $q$-Frobenius on $\widehat{K}^{\mathrm{perf}}$ which is compatible to the Frobenius $\sigma$ on $K$ by the system 
$\sigma = (\sigma \overset{\sigma}{\rightarrow} \sigma \overset{\sigma}{\rightarrow} \sigma \overset{\sigma}{\rightarrow} \cdots)$ 
of Frobenius. $\widehat{K}^{\mathrm{perf, ur}}$ is the $p$-adic completion of the maximally unramified extension of $\widehat{K}^{\mathrm{perf}}$. 
Then the $q$-Frobenius extends uniquely on $\widehat{K}^{\mathrm{perf, ur}}$ and it is denoted by the same symbol $\sigma$. 

\vspace*{2mm}

\noindent
{\sc Proof of Lemma \ref{ffix}.} (1) One can easily see that $K_\sigma$ is a complete discrete valuation field. 
If $a \in K_\sigma \cap R$, then $a$ satisfies the congruence $a^q\, \equiv\, a\, \, (\mathrm{mod}\, \mathbf m)$. 
Hence the residue field of $K_\sigma$ is finite of cardinal $\leq q$. 

(2) Since $(\widehat{K}^{\mathrm{perf, ur}})_\sigma$ is a finite extension of $\mathbb Q_p$ with the residue 
field of $q$ elements and the valuation ring of $(\widehat{K}^{\mathrm{perf, ur}})_\sigma$ coincides that of $\widehat{K}^{\mathrm{perf, ur}}$, 
the composite field $K' = (\widehat{K}^{\mathrm{perf, ur}})_\sigma K$ is a desired field. 
When the residue field $k$ is perfect, then $W(k)$ is the canonically subring of $K$. 
Comparing the ramification, one has a natural isomorphism $(K')_\sigma \otimes_{W(\mathbb F_q)}W(k) \rightarrow K'$. 

(3) Replacing $K$ and $L$ by finite extensions respectively, 
we may assume that $K$ admits a $q$-Frobenius such that $(\widehat{K}^{\mathrm{perf, ur}})_{\sigma^n} = K_{\sigma^n}$ by the proof of (2), 
the cardinal of the residue field of $L$ is $q^n$, and $L$ is a totally ramified extension of $K_{\sigma^n}$. 
Since the valuation group of $(\widehat{K}^{\mathrm{perf, ur}})_{\sigma^n}$ 
coincides with that of $\widehat{K}^{\mathrm{perf, ur}}$ and hence of $K$, the natural map $L \otimes_{K_{\sigma^n}}K \rightarrow LK$ 
is an isomorphism of fields. 
Then $K' = LK$ and $\sigma = \mathrm{id}_L\otimes\sigma$ are our desired field and $q$-Frobenius. \hspace*{\fill} $\Box$

\subsection{Frobenius on $\dagger$-spaces associated to affine smooth schemes}\label{FrobA} 
Let $X$ be an affine smooth integral scheme separated of finite type 
over $\mathrm{Spec}\, k$. Suppose there exist an affine smooth integral scheme $\mathcal X = \mathrm{Spec}\, A_X$ of finite type over $\mathrm{Spec}\, R$ 
such that $\mathcal X \times_{\mathrm{Spec}\, R}\mathrm{Spec}\, k = X$. 
Note that an affine smooth scheme separated of finite type over $\mathrm{Spec}\, k$ always admits a smooth lift over $\mathrm{Spec}\, R$ 
\cite[Th\'eor\`eme 6]{El73}.  
Suppose $A_X = R[x_1, \cdots, x_N]/I$. If $\widehat{A}_X$ and $A_X^\dagger$ are the $p$-adically formal completion of $A_X$ 
and the $p$-adically weak completion of $A_X$ (i.e., a weakly complete finitely generated (w.c.f.g.) algebra over $(R, \mathbf m)$), then they are defined by 
$$
    \begin{array}{lll}
    \widehat{A}_X &= &R[x_1, \cdots, x_N]\widehat{\, \, \, }/IR[x_1, \cdots, x_N]\widehat{\, \, \, }, \\
      A_X^\dagger &= &R[x_1, \cdots, x_N]^\dagger/IR[x_1, \cdots, x_N]^\dagger
      \end{array}
$$
respectively, where $R[x_1, \cdots, x_N]\widehat{\, \, \, }$ is the $\mathbf m$-adic completion of $R[x_1, \cdots, x_N]$, and 
$$
  \begin{array}{l}
     R[x_1, \cdots, x_N]^\dagger = \underset{\lambda\rightarrow 1^+}{\mathrm{lim}}R[x_1, \cdots, x_N]_\lambda, \\
     R[x_1, \cdots, x_N]_\lambda = \left\{ f \in R[x_1, \cdots, x_N]\widehat{\, \, \, }\, \left|\, 
     \begin{array}{l} f\, \, \mbox{\rm is convergent on the closed} \\
      \mbox{\rm closed ball $\max_i |x_i| \leq \lambda$.} \end{array} \right. \right\}. 
      \end{array}
$$
Then $A_X^\dagger \subset \widehat{A}_X$ are Noetherian \cite[Theorem]{Fu69} 
and integral domains. Indeed the localization $(A_X)_{\mathbf m}$ of $A_X$ at $\mathbf mA_X$ is 
analytically irreducible and analytically unramified, and the natural morphism $\widehat{A}_X \rightarrow \widehat{(A_X)_{\mathbf m}}$ is injective. 
In addition, $\widehat{A}_X$ (resp. $A_X^\dagger$) 
is furnished a continuous integrable derivation $d : \widehat{A}_X \rightarrow \widehat{A}_X \otimes_{A_X}\Omega^1_{A_X/R}$ 
(resp. $d : A_X^\dagger \rightarrow A_X^\dagger \otimes_{A_X}\Omega^1_{A_X/R}$) 
which is an extension of the $R$-derivation on $A_X$. 
It is known that $A_X^\dagger \subset \widehat{A}_X$ are independent of the choice of the lift $A_X$ up to continuous $R$-isomorphisms 
by the approximation theorem for w.c.f.g. algebras in \cite[Theorem 2]{Bo81}. 
We do not discuss the functorialities here (see \cite[Section 2]{vP86} details). 

A continuous endomorphism $\varphi$ on $\widehat{A}_X$ (resp. $A^\dagger_X$) is said to be a $q$-Frobenius with respect to $\sigma$ 
if it satisfies $\varphi(a) \equiv a^q\, (\mathrm{mod}\, \mathbf m)$ and $\varphi|_R = \sigma$. 
The $q$-Frobenius alway exists on $\widehat{A}_X$ (resp. $A^\dagger_X$) by formal smoothness (resp. and \cite[Theorem 2]{Bo81}) 
\cite[Section 2.4]{vP86} and it is faithfully flat since $X$ is regular. 

Let $\overline{\mathcal X}$ be the Zariski closure of the immersion $\mathcal X \rightarrow \mathbb A^N_R \subset \mathbb P^N_R$ 
determined by $x_1, \cdots, x_N$ 
where $\mathbb P^N_R$ is the $N$-dimensional projective space over $\mathrm{Spec}\, R$, 
$\overline{X}$ the closed fiber of $\overline{\mathcal X}$ with the canonical open immersion $j_X : X \rightarrow \overline{X}$, 
$\widehat{\overline{\mathcal X}}$ the $p$-adic formal completion of $\overline{\mathcal X}$, and 
$\widehat{\overline{\mathcal X}}_K$ the Raynaud's generic fiber of $\widehat{\overline{\mathcal X}}$ 
with the specialization morphism $\mathrm{sp} : \widehat{\overline{\mathcal X}}_K \rightarrow \widehat{\overline{\mathcal X}}$ 
of locally ringed $G$-spaces. For a locally closed subset $Z \subset \widehat{\overline{\mathcal X}}$, we put 
$]Z[_{\widehat{\overline{\mathcal X}}} = \mathrm{sp}(Z)$ to be the tube of $Z$ in $\widehat{\overline{\mathcal X}}_K$. 

Since $\widehat{\overline{\mathcal X}}$ is smooth around $X$, there exists a strict neighborhood $W_0$ 
of $]X[_{\widehat{\overline{\mathcal X}}}$ in $]\overline{X}[_{\widehat{\overline{\mathcal X}}}$ such that $W_0$ is a smooth affinoid space over $\mathrm{Spm}\, K$. 
The following lemma is an essential part of \cite[Proposition 2.5.5]{Be96}. 

\begin{lemma}\label{frrig} Let $\varphi$ be a $q$-Frobenius $\varphi$ on $A_X^\dagger$. 
By replacing $W_0$ by a sufficiently small strict neighborhood of 
$]X[_{\widehat{\overline{\mathcal X}}}$ in $]\overline{X}[_{\widehat{\overline{\mathcal X}}}$, for any strict neighborhood $W$ 
of $]X[_{\widehat{\overline{\mathcal X}}}$ in $W_0$, 
there exists a strict neighborhood $W'$ 
of $]X[_{\widehat{\overline{\mathcal X}}}$ in $W$ such that the $q$-Frobenius $\varphi$ induces 
a morphism $\varphi : W' \rightarrow W$ of rigid analytic spaces over $\mathrm{Spm}\, K$. 
\end{lemma}

Instead of recalling the definition of overconvergent $F$-isocrystals \cite[Chapter 2]{Be96} (see \cite[10.2, 12.1]{CT03} in the general case 
without assuming the global embedding of $\overline{X}$ into a smooth formal scheme), 
we give an equivalent condition of the definition of overconvergent $F$-isocrystals in affine  smooth cases. 

\begin{theorem}\label{affov} \cite[Corollaire 2.5.8]{Be96} Let $\varphi$ be a $q$-Frobenius $\varphi$ on $A_X^\dagger$. 
The global section functor $\Gamma(]\overline{X}[_{\widehat{\overline{\mathcal X}}}, -)$ 
induces an equivalence between 
the category of overconvergent $F$-isocrystals on $X/K$ with respect to $\sigma$ 
and the category of $A^\dagger_{X, K}$-modules $M^\dagger$ of finite type which is furnished 
with an integrable connection $\nabla : M^\dagger \rightarrow M^\dagger\otimes_{A_X}\Omega^1_{A_X/R}$ 
and a horizontal isomorphism $\varphi_{M^\dagger} : (\varphi^\ast M^\dagger, \varphi^\ast\nabla) \rightarrow (M^\dagger, \nabla)$ called Frobenius. 
Note that $M^\dagger$ is a projective $A^\dagger_{X, K}$-module. 
\end{theorem}

\section{Change of base fields and Frobenius}\label{coeffrob} 

In this appendix we recall the push forward functor and the pull back functor between 
the categories of $F$-isocrystals under the base extensions and changing Frobenius by its powers. 
Suppose $k$ is an arbitrary perfect field of characteristic $p$. 
Let $L$ be a finite extension of $K$ with the residue field $l$ 
such that $L$ admits an extension $\sigma$ of the $q$-Frobenius $\sigma$ on $K$ (note that we use the same notation $\sigma$). 
Let $X$ be a smooth scheme separated of finite type over $\mathrm{Spec}\, k$ with a completion $\overline{X}$ of $X$ over $\mathrm{Spec}\, k$, 
and put the scalar extension $X_l = X \times_{\mathrm{Spec}\, k}\mathrm{Spec}\, l$ and the projection $\pi_{l/k} : X_l \rightarrow X$. 
For a positive integer $n$, we define the pull back functor
$$
    \pi_{(L, \sigma^n)/(K, \sigma)}^\ast : \left(\begin{array}{l} \mbox{\rm overconvergent $F$-isocrystals} \\ 
    \mbox{\rm on $X/K$ with respect to $\sigma$}\end{array}\right) 
    \rightarrow \left(\begin{array}{l} \mbox{\rm overconvergent $F$-isocrystals} \\ 
    \mbox{\rm on $X_l/L$ with respect to $\sigma^n$}\end{array}\right) 
$$
and the push forward functor 
$$
    \pi_{(L, \sigma^n)/(K, \sigma)\, \ast} : \left(\begin{array}{l} \mbox{\rm overconvergent $F$-isocrystals} \\ 
    \mbox{\rm on $X_l/L$ with respect to $\sigma^n$}\end{array}\right) 
    \rightarrow \left(\begin{array}{l} \mbox{\rm overconvergent $F$-isocrystals} \\ 
    \mbox{\rm on $X/K$ with respect to $\sigma$}\end{array}\right) 
$$
as follows. Locally on $X$, there exists a smooth lift $\mathcal X = \mathrm{Spec}\, A_X$ over $\mathrm{Spec}\, R$ 
and a $q$-Frobenius $\varphi$ on the $p$-adically  weak completion $A_X^\dagger$ of $A_X$ which is compatible to $\sigma$ (Appendix \ref{FrobA}). 
Then the functors are defined by 
$$
\pi_{(L, \sigma^n)/(K, \sigma)}^\ast(\mathcal M^\dagger, \nabla_{\mathcal M^\dagger}, F_{\mathcal M^\dagger}) 
= (L\otimes_K \mathcal M^\dagger, \mathrm{id}_L\otimes \nabla_{\mathcal M^\dagger}, \mathrm{id}_L\otimes_{\sigma^n} (\sigma^n \otimes F_{\mathcal M^\dagger}^n))
$$
where $\otimes_{\sigma^n}$ means the scalar extension by $\sigma^n : L \rightarrow L$, 
and 
$$
  \begin{array}{l}
    \pi_{(L, \sigma^n)/(K, \sigma)\, \ast}(\mathcal N^\dagger, \nabla_{\mathcal M^\dagger}, F_{\mathcal N^\dagger}) 
     = (\oplus_{i=0}^{n-1}(\varphi^i)^\ast(\pi_\ast \mathcal N^\dagger, \pi_\ast \nabla), F_{\pi_{(L, \sigma^n)/(K, \sigma)\, \ast}(\mathcal N^\dagger)}) \\
     F_{\pi_{(L, \sigma^n)/(K, \sigma)\, \ast}(\mathcal N^\dagger)}(a_0 \otimes_{\sigma} m_0, \cdots, a_{n-1} \otimes_{\sigma} m_{n-1}) \\
     \hspace*{35mm} = (F_{\mathcal N^\dagger}(a_{n-1}\otimes_\sigma m_{n-1}), a_0 \otimes_{\sigma}m_0, \cdots, a_{n-2} \otimes_{\sigma}m_{n-2}). 
     \end{array}
$$
Then $\pi_{(L, \sigma^n)/(K, \sigma)}^\ast$ is a left adjoint of $\pi_{(L, \sigma^n)/(K, \sigma)\, \ast}$. 
For an overconvergent $F$-isocrystal $\mathcal M^\dagger$ on $X/K$, 
the adjoint $\mathrm{ad}_{(L, \sigma^n)/(K, \sigma)} : \mathcal M^\dagger 
\rightarrow \pi_{(L, \sigma^n)/(K, \sigma)\, \ast}\pi_{(L, \sigma^n)/(K, \sigma), n}^\ast\mathcal M^\dagger$ 
of $\mathrm{id} : \pi_{(L, \sigma^n)/(K, \sigma)}^\ast\mathcal M^\dagger \rightarrow \pi_{(L, \sigma^n)/(K, \sigma)}^\ast\mathcal M^\dagger$ is given by 
$$
      \mathrm{ad}_{(L, \sigma^n)/(K, \sigma)}(m) 
      = (m, F_{\mathcal M^\dagger}^{-1}(m), \cdots, F_{\mathcal M^\dagger}^{-n+1}(m))
$$
where $F_{\mathcal M^\dagger}^{-i} : \mathcal M^\dagger \rightarrow (\varphi^i)^\ast\mathcal M^\dagger$ is the 
inverse of $F_{\mathcal M^\dagger}^i 
= F_{\mathcal M^\dagger}\circ \cdots \circ (\varphi^{i-1})^\ast(F_{\mathcal M^\dagger}) : (\varphi^i)^\ast\mathcal M^\dagger \rightarrow \mathcal M^\dagger$. 
We also define a trace map 
$\mathrm{Tr}_{(L, \sigma^n)/(K, \sigma)} : \pi_{(L, \sigma^n)/(K, \sigma)\, \ast}\pi_{(L, \sigma^n)/(K, \sigma)}^\ast\mathcal M^\dagger \rightarrow \mathcal M^\dagger$ by 
$$
      \mathrm{Tr}_{(L, \sigma^n)/(K, \sigma)}(a_0 \otimes m_0, a_1 \otimes m_1, \cdots, a_{n-1} \otimes m_{n-1}) 
      = \sum_{i=0}^{n-1} \mathrm{tr}_{L/K}(\sigma^i(a_i))F_{\mathcal M^\dagger}^i(1 \otimes m_i), 
$$
where $\mathrm{tr}_{L/K} : L \rightarrow K$ is the trace of finite extension of fields. 
Since $\mathrm{tr}_{L/K}$ commutes with the $q$-Frobenius $\sigma$, 
the trace map $\mathrm{Tr}_{(L, \sigma^n)/(K, \sigma)}$ is well-defined and 
commutes with connections and Frobenius. Hence it is a morphism 
of overconvergent $F$-isocrystals on $X/K$ with respect to $\sigma$. 

\begin{lemma}\label{adtr} 
\begin{enumerate}
\item As morphisms of (over)convergent $F$-isocrystals on $X/K$ with respect to $\sigma$, the following identity holds:
$$
\mathrm{Tr}_{(L, \sigma^n)/(K, \sigma)}\circ\mathrm{ad}_{(L, \sigma^n)/(K, \sigma)} = n\mathrm{deg}(L/K)\mathrm{id}.
$$
\item Let $K \subset L \subset M$ be a sequence of finite extensions with the compatible $q$-Frobenius $\sigma$. 
Then the associativities below hold:
$$
   \begin{array}{l}
   \pi_{(M, \sigma^{mn})/(K, \sigma)}^\ast = \pi_{(M, \sigma^{mn})/(L, \sigma^m)}^\ast\circ \pi_{(L, \sigma^m)/(K, \sigma)}^\ast, \\
\pi_{(M, \sigma^{mn})/(K, \sigma)\, \ast} = \pi_{(L, \sigma^m)/(K, \sigma)\, \ast}\circ \pi_{(M, \sigma^{mn})/(L, \sigma^m)\, \ast}, \\
\mathrm{ad}_{(M, \sigma^{mn})/(K, \sigma)} = \pi_{(L, \sigma^m)/(K, \sigma)\, \ast}(\mathrm{ad}_{(M, \sigma^{mn})/(L, \sigma^m)})
\circ\mathrm{ad}_{(L, \sigma^m)/(K, \sigma)} \\
\mathrm{Tr}_{(M, \sigma^{mn})/(K, \sigma)} = \mathrm{Tr}_{(L, \sigma^m)/(K, \sigma)}\circ 
\pi_{(L, \sigma^m)/(K, \sigma)\, \ast}(\mathrm{Tr}_{(M, \sigma^{mn})/(L, \sigma^m)}). 
\end{array}
$$ 
\end{enumerate} 
\end{lemma}

\begin{remark} In this appendix we define the functors $\pi_{(L, \sigma^n)/(K, \sigma)}^\ast$ and $\pi_{(L, \sigma^n)/(K, \sigma)\, \ast}$ 
only for smooth schemes. Using the method of hypercoverings \cite[2.3.2 (iii), 2.3.7]{Be96} \cite[10.2, 12.1]{CT03}, one can define 
the same for arbitrary schemes locally of finite type over $\mathrm{Spec}\, k$. 
\end{remark}


\begin{thebibliography}{}

\bibitem{Ab18b} T.Abe, 
Langlands correspondence for isocrystals and the existence of crystalline companions for curves. 
J. Amer. Math. Soc. \textbf{31}, no.4, 921-1057 (2018). 

\bibitem{Ab18} T.Abe, 
Langlands program for $p$-adic coefficients and the petits camarades conjecture. 
J.Reine Angew. Math. \textbf{734}, 59-69 (2018). 


\bibitem{AC18} T.Abe; D.Caro, Theory of weights in $p$-adic cohomology. 
Amer. J. Math. \textbf{140}, no.4, 879-975 (2018).

\bibitem{AE16} T.Abe; H.Esnault, 
A Lefschetz theorem for overconvergent isocrystals with Frobenius structure. 
Ann. Ec. Nor. Sup. 4e s\'erie, \textbf{52}, 1243 -1264 (2019). 

\bibitem{AD18} E.Ambrosi; M.D'Addezio, 
Maximal tori of monodromy groups of $F$-isocrystals and an application to abelian varieties. 
arXiv:1811.08423. 

\bibitem{Be96} P.Berthelot, 
Cohomologie rigide et cohomologie rigide \`a support propre, part 1. 
Pr\'epublication IRMAR 96-03, available at \url{https://perso.univ-rennes1.fr/pierre.berthelot/}.

\bibitem{Bo81} S.Bosch, 
A rigid analytic version of M.Artin's theorem on analytic equations. 
Math. Ann. \textbf{255}, 395-404 (1981).


\bibitem{CT03}  B.Chiarellotto; N.Tsuzuki, 
Cohomological descent of rigid cohomology for etale coverings. 
Rendiconti di Padova \textbf{109}, 63-215 (2003).

\bibitem{CT09} B.Chiarellotto; N.Tsuzuki, 
Logarithmic growth and Frobenius filtrations for solutions of $p$-adic differential equations. 
J. Inst. Math. Jussieu \textbf{8}, no. 3, 465-505 (2009). 

\bibitem{CT11} B.Chiarellotto; N.Tsuzuki, 
Log-growth filtration and Frobenius slope filtration of $F$-isocrystals at the generic and special points. 
Doc. Math. \textbf{16}, 33-69 (2011). 

\bibitem{Ch83} G.Christol, 
\textit{Modules diff\'erentiels et \'equations diff\'erentiels $p$-adiques.} 
Queen's Papers 
In Pure and Applied Mathematics, \textbf{66}, Queen's University, Kingston, 1983.

\bibitem{Ch84} G.Christol, 
Un  th\'eor\`eme de transfert pour les disques singuli\`eres reguli\'eres, Cohomologie $p$-adique. 
Ast\'erisque, \textbf{119-120}, SMF, 151-168, 1984.

\bibitem{Cr87}  R.Crew, 
$F$-isocrystals and $p$-adic representations. Algebraic geometry, Bowdoin, 
1985 (Brunswick, Maine, 1985), 111-138, Proc. Sympos. Pure Math., \textbf{46}, Part 2, Amer. Math. Soc., 
Providence, RI, 1987.

\bibitem{Cr92} R.Crew, 
$F$-isocrystals and their monodromy groups. Ann. Scient. \'Ec. Norm. Sup. \textbf{25}, 429-464 (1992).

\bibitem{dJ98} J.A.De Jong, 
Homomorphisms of Barsotti-Tate groups and crystals in positive characteristic. 
Invent. Math. \textbf{134}, no. 2, 301-333 (1998).

\bibitem{DA20} M.D'Addezio, 
Parabolicity conjecture of $F$-isocrystals. arXiv:2012.12879. 

\bibitem{De81} P.Deligne, 
La conjecture de Weil, II. Publ. Math. IH\'ES \textbf{52}, 313-428 (1981).

\bibitem{Dw73} B.Dwork, 
Normalized period matrices. II. Ann. of Math. (2) \textbf{98}, 1-57 (1973). 

\bibitem{Dw74} B.Dwork, 
Bessel functions as $p$-adic functions of the argument. 
Duke Math. J. \textbf{41}, no.4, 711-738 (1974).

\bibitem{El73} R.Elkik, 
Solutions d'\'equations \`a coefficients dans un anneau hens\'elien. Ann. Sci. \'Ecole Norm.
Sup. \textbf{6}, 553-603 (1973).

\bibitem{EL93} J.-Y.\'Etesse; B. Le Stum,
Fonctions $L$ associ\'ees aux $F$-isocristaux surconvergents, I: Interpr\'etation
cohomologique. Math. Ann. \textbf{296}, 557-576 (1993).

\bibitem{Fu69} W.Fulton, 
A note on weakly complete algebras. Bull. Amer. Math. Soc. \textbf{75},591-593 (1969). 

\bibitem{Gr} A.Grothendieck, 
\textit{Rev\^{e}tements \'etales et groupe fondamental (SGA1).} Lecture Notes in Math. \textbf{224}, Springer-Velag (1971). 

\bibitem{Ha17} C.D.Haessig, 
$L$-functions of symmetric powers of Kloosterman sums (unit root $L$-functions and $p$-adic estimates). 
Math. Ann. \textbf{369}, no. 1-2, 17-47 (2017).

\bibitem{Ka72} N.M.Katz, 
$p$-adic properties of modular schemes and modular forms. 
Modular Functions of One Variable, III (Proc. Internat. Summer School, Univ. Antwerp, Antwerp, 
1972). Lecture Notes in Math. \textbf{350},  Springer, Berlin,  69-190, 1973.

\bibitem{Ka99} N.M.Katz, 
Space filling curves over finite fields. 
Math. Res. Lett. \textbf{6}, 613-624 (1999). 

\bibitem{Ke01} K.S.Kedlaya, 
Power series and $p$-adic algebraic closures.  
J. Number theory \textbf{89}, 324-339 (2001). 

\bibitem{Ke05} K.S.Kedlaya, 
Slope filtrations revisited. Doc. Math. \textbf{10}, 447-525 (2005); 
errata, ibid. \textbf{12}, 361-362 (2007). 

\bibitem{Ke06} K.S.Kedlaya, 
Finiteness of rigid cohomology with coefficients. Duke Math. J. \textbf{134}, 15-97 (2006).

\bibitem{Ke06b} K.S.Kedlaya, 
Fourier transforms and $p$-adic ``Weil II". Compos. Math. \textbf{142}, 1426-1450 (2006).


\bibitem{Ke07} K.S.Kedlaya, 
Semistable reduction for overconvergent $F$-isocrystals, I: Unipotence and logarithmic extensions. 
Compos. Math. \textbf{143}, 1164-1212 (2007).

\bibitem{Ke08} K.S.Kedlaya, 
Slope filtration for relative Frobenius. 
Ast\'erisque \textbf{319}, 259-301 (2008).

\bibitem{Ke08b} K.S.Kedlaya, 
Semistable reduction for overconvergent $F$-isocrystals, II: A valuation-theoretic approach. 
Compos. Math. \textbf{144}, 657-672 (2008).

\bibitem{Ke16} K.S.Kedlaya,
Notes on isocrystals. 
arXiv:1606.01321.

\bibitem{Ke19} K.S.Kedlaya,
\'Etale and crystalline companions, I.  
arXiv:1811.00204v2. 	

\bibitem{Kr16} J.Kramer-Miller, 
The monodromy of $F$-isocrystals with log-decay. 
arXiv:1612.01164. 

\bibitem{Kr19} J.Kramer-Miller, 
Slope filtrations of F-isocrystals and logarithmic decay. 
Math. Res. Lett. \textbf{28}, no.1, 107-125 (2021). 

\bibitem{Kr18} J.Kramer-Miller, 
The monodromy of unit-root $F$-isocrystals with geometric origin. 
to appear in Compos. Math. 

\bibitem{MW68} P.Monsky; G.Washnitzer, 
Formal Cohomology I.   Ann. of  Math.  \textbf{88}, 181-217 (1968). 

\bibitem{Oh18} S.Ohkubo,
Logarithmic growth filtrations for $(\varphi, \nabla)$-modules over the bounded Robba ring.  
 Compos. Math. \textbf{157}, no. 6, 1265-1301 (2021).

\bibitem{Ro75} P.Robba, 
On the index of $p$-adic differential operators. I. Annals Math. \textbf{101} ,
280-316 (1975).

\bibitem{IS} Shioda, T., Inose, H, On singular $K3$ surfaces. 
Complex analysis and algebraic geometry (W.L. Bailey and T.Shioda eds.), 119-136, Cambridge 
University Press 1977. 

\bibitem{SB85} J.Stienstra; F.Beukers, 
On the Picard-Fuchs equation and the formal Brauer group of certain elliptic K3-surfaces. 
Math. Ann. \textbf{271}, no. 2, 269-304 (1985).

\bibitem{Ts96} N.Tsuzuki, 
The overconvergence of morphisms of \'etale $\varphi$-$\nabla$-spaces on a local field. 
Compos. Math. \textbf{103}, no. 2, 227-239 (1996).

\bibitem{Ts98} N.Tsuzuki, 
Finite local monodromy of overconvergent unit-root $F$-isocrystals on a curve. 
Amer. J. Math. \textbf{120}, no. 6, 1165-1190 (1998).

\bibitem{Ts98b} N.Tsuzuki, 
Slope filtration of quasi-unipotent overconvergent $F$-isocrystals. 
Ann. Inst. Fourier (Grenoble) \textbf{48}, no. 2, 379-412 (1998).

\bibitem{Ts02} N.Tsuzuki, 
Morphisms of $F$-isocrystals and the finite monodromy theorem for unit-root $F$-isocrystals. 
Duke Math. J. \textbf{111}, no. 3, 385-418 (2002). 

\bibitem{Ts19} N.Tsuzuki, 
Constancy of Newton polygons of $F$-isocrystals on Abelian varieties and isotriviality of families of curves. 
 J. Inst. Math. Jussieu \textbf{20}, no. 2, 587-625  (2021). 
\bibitem{vP86} M.Van der Put, 
The cohomology of Monsky and Washnitzer. M\'emoires de la S.M.F. 2e s\'erie, \textbf{23}, 33-59 (1986). 

\bibitem{Wa99} D.Wan, 
Dwork's conjecture of unit root zeta functions. Annals of Math. \textbf{150}, 867-927 (1999). 

\bibitem{Wa00a} D.Wan, 
Higher rank case of Dwork's conjecture. J. Amer. Math. Soc. \textbf{13}, 807-852 (2000). 

\bibitem{Wa00b} D.Wan, 
Rank one case of Dwork's conjecture. J. Amer. Math. Soc. \textbf{13}, 853-908 (2000). 
\end{thebibliography}
\end{document}